\documentclass{article}

\usepackage[leqno,intlimits]{amsmath}
\usepackage{amssymb}
\usepackage{amsthm}
\usepackage{hyperref}
\usepackage{url}
\usepackage[hyphenbreaks]{breakurl} % For breaking long urls in the references.

\usepackage[margin=1in]{geometry}

\allowdisplaybreaks[1]

\newcommand{\E}{\mathbb{E}} % Expectation E
\newcommand{\p}{\mathbb{P}} % Probability P
\newcommand{\eps}{\varepsilon} % Epsilon
\DeclareMathOperator{\Var}{Var} % Variance Var
\DeclareMathOperator{\Dist}{{\bf D}}

\DeclareMathOperator{\NONMANIP}{NONMANIP}
\DeclareMathOperator{\Maj}{Maj}
\DeclareMathOperator{\maj}{maj}
\DeclareMathOperator{\Inf}{Inf}
\DeclareMathOperator{\tp}{top}

\DeclareMathOperator{\Lg}{Lg}
\DeclareMathOperator{\Sm}{Sm}
\DeclareMathOperator{\LD}{LD}
\newcommand{\adj}[2]{\ensuremath{[#1:#2]}}

\newtheorem{theorem}{Theorem}[section]
\newtheorem{lemma}[theorem]{Lemma}

\newtheorem{corollary}[theorem]{Corollary}

\newtheorem{definition}{Definition}  % Note, this italicizes everything

\begin{document} %%%
%%%%%%%%%%%%%%%%%%%%

%%%%%%%%%%%%%
%%% Title %%%
%%%%%%%%%%%%%

\title{A quantitative Gibbard-Satterthwaite theorem without neutrality}
\author{
	Elchanan Mossel
	\thanks{University of California, Berkeley and Weizmann Institute of Science; \texttt{mossel@stat.berkeley.edu}; supported by NSF CAREER award (DMS 0548249) and by DOD ONR grant N000141110140.}
	\and
	Mikl\'os Z. R\'acz
	\thanks{University of California, Berkeley; \texttt{racz@stat.berkeley.edu}; supported by a UC Berkeley Graduate Fellowship and by NSF DMS 0548249 (CAREER).}
}
\maketitle

%%%%%%%%%%%%%%%%
%%% Abstract %%%
%%%%%%%%%%%%%%%%

\begin{abstract}
Recently, quantitative versions of the Gibbard-Satterthwaite theorem were proven for $k=3$ alternatives by Friedgut, Kalai, Keller and Nisan and for neutral functions on $k \geq 4$ alternatives by Isaksson, Kindler and Mossel.

We prove a quantitative version of the Gibbard-Satterthwaite theorem for general social choice functions for any number $k \geq 3$ of alternatives. In particular we show that for a social choice function $f$ on $k \geq 3$ alternatives and $n$ voters, which is $\eps$-far from the family of nonmanipulable functions, a uniformly chosen voter profile is manipulable with probability at least inverse polynomial in $n$, $k$, and $\eps^{-1}$.

Removing the neutrality assumption of previous theorems is important for multiple reasons. For one, it is known that there is a conflict between anonymity and neutrality, and since most common voting rules are anonymous, they cannot always be neutral. Second, virtual elections are used in many applications in artificial intelligence, where there are often restrictions on the outcome of the election, and so neutrality is not a natural assumption in these situations.

Ours is a unified proof which in particular covers all previous cases established before. The proof crucially uses reverse hypercontractivity in addition to several ideas from the two previous proofs. Much of the work is devoted to understanding functions of a single voter, and in particular we also prove a quantitative Gibbard-Satterthwaite theorem for one voter.
\end{abstract}

%%%%%%%%%%%%%%%%
%%% Document %%%
%%%%%%%%%%%%%%%%

%%%%%%%%%%%%%%%%%%%%%%%%%%
\section{Introduction} %%%
%%%%%%%%%%%%%%%%%%%%%%%%%%

One of the main goals in social choice theory is to come up with ``good'' voting systems, which satisfy a few natural requirements. This problem is increasingly relevant in the area of artificial intelligence and computer science as well, where virtual elections are now an established tool in preference aggregation (see the survey by Faliszewski and Procaccia~\cite{faliszewski2010ai}). Many of the results in the study of social choice are negative: it is impossible to design a voting system that satisfies a few desired properties all at once. The first realization of an apparent problem is due to Condorcet, who, at the end of the $18^{\text{th}}$ century, noticed the following paradox: when ranking three candidates, $a$, $b$, and $c$, it may happen that a majority of voters prefer $a$ over $b$, a majority prefers $b$ over $c$, and a majority prefers $c$ over $a$, thus producing an ``irrational'' circular ranking of the candidates. Arrow's impossibility theorem~\cite{Arrow:50,Arrow:63} showed that this paradox holds under very natural assumptions, thus marking the basis of modern social choice theory.

A naturally desirable property of a voting system is \emph{strategyproofness} (a.k.a.\ nonmanipulability): no voter should benefit from voting strategically, i.e.\ voting not according to her true preferences. However, Gibbard~\cite{gibbard1973manipulation} and Satterthwaite~\cite{satterthwaite1975strategy} showed that no reasonable voting system can be strategyproof. Before stating their result, let us specify the problem more formally. 

We consider $n$ voters electing a winner among $k$ alternatives. The voters specify their opinion by ranking the alternatives, and the winner is determined according to some predefined \emph{social choice function} (SCF) $f : S_k^n \to \left[k\right]$ of all the voters' rankings, where $S_k$ denotes the set of all possible total orderings of the $k$ alternatives. We call a collection of rankings by the voters a \emph{ranking profile}. We say that a SCF is \emph{manipulable} if there exists a ranking profile where a voter can achieve a more desirable outcome of the election according to her true preferences by voting in a way that does not reflect her true preferences (see Definition \ref{def:manip} for a more detailed definition).

The Gibbard-Satterthwaite theorem states that any SCF which is not a dictatorship (i.e.\ not a function of a single voter), and which allows at least three alternatives to be elected, is manipulable. This has contributed to the realization that it is unlikely to expect truthfulness in voting. Consequently, there have been many branches of research devoted to understanding the extent of manipulability of voting systems, and to finding ways of circumventing the negative results. 

One approach, introduced by Bartholdi, Tovey and Trick~\cite{bartholdi1989computational}, suggests computational complexity as a barrier against manipulation: if it is computationally hard for a voter to manipulate, then she would just tell the truth (we refer to the survey by Faliszewski and Procaccia~\cite{faliszewski2010ai} for a detailed history of the surrounding literature). This is a worst-case approach, and while worst-case hardness of manipulation is a desirable property for a SCF to have, this does not tell us anything about \emph{typical} instances of the problem---is it easy or hard to manipulate \emph{on average}?

A recent line of research with an average-case algorithmic approach has suggested that manipulation is indeed easy on average; see e.g.\ Kelly~\cite{kelly1993almost}, Conitzer and Sandholm~\cite{conitzer2006nonexistence}, and Procaccia and Rosenschein~\cite{procaccia2007junta} for results on certain restricted classes of SCFs (see also the survey~\cite{faliszewski2010ai}).

A different approach was taken by Friedgut, Kalai, Keller and Nisan~\cite{friedgut2008elections,friedgut2011quantitative}, who looked at the fraction of ranking profiles that are manipulable. To put it differently: assuming each voter votes independently and uniformly at random (known as the \emph{impartial culture assumption} in the social choice literature), what is the probability that a ranking profile is manipulable? Is it perhaps exponentially small (in the parameters $n$, $k$), or is it nonnegligible? Of course, if the SCF is nonmanipulable then this probability is zero. Similarly, if the SCF is ``close'' to being nonmanipulable in some sense, then this probability can be small. We say that a SCF $f$ is $\eps$-far from the family of nonmanipulable functions, if one must change the outcome of $f$ on at least an $\eps$-fraction of the ranking profiles in order to transform $f$ into a nonmanipulable function. Friedgut et al.\ conjectured that if $k \geq 3$ and the SCF $f$ is $\eps$-far from the family of nonmanipulable functions, then the probability of a ranking profile being manipulable is bounded from below by a polynomial in $1/n$, $1/k$, and $\eps$. Moreover, they conjectured that a random manipulation will succeed with nonnegligible probability, suggesting that manipulation by computational agents in this setting is easy.

Friedgut et al.\ proved their conjecture in the case of $k=3$ alternatives, showing a lower bound of $C \eps^{6}/n$ in the general setting, and $C' \eps^{2}/n$ in the case when the SCF is neutral (invariant under changes made to the names of the alternatives), where $C,C'$ are constants. Note that this result does not have any computational consequences, since when there are only $k=3$ alternatives, a computational agent may easily try all possible permutations of the alternatives to find a manipulation (if one exists). Several follow-up works have since extended this result. First, Xia and Conitzer~\cite{xia2008sufficient} used the proof technique of Friedgut et al.\ to extend their result to a constant number of alternatives, assuming several additional technical assumptions. However, this still does not have any computational consequences, since the result holds only for a constant number of alternatives. Dobzinski and Procaccia~\cite{dobzinski2008frequent} proved the conjecture in the case of two voters under the assumption that the SCF is Pareto optimal. Finally, the latest work is due to Isaksson, Kindler and Mossel~\cite{isaksson2010geometry}, who proved the conjecture in the case of $k \geq 4$ alternatives with only the added assumption of \emph{neutrality}. Moreover, they showed that a random manipulation which replaces four adjacent alternatives in the preference order of the manipulating voter by a random permutation of them succeeds with nonnegligible probability. Since this result is valid for any number of ($k \geq 4$) alternatives, it does have computational consequences, implying that for neutral SCFs, manipulation by computational agents is easy on average.

In this paper we remove the neutrality condition and resolve the conjecture of Friedgut et al.: if $k \geq 3$ and the SCF $f$ is $\eps$-far from the family of nonmanipulable functions, then the probability of a ranking profile being manipulable is bounded from below by a polynomial in $1/n$, $1/k$, and $\eps$. We continue by first presenting our results, then discussing their implications, and finally we conclude this section by commenting on the techniques used in the proof.

%%%%%%%%%%%%%%%%%%%%%%%%%%%%
\subsection{Basic setup} %%%
%%%%%%%%%%%%%%%%%%%%%%%%%%%%

% 
% The voters specify their opinion by ranking the alternatives, and the winner is determined according to some pre-defined \emph{social choice function} (SCF) $f : S_k^n \to \left[k\right]$ of all the voters' rankings, where $S_k$ denotes the set of all possible permutations (total orderings) of the $k$ alternatives. We call a collection of rankings by the voters a \emph{ranking profile}.%We denote by $\sigma = \left( \sigma_1, \dots, \sigma_n \right)$ the collection of rankings by the voters and call this a \emph{ranking profile}.

Recall that our basic setup consists of $n$ voters electing a winner among $k$ alternatives via a SCF $f : S_k^n \to \left[k\right]$. We now define manipulability in more detail:
\begin{definition}[Manipulation points]\label{def:manip}
Let $\sigma \in S_k^n$ be a ranking profile. Write $a \stackrel{\sigma_i}{>} b$ to denote that alternative $a$ is preferred over $b$ by voter $i$. A SCF $f \colon S_k^n \to [k]$ is \emph{manipulable} at the ranking profile $\sigma \in S_k^n$ if there exists a $\sigma' \in S_k^n$ and an $i\in \left[n\right]$ such that $\sigma$ and $\sigma'$ only differ in the $i^{\text{th}}$ coordinate and
\begin{equation*}
    f(\sigma') \stackrel{\sigma_i}{>} f(\sigma).
\end{equation*}
In this case we also say that $\sigma$ is a {\em manipulation point} of $f$, and that $(\sigma,\sigma')$ is {\em a manipulation pair} for $f$. We say that $f$ is {\em manipulable} if it is manipulable at some point $\sigma$. We also say that $\sigma$ is an $r$-{\em manipulation point} of $f$ if $f$ has a manipulation pair $(\sigma,\sigma')$ such that $\sigma'$ is obtained from $\sigma$ by permuting (at most) $r$ adjacent alternatives in one of the coordinates of $\sigma$. (We allow $r > k$---any manipulation point is an $r$-manipulation point for $r > k$.)

Let $M \left( f \right)$ denote the set of manipulation points of the SCF $f$, and for a given $r$, let $M_r \left( f \right)$ denote the set of $r$-manipulation points of $f$. When the SCF is obvious from the context, we write simply $M$ and $M_r$.
\end{definition}
Gibbard and Satterthwaite proved the following theorem.
\begin{theorem}[Gibbard-Satterthwaite~\cite{gibbard1973manipulation,satterthwaite1975strategy}]\label{thm:GS}
Any SCF $f \colon S_k^n \to \left[k\right]$ which takes at least three values and is not a dictator (i.e.\ not a function of only one voter) is manipulable.
\end{theorem}
This theorem is tight in the sense that \emph{monotone} SCFs which are dictators or only have two possible outcomes are indeed nonmanipulable (a function is non-monotone, and clearly manipulable, if for some ranking profile a voter can change the outcome from, say, $a$ to $b$ by moving $a$ ahead of $b$ in her preference). It is useful to introduce a refined notion of a dictator before defining the set of nonmanipulable SCFs.
\begin{definition}[Dictator on a subset]\label{def:dict_subset}
For a subset of alternatives $H \subseteq \left[k\right]$, let $\tp_H$ be the SCF on one voter whose output is always the top ranked alternative among those in $H$.
\end{definition}
\begin{definition}[Nonmanipulable SCFs]
We denote by $\NONMANIP \equiv \NONMANIP \left( n, k \right)$ the set of nonmanipulable SCFs, which is the following:
\begin{align*}
\NONMANIP \left( n, k\right) &= \left\{ f : S_k^n \to \left[k\right] \mid f \left( \sigma \right) = \tp_H \left( \sigma_i \right) \text{ for some } i \in \left[n\right], H \subseteq \left[k\right], H \neq \emptyset \right\} \bigcup\\
&\bigcup \left\{ f : S_k^n \to \left[k\right] \mid f \text{ is a monotone function taking on exactly two values} \right\}.
\end{align*}
When the parameters $n$ and $k$ are obvious from the context, we omit them.
\end{definition}
Another important class of functions, which is larger than $\NONMANIP$, but which has a simpler description, is the following.
\begin{definition}
Define, for parameters $n$ and $k$ that remain implicit (when used the parameters will be obvious from the context):
\[
  \overline{\NONMANIP} = \{f \colon S_k^n \to [k] \mid f \text{ only depends on one coordinate or takes at most two values}\}.
\]
\end{definition}
The notation should be thought of as ``closure'' rather than ``complement''. We remark that in \cite{isaksson2010geometry} the set $\overline{\NONMANIP}$ is denoted by $\NONMANIP$---but these two sets of functions should not be confused.

As discussed previously, our goal is to study manipulability from a quantitative viewpoint, and in order to do so we need to define the distance between SCFs.
\begin{definition}[Distance between SCFs]
The distance $\Dist(f,g)$ between two SCFs $f,g \colon S_k^n \to \left[ k \right]$ is defined as the fraction of inputs on which they differ: $\Dist\left(f,g\right) = \p\left(f\left(\sigma\right) \neq g\left(\sigma\right)\right)$, where $\sigma \in S_k^n$ is uniformly selected. For a class $G$ of SCFs, we write $\Dist\left(f,G\right) = \min_{g \in G} \Dist\left(f,g\right)$.
\end{definition}

% The quantitative Gibbard-Satterthwaite theorems of Friedgut, Kalai, Keller and Nisan~\cite{friedgut2008elections,friedgut2011quantitative}, and Isaksson, Kindler and Mossel \cite{isaksson2010geometry} involve the distance of a SCF from $\overline{\NONMANIP}$. Any SCF that is not in $\overline{\NONMANIP}$ is manipulable (by the Gibbard-Satterthwaite theorem), but as some SCFs in $\overline{\NONMANIP}$ are manipulable as well, ideally a quantitative Gibbard-Satterthwaite theorem would involve the distance of a SCF from the set of (truly) nonmanipulable SCFs, $\NONMANIP$.

The concepts of anonymity and neutrality of SCFs will be important to us, so we define them here.
\begin{definition}[Anonymity]\label{def:anonymity}
 A SCF is \emph{anonymous} if it is invariant under changes made to the names of the voters. More precisely, a SCF $f : S_k^n \to \left[k\right]$ is anonymous if for every $\sigma = \left( \sigma_1, \dots, \sigma_n \right) \in S_k^n$ and every $\pi \in S_n$,
\[
 f\left( \sigma_1, \dots, \sigma_n \right) = f\left( \sigma_{\pi\left( 1 \right)}, \dots, \sigma_{\pi \left( n \right)} \right).
\]
\end{definition}
\begin{definition}[Neutrality]\label{def:neutrality}
 A SCF is \emph{neutral} if it is invariant under changes made to the names of the alternatives. More precisely, a SCF $f : S_k^n \to \left[k\right]$ is neutral if for every $\sigma = \left( \sigma_1, \dots, \sigma_n \right) \in S_k^n$ and every $\pi \in S_k$,
\[
 f \left( \pi \circ \sigma_1, \dots, \pi \circ \sigma_n \right) = \pi \left( f \left( \sigma \right) \right).
\]
\end{definition}

%%%%%%%%%%%%%%%%%%%%%%%%%%%%%%%%%%%%%%%%%%%%%%%%%%%
\subsection{Our main result}\label{sec:results} %%%
Our main result, which resolves the conjecture of Friedgut et al.~\cite{friedgut2008elections,friedgut2011quantitative}, is the following.
\begin{theorem}\label{cor:k_refined_truenonmanip}
Suppose we have $n \geq 1$ voters, $k \geq 3$ alternatives, and a SCF $f : S_k^n \to \left[k\right]$ satisfying $\Dist \left( f, \NONMANIP \right) \geq \eps$. Then
\begin{equation}\label{eq:manip_refined_true}
\p\left( \sigma \in M \left( f \right) \right)\\
\geq \p \left( \sigma \in M_4 \left( f \right) \right)\geq p \left( \eps, \frac{1}{n}, \frac{1}{k} \right)
\end{equation}
for some polynomial $p$, where $\sigma \in S_k^n$ is selected uniformly. In particular, we show a lower bound of $\frac{\eps^{15}}{10^{39} n^{67} k^{166}}$.

An immediate consequence is that
\[
\p \left( \left( \sigma, \sigma' \right) \text{ is a manipulation pair for } f \right) \geq q \left( \eps, \frac{1}{n}, \frac{1}{k} \right)
\]
for some polynomial $q$, where $\sigma \in S_k^n$ is uniformly selected, and $\sigma'$ is obtained from $\sigma$ by uniformly selecting a coordinate $i \in \left\{1, \dots, n \right\}$, uniformly selecting $j \in \left\{1, \dots, n-3 \right\}$, and then uniformly randomly permuting the following four adjacent alternatives in $\sigma_i$: $\sigma_i \left( j \right), \sigma_i \left( j + 1 \right), \sigma_i \left( j + 2 \right)$, and $\sigma_i \left( j + 3 \right)$. In particular, the specific lower bound for $\p \left( \sigma \in M_4 \left( f \right) \right)$ implies that we can take $q \left( \eps, \frac{1}{n}, \frac{1}{k} \right) = \frac{\eps^{15}}{10^{41} n^{68} k^{167}}$.
\end{theorem}

%%%%%%%%%%%%%%%%%%%%%%%%%%%%%%%%%%%%%%%%%%%%%%%%%
\subsection{Discussion}\label{sec:discussion} %%%
%%%%%%%%%%%%%%%%%%%%%%%%%%%%%%%%%%%%%%%%%%%%%%%%%

Our results cover all previous cases for which a quantitative Gibbard-Satterthwaite theorem has been established before. In particular, the main novelty is that neutrality of the SCF is not assumed, and therefore our results hold for nonneutral SCFs as well, thereby solving the main open problem of Friedgut, Kalai, Keller and Nisan~\cite{friedgut2011quantitative}, and Isaksson, Kindler and Mossel \cite{isaksson2010geometry}. The main message of our results is that the approach of masking manipulation behind computational hardness cannot hide manipulations completely even in the nonneutral setting.

\medskip

\noindent\textbf{Importance of nonneutrality.} While neutrality seems like a very natural assumption, there are multiple reasons why removing this assumption is important:
\begin{itemize}
\item \textbf{Anonymity vs.\ neutrality.} It is known that there is a conflict between anonymity and neutrality (recall Definitions~\ref{def:anonymity} and~\ref{def:neutrality}). In particular, there are some combinations of $n$ and $k$ when there exists no SCF which is both anonymous and neutral.
\begin{theorem}\cite[Chapter 2.4.]{moulin1983strategy}
There exists a SCF on $n$ voters and $k$ alternatives which is both anonymous and neutral if and only if $k$ cannot be written as the sum of (non-trivial) divisors of $n$.
\end{theorem}

The difficulty comes from rules governing tie-breaking. Consider the following example: suppose $n=k=2$, i.e.\ we have two voters, voter 1 and voter 2, and two alternatives, $a$ and $b$. Suppose further (w.l.o.g.) that when voter 1 prefers $a$ over $b$ and voter 2 prefers $b$ over $a$ then the outcome is $a$. What should the outcome be when voter 1 prefers $b$ over $a$ and voter 2 prefers $a$ over $b$? By anonymity the outcome should be $a$ for this configuration as well, but by neutrality the outcome should be $b$.

Most common voting rules (plurality, Borda count, etc.) break ties in an anonymous way, and therefore they cannot be neutral as well (or can only be neutral for special values of $n$ and $k$). See Moulin~\cite[Chapter 2.4.]{moulin1983strategy} for more on anonymity and neutrality.

\item \textbf{Nonneutrality in virtual elections.} As mentioned before, voting manipulation is a serious issue in artificial intelligence and computer science as well, where virtual elections are becoming more and more popular as a tool in preference aggregation (see the survey by Faliszewski and Procaccia \cite{faliszewski2010ai}). For example, consider web (meta-)search engines (see e.g.\ Dwork et al.\ \cite{dwork2001rank}), where one inputs a query and the possible outcomes (``alternatives'') are the web pages (with the various search engines acting as ``voters''). Here, due to various restrictions, neutrality is not a natural assumption. For example, there can be language-related restrictions: if one searches in English then the top-ranked webpage will also be in English; or safety-related restrictions: if one searches in child-safe mode, then the top-ranked webpage cannot have adult content. These restrictions imply that the appropriate aggregating function cannot be neutral.

\item \textbf{Nonneutrality in real-life elections.} Although not a common occurrence, there have been cases in real-life elections when a candidate is on the ballot, but is actually ineligible---she cannot win the election no matter what. In such a case the SCF is necessarily nonneutral.

In a recent set of local elections in Philadelphia there were actually three such occurences~\cite{cbsphilly}: one of the candidates for the one open Municipal Court slot was not a lawyer, which is a prerequisite for someone elected to this position; another judicial candidate received a court order to leave the race; finally, in the race for a district seat in Philadelphia, one of the candidates had announced that he is abandoning his candidacy; yet all three of them remained on the respective ballots.

A more curious story is that of the New York State Senate elections in 2010, where the name of a dead man appeared on the ballot (he received 828 votes)~\cite{bronxdead}.
\end{itemize}

\noindent\textbf{A quantitative Gibbard-Satterthwaite theorem for one voter.} %As mentioned in Section \ref{sec:results}, a major part of the work in proving Theorem \ref{thm:k_refined} is devoted to understanding functions of a single voter, and hence we chose to formulate the result of these arguments separately in Theorem \ref{thm:quant_GS_1voter}, a quantitative Gibbard-Satterthwaite theorem for one voter. We note that this is a new result, which has not been studied in the literature before. 
A major part of the work in proving Theorem~\ref{cor:k_refined_truenonmanip} is devoted to understanding functions of a single voter, essentially proving a quantitative Gibbard-Satterthwaite theorem for one voter. This can be formulated as follows.
\begin{theorem}\label{thm:quant_GS_1voter}
Suppose $f : S_k \to \left[ k \right]$ is a SCF on $n=1$ voter and $k \geq 3$ alternatives which satisfies $\Dist \left( f, \NONMANIP \right) \geq \eps$. Then
\begin{equation}\label{eq:main}
\p \left( \sigma \in M \left( f \right) \right) \geq \p \left( \sigma \in M_3 \left( f \right) \right) \geq p \left( \eps, \frac{1}{k} \right),
\end{equation}
for some polynomial $p$, where $\sigma \in S_k$ is selected uniformly. In particular, we show a lower bound of $\frac{\eps^3}{10^5 k^{16}}$.
\end{theorem}

% Dobzinski and Procaccia \cite{dobzinski2008frequent} proved a quantitative Gibbard-Satterthwaite theorem for two voters, but they assume that the SCF is \emph{Pareto-optimal}, i.e.\ if all voters rank alternative $a$ above $b$, then $b$ is not elected. This is a reasonably weak assumption; however, it implies that every alternative is elected in at least a $1/k^2$ fraction of the ranking profiles (when both voters rank this alternative at the top), an assumption which is not needed and we do not assume. Moreover, there is only one SCF on one voter which is Pareto-optimal (the one which always outputs the top-ranked alternative), so we cannot make this assumption when we want to understand functions of a single voter.

We note that this is a new result, which has not been studied in the literature before. 

Dobzinski and Procaccia \cite{dobzinski2008frequent} proved a quantitative Gibbard-Satterthwaite theorem for two voters, assuming that the SCF is \emph{Pareto optimal}, i.e.\ if all voters rank alternative $a$ above $b$, then $b$ is not elected. The assumption of Pareto optimality is natural in the context of classical social choice, but it is a very strong assumption in the context of quantitative social choice. For one, it implies that every alternative is elected with probability at least $1/k^2$. Second, for one voter, there exists a unique Pareto optimal SCF, while the number of nonmanipulable SCFs is exponential in $k$. The assumption also prevents applying the result of Dobzinski and Procaccia to SCFs obtained from a SCF on many voters when the votes of all voters but two are fixed (since even if the original SCF is Pareto optimal, the restricted function may not be so). In our proof we often deal with such restricted SCFs (where the votes of all but one or two voters are fixed), and this is also what led us to our quantitative Gibbard-Satterthwaite theorem for one voter.

\medskip

\noindent\textbf{On $\mathbf{\NONMANIP}$ versus $\mathbf{\overline{\NONMANIP}}$.} The quantitative Gibbard-Satterthwaite theorems of Friedgut, Kalai, Keller and Nisan~\cite{friedgut2008elections,friedgut2011quantitative}, and Isaksson, Kindler and Mossel \cite{isaksson2010geometry} involve the distance of a SCF from $\overline{\NONMANIP}$. Any SCF that is not in $\overline{\NONMANIP}$ is manipulable (by the Gibbard-Satterthwaite theorem), but as some SCFs in $\overline{\NONMANIP}$ are manipulable as well, ideally a quantitative Gibbard-Satterthwaite theorem would involve the distance of a SCF from the set of (truly) nonmanipulable SCFs, $\NONMANIP$. Theorem~\ref{cor:k_refined_truenonmanip} addresses this concern, as it involves the distance of a SCF from $\NONMANIP$. This is done via the following reduction theorem that implies that whenever one has a quantitative Gibbard-Satterthwaite theorem involving $\Dist \left( f, \overline{\NONMANIP} \right)$, this can be turned into a quantitative Gibbard-Satterthwaite theorem involving $\Dist \left( f, \NONMANIP \right)$.

\begin{theorem}\label{thm:TRUENONMANIP}
Suppose $f$ is a SCF on $n$ voters and $k \geq 3$ alternatives for which $\Dist \left( f, \overline{\NONMANIP} \right) \leq \alpha$. Then
either
\begin{equation}\label{eq:true_NONMANIP}
\Dist \left( f, \NONMANIP \right) < 100 n^4 k^8 \alpha^{1/3}
\end{equation}
or
\begin{equation}\label{eq:many_manip}
\p \left( \sigma \in M \left( f \right) \right) \geq \p \left( \sigma \in M_3 \left( f \right) \right) \geq \alpha.
\end{equation}
\end{theorem}

The proof of this result also uses Theorem~\ref{thm:quant_GS_1voter}, our quantitative Gibbard-Satterthwaite theorem for one voter.

\medskip

\noindent\textbf{A note on our quantitative bounds.} The lower bounds on the probability of manipulation derived in Theorems~\ref{cor:k_refined_truenonmanip}, \ref{thm:quant_GS_1voter}, and various results along the way, are not tight. Moreover, we do not believe that our techniques allow us to obtain tight bounds. Consequently, we did not try to optimize these bounds, but rather focused on the qualitative result: obtaining polynomial bounds.

%%%%%%%%%%%%%%%%%%%%%%%%%%%%%%%%%%%%%%%%%%%%%%%%%%%%%%%%%%%%%
\subsection{Proof techniques}\label{sec:proof_techniques} %%%
%%%%%%%%%%%%%%%%%%%%%%%%%%%%%%%%%%%%%%%%%%%%%%%%%%%%%%%%%%%%%

In our proof we combine ideas from both Friedgut, Kalai, Keller and Nisan~\cite{friedgut2008elections,friedgut2011quantitative} and Isaksson, Kindler and Mossel~\cite{isaksson2010geometry}, and in addition we use a reverse hypercontractivity lemma that was applied in the proof of a quantitative version of Arrow's theorem by Mossel~\cite{mossel2009quantitative}. (Reverse hypercontractivity was originally proved and discussed by Borell~\cite{borell1982positivity}, and was first applied by Mossel, O'Donnell, Regev, Steif and Sudakov~\cite{mossel2006non}.) Our techniques most closely resemble those of Isaksson et al.~\cite{isaksson2010geometry}; here the authors used a variant of the canonical path method to show the existence of a large interface where three bodies touch. Our goal is also to come to this conclusion, but we do so via different methods. %We present the main ideas in the case of $k=3$ alternatives.

%By our reduction theorem (Theorem~\ref{thm:TRUENONMANIP}) it is enough to assume that $\Dist \left( f, \overline{\NONMANIP} \right) \geq \eps$...

We first present our techniques that achieve a lower bound for the probability of manipulation that involves factors of $\frac{1}{k!}$ (see Theorem~\ref{thm:k_bdd} in Section~\ref{sec:k_bdd}), and then describe how a refined approach leads to a lower bound which has inverse polynomial dependence on $k$ (see Theorem~\ref{thm:k_refined} in Section~\ref{sec:manip_ref}).

\medskip

\noindent\textbf{Rankings graph and applying the original Gibbard-Satterthwaite theorem.} As in Isaksson et al.~\cite{isaksson2010geometry}, think of the graph $G = \left( V, E \right)$ having vertex set $V = S_k^n$, the set of all ranking profiles, and let $\left( \sigma, \sigma' \right) \in E$ if and only if $\sigma$ and $\sigma'$ differ in exactly one coordinate. The SCF $f : S_k^n \to \left[k \right]$ naturally partitions $V$ into $k$ subsets. Since every manipulation point must be on the boundary between two such subsets, we are interested in the size of such boundaries.

For two alternatives $a$ and $b$, and voter $i$, denote by $B_i^{a,b}$ the boundary between $f^{-1} \left( a \right)$ and $f^{-1}\left( b \right)$ in voter $i$. A lemma from Isaksson et al.~\cite{isaksson2010geometry} tells us that at least two of the boundaries are large; in the following assume that these are $B_1^{a,b}$ and $B_2^{a,c}$. Now if a ranking profile $\sigma$ lies on \emph{both} of these boundaries, then applying the original Gibbard-Satterthwaite theorem to the restricted SCF on two voters where we fix all coordinates of $\sigma$ except the first two, we get that there must exist a manipulation point which agrees with $\sigma$ in all but the first two coordinates. Consequently, if we can show that the \emph{intersection} of the boundaries $B_1^{a,b}$ and $B_2^{a,c}$ is large, then we have many manipulation points.
%Consequently, if we can show that the number of ranking profiles that lie on \emph{both} $B_1^{a,b}$ and $B_2^{a,c}$ is large, then we have many manipulation points.

\medskip

\noindent\textbf{Fibers and reverse hypercontractivity.} In order to have more ``control'' over what is happening at the boundaries, we partition the graph further---this idea is due to Friedgut et al.~\cite{friedgut2008elections,friedgut2011quantitative}. Given a ranking profile $\sigma$ and two alternatives $a$ and $b$, $\sigma$ induces a \emph{vector of preferences} $x^{a,b}\left( \sigma \right) \in \left\{-1,1\right\}^{n}$ between $a$ and $b$. For a vector $z^{a,b} \in \left\{-1,1\right\}^n$ we define the \emph{fiber with respect to preferences between $a$ and $b$}, denoted by $F\left( z^{a,b} \right)$, to be the set of ranking profiles for which the vector of preferences between $a$ and $b$ is $z^{a,b}$. We can then partition the vertex set $V$ into such fibers, and work inside each fiber separately. Working inside a specific fiber is advantageous, because it gives us the extra knowledge of the vector of preferences between $a$ and $b$.

% In particular, we can also partition the \emph{boundaries} according to the fibers (boundaries consist of edges, but for each edge we can consider one its endpoints, e.g.\ for $B_1^{a,b}$, consider the endpoints at which the outcome of $f$ is $a$). Working inside a specific fiber is advantageous, because it gives us the extra knowledge of the vector of preferences between $a$ and $b$.

We distinguish two types of fibers: large and small. We say that a fiber w.r.t.\ preferences between $a$ and $b$ is \emph{large} if almost all of the ranking profiles in this fiber lie on the boundary $B_1^{a,b}$, and \emph{small} otherwise. Now since the boundary $B_1^{a,b}$ is large, either there is big mass on the large fibers w.r.t.\ preferences between $a$ and $b$ or big mass on the small fibers. This holds analogously for the boundary $B_2^{a,c}$ and fibers w.r.t.\ preferences between $a$ and $c$.

Consider the case when there is big mass on the large fibers of both $B_1^{a,b}$ and $B_2^{a,c}$. Notice that for a ranking profile $\sigma$, being in a fiber w.r.t.\ preferences between $a$ and $b$ only depends on the vector of preferences between $a$ and $b$, $x^{a,b}\left( \sigma \right)$, which is a uniform bit vector. Similarly, being in a fiber w.r.t.\ preferences between $a$ and $c$ only depends on $x^{a,c}\left( \sigma \right)$. Moreover, we know the exact correlation between the coordinates of $x^{a,b} \left( \sigma \right)$ and $x^{a,c} \left( \sigma \right)$, and it is in exactly this setting where \emph{reverse hypercontractivity} applies (see Lemma \ref{lem:invhypcontr} for a precise statement), and shows that the \emph{intersection} of the large fibers of $B_1^{a,b}$ and $B_2^{a,c}$ is also large. Finally, by the definition of a large fiber it follows that the intersection of the \emph{boundaries} $B_1^{a,b}$ and $B_2^{a,c}$ is large as well, and we can finish the argument using the Gibbard-Satterthwaite theorem as above.

To deal with the case when there is big mass on the \emph{small} fibers of $B_1^{a,b}$ we use various isoperimetric techniques, including the canonical path method developed for this problem by Isaksson et al.~\cite{isaksson2010geometry}. In particular, we use the fact that for a small fiber for $B_1^{a,b}$, the size of the boundary of $B_1^{a,b}$ in the small fiber is comparable to the size of $B_1^{a,b}$ in the small fiber itself, up to polynomial factors.

\medskip

\noindent\textbf{A refined geometry.} Using this approach with the rankings graph above, our bound includes $\frac{1}{k!}$ factors (see Theorem \ref{thm:k_bdd} in Section \ref{sec:k_bdd}). In order to obtain inverse polynomial dependence on $k$ (as in Theorem \ref{thm:k_refined} in Section~\ref{sec:manip_ref}), we use a refined approach, similar to that in Isaksson et al.~\cite{isaksson2010geometry}. Instead of the rankings graph outlined above, we use an underlying graph with a different edge structure: $\left( \sigma, \sigma' \right) \in E$ if and only if $\sigma$ and $\sigma'$ differ in exactly one coordinate, and in this coordinate they differ by a single adjacent transposition. In order to prove the refined result, we need to show that the geometric and combinatorial quantities such as boundaries and manipulation points are roughly the same in the refined graph as in the original rankings graph. In particular, this is where we need to analyze carefully functions of one voter, and ultimately prove a quantitative Gibbard-Satterthwaite theorem for one voter.

\medskip

\noindent A detailed outline of the proof of Theorem \ref{thm:k_bdd} can be found in Section \ref{sec:proof_outline_1}, while an overview of the changes in the proof to get to Theorem \ref{thm:k_refined}, in which the lower bound has inverse polynomial dependence on $k$, is given in Section \ref{sec:k_refined_overview}. Both of these theorems assume that $\Dist \left( f, \overline{\NONMANIP} \right) \geq \eps$; the reduction theorem (Theorem~\ref{thm:TRUENONMANIP}) then allows us to weaken this assumption to $\Dist \left( f, \NONMANIP \right) \geq \eps$ (the proof of Theorem~\ref{thm:TRUENONMANIP} can be found at the end of the paper, in Section~\ref{sec:TRUENONMANIP}).

%%%%%%%%%%%%%%%%%%%%%%%%%%%%%%%%%%%%%%%%%%
\subsection{Organization of the paper} %%%
%%%%%%%%%%%%%%%%%%%%%%%%%%%%%%%%%%%%%%%%%%

The rest of the paper is outlined as follows. We introduce necessary preliminaries (definitions and previous technical results) in Section \ref{sec:prelims}. We then proceed by proving Theorem \ref{thm:k_bdd} in Section \ref{sec:k_bdd}, which is weaker than Theorem~\ref{cor:k_refined_truenonmanip} in two aspects: first, the condition $\Dist \left( f, \NONMANIP \right) \geq \eps$ is replaced with the stronger condition $\Dist \left( f, \overline{\NONMANIP} \right) \geq \eps$, and second, we allow factors of $\frac{1}{k!}$ in our lower bounds for $\p \left( \sigma \in M \left( f \right) \right)$. We continue by explaining the necessary modifications we have to make in the refined setting to get inverse polynomial dependence on $k$ in Section \ref{sec:k_refined_overview}. Additional preliminaries necessary for the proofs of Theorems~\ref{thm:k_refined}, \ref{thm:quant_GS_1voter} and \ref{thm:TRUENONMANIP} are in Section~\ref{sec:k_refined_prelims}, while the remaining sections contain the proofs of these theorems. We prove Theorem~\ref{thm:quant_GS_1voter} in Section~\ref{sec:1voter}, Theorem~\ref{thm:k_refined} in Section~\ref{sec:manip_ref}, and Theorem~\ref{thm:TRUENONMANIP} and Theorem~\ref{cor:k_refined_truenonmanip} in Section~\ref{sec:TRUENONMANIP}. Finally we conclude with some open problems in Section~\ref{sec:open}.

%%%%%%%%%%%%%%%%%%%%%%%%%%%%%%%%%%%%%%%
%%%%%%%%%%%%%%%%%%%%%%%%%%%%%%%%%%%%%%%
%%%%%%%%%%%%%%%%%%%%%%%%%%%%%%%%%%%%%%%

%%%%%%%%%%%%%%%%%%%%%%%%%%%%%%%%%%%%%%%%%%%%%%%%%%%%%%%%%%%%%%%%%%%%%%%%%%%%%%%%%%%%%%%%%%
\section{Preliminaries: definitions and previous technical results}\label{sec:prelims} %%%
%%%%%%%%%%%%%%%%%%%%%%%%%%%%%%%%%%%%%%%%%%%%%%%%%%%%%%%%%%%%%%%%%%%%%%%%%%%%%%%%%%%%%%%%%%

In this section we introduce some definitions and previous technical results that we use throughout the paper. In particular, we introduce everything that is needed to prove our first theorem, Theorem \ref{thm:k_bdd}, and in Section \ref{sec:k_refined_prelims} we present additional preliminaries needed for the rest of the theorems.

%%%%%%%%%%%%%%%%%%%%%%%%%%%%%%%%%%%%%%%%%%
\subsection{Boundaries and influences} %%%
%%%%%%%%%%%%%%%%%%%%%%%%%%%%%%%%%%%%%%%%%%

For a general graph $G = \left( V, E \right)$, and a subset of the vertices $A \subseteq G$, we define the \emph{edge boundary} of $A$ as
\[
\partial_{e} \left( A \right) = \left\{ \left( u, v \right) \in E : u \in A, v \notin A \right\}.
\]
We also define the \emph{boundary} (or vertex boundary) of  a subset of the vertices $A \subseteq G$ to be the set of vertices in $A$ which have a neighbor that is not in $A$:
\[
\partial \left( A \right) = \left\{ u \in A : \text{ there exists } v \notin A \text{ such that } \left( u, v \right) \in E \right\}.
\]
If $u \in \partial \left( A \right)$, we also say that $u$ is \emph{on} the edge boundary of $A$.

As discussed in Section \ref{sec:proof_techniques}, we can view the ranking profiles (which are elements of $S_k^n$) as vertices of a graph---the rankings graph---where two vertices are connected by an edge if they differ in exactly one coordinate. The SCF $f$ naturally partitions the vertices of this graph into $k$ subsets, depending on the value of $f$ at a given vertex. Clearly, a manipulation point can only be on the edge boundary of such a subset, and so it is important to study these boundaries. In this spirit, we introduce the following definitions.

\begin{definition}[Boundaries]
For a given SCF $f$ and a given alternative $a \in \left[k \right]$, we define
\[
H^{a} \left( f \right) = \left\{ \sigma \in S_k^n : f\left( \sigma \right) = a \right\},
\]
the set of ranking profiles where the outcome of the vote is $a$. The edge boundary of this set is denoted by $B^a \left( f \right)$ : $B^a \left( f \right) = \partial_e \left( H^a \left( f \right) \right)$. This boundary can be partitioned: we say that the edge boundary of $H^a \left( f \right)$ in the direction of the $i^{\text{th}}$ coordinate is
\[
B_i^{a} \left( f \right) = \left\{ \left( \sigma, \sigma' \right) \in B^a \left( f \right) :  \sigma_i \neq \sigma'_i \right\}.
\]
The boundary $B^a \left( f \right)$ can be therefore written as $B^a \left( f \right) = \cup_{i=1}^{n} B_i^a \left( f \right)$. We can also define the boundary between two alternatives $a$ and $b$ in the direction of the $i^{\text{th}}$ coordinate:
\[
B_i^{a,b} \left( f \right) = \left\{ \left( \sigma, \sigma' \right) \in B_i^a \left( f \right) : f\left( \sigma' \right) = b \right\}.
\]
We also say that $\sigma \in B_i^a \left( f \right)$ is \emph{on} the boundary $B_i^{a,b} \left( f \right)$ if there exists $\sigma'$ such that $\left( \sigma, \sigma' \right) \in B_i^{a,b} \left( f \right)$.
\end{definition}

\begin{definition}[Influences]
We define the \emph{influence} of the $i^{\text{th}}$ coordinate on $f$ as
\[
\Inf_i \left( f \right) = \p \left( f\left( \sigma \right) \neq f \left( \sigma^{\left( i \right)} \right) \right) = \p \left( \left( \sigma, \sigma^{\left( i \right)} \right) \in \cup_{a = 1}^{k} B_i^a \left( f \right) \right),
\]
where $\sigma$ is uniform on $S_{k}^{n}$ and $\sigma^{\left( i \right)}$ is obtained from $\sigma$ by rerandomizing the $i^{\text{th}}$ coordinate. Similarly, we define the influence of the $i^{\text{th}}$ coordinate with respect to a single alternative $a \in \left[k \right]$ or a pair of alternatives $a,b \in \left[k\right]$ as
\[
\Inf_i^{a} \left( f \right) = \p \left( f\left( \sigma \right) = a, f \left( \sigma^{\left( i \right)} \right) \neq a \right) = \p \left( \left( \sigma, \sigma^{\left( i \right)} \right) \in B_i^a \left( f \right) \right) ,
\]
and
\[
\Inf_i^{a,b} \left( f \right) = \p \left( f\left( \sigma \right) = a, f \left( \sigma^{\left( i \right)} \right) = b \right) = \p \left( \left( \sigma, \sigma^{\left( i \right)} \right) \in B_i^{a,b} \left( f \right) \right),
\]
respectively.
\end{definition}
Clearly
\[
\Inf_i \left( f \right) = \sum_{a=1}^{k} \Inf_i^a \left( f \right) = \sum_{a,b \in \left[k \right]: a\neq b} \Inf_{i}^{a,b} \left( f \right).
\]

Most of the time the specific SCF $f$ will be clear from the context, in which case we omit the dependence on $f$, and write simply $B^a \equiv B^a \left( f \right)$, $B_i^a \equiv B_i^a \left( f \right)$, etc.

%%%%%%%%%%%%%%%%%%%%%%%%%%%%%%%%%
\subsection{Large boundaries} %%%
%%%%%%%%%%%%%%%%%%%%%%%%%%%%%%%%%

The following lemma from Isaksson, Kindler and Mossel \cite[Lemma 3.1.]{isaksson2010geometry} shows that there are some boundaries which are large (in the sense that they are only inverse polynomially small in $n$, $k$ and $\eps^{-1}$)---our task is then to find many manipulation points on these boundaries.
\begin{lemma} \label{lem:boundaries1}
  Fix $k \ge 3$ and
  $f \colon S_k^n \to \left[k\right]$
  satisfying $\Dist(f, \overline{\NONMANIP}) \ge \eps$.
%  \question{Does it hold if the distance to \emph{non-manipulable} functions is $ \ge \eps$?}
  Then there exist distinct $i,j \in [n]$
  and $\{a,b\},\{c,d\} \subseteq [k]$ such that $c \notin \{a,b\}$ and
  \begin{equation}
    \Inf_i^{a,b} (f) \ge \frac{2\eps}{n k^2 (k-1)}
    \quad \text{ and } \quad
    \Inf_j^{c,d} (f) \ge \frac{2\eps}{n k^2 (k-1)}.
  \end{equation}
\end{lemma}

%%%%%%%%%%%%%%%%%%%%%%%%%%%%%%%%%%%%%%%%%%%%%%
\subsection{General isoperimetric results} %%%
%%%%%%%%%%%%%%%%%%%%%%%%%%%%%%%%%%%%%%%%%%%%%%

Our rankings graph is the Cartesian product of $n$ complete graphs on $k!$ vertices. We therefore use isoperimetric results on products of graphs---see \cite{harper2004global} for an overview. In particular, the edge-isoperimetric problem on the product of complete graphs was originally solved by Lindsey in 1964 \cite{lindsey1964assignment}:
\begin{theorem}[Lindsey \cite{lindsey1964assignment}]\label{thm:isoLindsey}
The edge-isoperimetric problem on $K_{n_1} \times K_{n_2} \times \dots \times K_{n_d}$, a product of complete graphs with $n_1 \leq n_2 \leq \dots \leq n_d$, has lexicographic nested solutions.
\end{theorem}
\begin{corollary}\label{cor:isoLindsey}
If $A \subseteq K_{k} \times \dots \times K_{k}$ ($n$ copies) and $\left| A \right| \leq \left( 1 - \frac{1}{k} \right) k^n$, then $\left| \partial_{e} \left( A \right) \right| \geq \left| A \right|$.
\end{corollary}

%%%%%%%%%%%%%%%%%%%%%%%%%%%%%%%%%%%%%%%%%
\subsection{Fibers}\label{sec:fibers} %%%
%%%%%%%%%%%%%%%%%%%%%%%%%%%%%%%%%%%%%%%%%

In our proof we need to partition the graph even further---this idea is due to Friedgut, Kalai, Keller, and Nisan~\cite{friedgut2008elections,friedgut2011quantitative}.
\begin{definition}\label{def:xab}
For a ranking profile $\sigma \in S_k^n$ define the vector
\[
x^{a,b} \equiv x^{a,b} \left( \sigma \right) = \left( x_1^{a,b} \left( \sigma \right), \dots, x_n^{a,b} \left( \sigma \right) \right)
\]
of preferences between $a$ and $b$, where $x_i^{a,b} \left( \sigma \right) = 1$ if $a \stackrel{\sigma_i}{>} b$, i.e.\ voter $i$ prefers $a$ over $b$, and $x_i^{a,b} \left( \sigma \right) = -1$ otherwise.
\end{definition}
\begin{definition}[Fibers]\label{def:fibers}
For a pair of alternatives $a,b \in \left[k \right]$ and a vector $z^{a,b} \in \left\{-1,1\right\}^n$, write
\[
F \left( z^{a,b} \right) := \left\{ \sigma : x^{a,b} \left( \sigma \right) = z^{a,b} \right\}.
\]
We call the $F \left( z^{a,b} \right)$ \emph{fibers} with respect to preferences between $a$ and $b$.
\end{definition}
So for any pair of alternatives $a,b$, we can partition the ranking profiles according to its fibers:
\[
S_k^n = \bigcup_{z^{a,b} \in \left\{ - 1, 1 \right\}^{n}} F \left( z^{a,b} \right).
\]

Given a SCF $f$, for any pair of alternatives $a,b \in \left[ k \right]$ and $i \in \left[ n \right]$, we can also partition the boundary $B_i^{a,b} \left( f \right)$ according to its fibers. There are multiple, slightly different ways of doing this, but for our purposes the following definition is most useful. Define
\[
B_i \left( z^{a,b} \right) := \left\{ \sigma \in F \left( z^{a,b} \right) : f \left( \sigma \right) = a, \text{ and there exists } \sigma' \text{ s.t.\ } \left( \sigma, \sigma' \right) \in B_i^{a,b} \right\},
\]
where we omit the dependence of $B_i \left( z^{a,b} \right)$ on $f$. So $B_i \left( z^{a,b} \right) \subseteq F \left( z^{a,b} \right)$ is the set of vertices on the given fiber for which the outcome is $a$ and which lies on the boundary between $a$ and $b$ in direction $i$. We call the sets of the form $B_i \left( z^{a,b} \right)$ \emph{fibers for the boundary $B_i^{a,b}$} (again omitting the dependence on $f$ of both sets).

We now distinguish between small and large fibers for the boundary $B_i^{a,b}$.
\begin{definition}[Small and large fibers]\label{def:lg_fibers}
We say that the fiber $B_i \left( z^{a,b} \right)$ is \emph{large} if
\begin{equation}\label{eq:lg_fbr_def}
\p \left( \sigma \in B_i \left( z^{a,b} \right) \, \middle| \, \sigma \in F \left( z^{a,b} \right) \right) \geq 1 - \frac{\eps^{3}}{4 n^3 k^9},
\end{equation}
and \emph{small} otherwise.

We denote by $\Lg \left( B_i^{a,b} \right)$ the union of large fibers for the boundary $B_i^{a,b}$, i.e.\
\[
\Lg \left( B_i^{a,b} \right) := \left\{ \sigma : B_i \left( x^{a,b} \left( \sigma \right) \right) \text{ is a large fiber, and } \sigma \in B_i \left( x^{a,b} \left( \sigma \right) \right) \right\}
\]
and similarly, we denote by $\Sm \left( B_i^{a,b} \right)$ the union of small fibers.
\end{definition}
We remark that what is important is that the fraction appearing on the right hand side of \eqref{eq:lg_fbr_def} is a polynomial of $\frac{1}{n}$, $\frac{1}{k}$ and $\eps$---the specific polynomial in this definition is the end result of the computation in the proof.

Finally, for a voter $i$ and a pair of alternatives $a,b \in \left[k\right]$, we define
\[
F_i^{a,b} := \left\{ \sigma : B_i \left( x^{a,b} \left( \sigma \right) \right) \text{ is a large fiber} \right\}.
\]
So this means that
\begin{equation}\label{eq:Aiab}
\p \left( \sigma \in \cup_{z^{a,b}} B_i \left( z^{a,b} \right) \, \middle| \, \sigma \in F_i^{a,b} \right) \geq 1 - \frac{\eps^{3}}{4 n^3 k^9}.
\end{equation}

%%%%%%%%%%%%%%%%%%%%%%%%%%%%%%%%%%%%%%%%%
\subsection{Boundaries of boundaries} %%%
%%%%%%%%%%%%%%%%%%%%%%%%%%%%%%%%%%%%%%%%%

Finally, we also look at boundaries of boundaries. In particular, for a given vector $z^{a,b}$ of preferences between $a$ and $b$, we can think of the fiber $F \left( z^{a,b} \right)$ as a subgraph of the original rankings graph. When we write $\partial \left( B_i \left( z^{a,b} \right) \right)$, we mean the boundary of $B_i \left( z^{a,b} \right)$ in the subgraph of the rankings graph induced by the fiber $F \left( z^{a,b} \right)$. That is,
\[
\partial \left( B_i \left( z^{a,b} \right) \right) = \{ \sigma \in B_i \left( z^{a,b} \right) : \exists\ \pi \in F\left( z^{a,b} \right) \setminus B_i \left( z^{a,b} \right) \text{ s.t. } \sigma \text{ and } \pi \text{ differ in exactly one coordinate} \}.
\]

%%%%%%%%%%%%%%%%%%%%%%%%%%%%%%%%%%%%%%%%%%%
\subsection{Reverse hypercontractivity} %%%
%%%%%%%%%%%%%%%%%%%%%%%%%%%%%%%%%%%%%%%%%%%

We use the following lemma about reverse hypercontractivity from Mossel \cite{mossel2009quantitative}.
\begin{lemma}\label{lem:invhypcontr}
Let $x,y \in \left\{-1,1\right\}^n$ be distributed uniformly and $\left\{ \left( x_i, y_i \right) \right\}_{i=1}^{n}$ are independent. Assume that $\E \left( x_i \right) = \E \left( y_i \right) = 0$ for all $i$ and that $\left| \E\left( x_i y_i \right) \right| \leq \rho$. Let $B_1, B_2 \subset \left\{ -1, 1 \right\}^n$ be two sets and assume that
\[
\p \left( B_1 \right) \geq e^{-\alpha^2}, \qquad \p \left( B_2 \right) \geq e^{-\beta^2}.
\]
Then
\[
\p \left( x \in B_1, y \in B_2 \right) \geq \exp \left( - \frac{\alpha^2 + \beta^2 + 2\rho \alpha \beta}{1 - \rho^2} \right).
\]
In particular, if $\p \left( B_1 \right) \geq \eps$ and $\p \left( B_2 \right) \geq \eps$, then
\[
\p \left( x \in B_1, y \in B_2 \right) \geq \eps^{\frac{2}{1-\rho}}.
\]
\end{lemma}

%%%%%%%%%%%%%%%%%%%%%%%%%%%%%%%%%%%%%%%%%%%%%%%%%%%%%%%%%%%%%%%%%%%%%%%%%%%%%%%%%
\subsection{Dictators and miscellaneous definitions}\label{sec:dict_misc_def} %%%
%%%%%%%%%%%%%%%%%%%%%%%%%%%%%%%%%%%%%%%%%%%%%%%%%%%%%%%%%%%%%%%%%%%%%%%%%%%%%%%%%

For a ranking profile $\sigma = \left( \sigma_1, \dots, \sigma_n \right)$ we sometimes write $\sigma_{-i}$ for the collection of all coordinates except the $i^{\text{th}}$ coordinate, i.e.\ $\sigma = \left( \sigma_i, \sigma_{-i} \right)$. Furthermore, we sometimes distinguish two coordinates, e.g.\ we write $\sigma = \left( \sigma_1, \sigma_i, \sigma_{-\left\{1, i\right\}} \right)$.

\begin{definition}[Induced SCF on one coordinate]\label{def:induced_SCF}
Let $f_{\sigma_{-i}}$ denote the SCF on one voter induced by $f$ by fixing all voter preferences except the $i^{\text{th}}$ one according to $\sigma_{-i}$. I.e.,
\[
f_{\sigma_{-i}} \left( \cdot \right) := f\left( \cdot, \sigma_{-i} \right).
\]
\end{definition}
Recall Definition \ref{def:dict_subset} of a dictator on a subset.
\begin{definition}[Ranking profiles giving dictators on a subset]
For a coordinate $i$ and a subset of alternatives $H \subseteq \left[k \right]$, define
\[
D_i^H := \left\{ \sigma_{-i} : f_{\sigma_{-i}} \left( \cdot \right) \equiv \tp_H \left( \cdot \right) \right\}.
\]

Also, for a pair of alternatives $a$ and $b$, define
\[
D_i \left( a,b \right) := \bigcup_{H : \left\{a,b \right\} \subseteq H, \left|H \right| \geq 3} D_i^H.
\]
\end{definition}

%%%%%%%%%%%%%%%%%%%%%%%%%%%%%%%%%%%%%%%%%%%%%%%%%%%%%%%%%%%%%%%%%%%%%%%%%%%%%%%%%%%%%%%%%%%%%%%%%%%
\section{Inverse polynomial manipulability for a fixed number of alternatives}\label{sec:k_bdd} %%%
%%%%%%%%%%%%%%%%%%%%%%%%%%%%%%%%%%%%%%%%%%%%%%%%%%%%%%%%%%%%%%%%%%%%%%%%%%%%%%%%%%%%%%%%%%%%%%%%%%%

Our goal in this section is to demonstrate the proof techniques described in Section~\ref{sec:proof_techniques}. We prove here the following theorem (Theorem~\ref{thm:k_bdd} below), which is weaker than our main theorem, Theorem~\ref{cor:k_refined_truenonmanip}, in two aspects: first, the condition $\Dist \left( f, \NONMANIP \right) \geq \eps$ is replaced with the stronger condition $\Dist \left( f, \overline{\NONMANIP} \right) \geq \eps$, and second, we allow factors of $\frac{1}{k!}$ in our lower bounds for $\p \left( \sigma \in M \left( f \right) \right)$. The advantage is that the proof of this statement is relatively simpler. We move on to getting a lower bound with polynomial dependence on $k$ in the following sections, and finally we replace the condition $\Dist \left( f, \overline{\NONMANIP} \right) \geq \eps$ with $\Dist \left( f, \NONMANIP \right) \geq \eps$ in Section~\ref{sec:TRUENONMANIP}.

\begin{theorem}\label{thm:k_bdd}
Suppose we have $n \geq 2$ voters, $k \geq 3$ alternatives, and a SCF $f : S_k^n \to \left[k\right]$ satisfying $\Dist \left( f, \overline{\NONMANIP} \right) \geq \eps$. Then
\begin{equation}\label{eq:manip_general}
\p \left( \sigma \in M \left( f \right) \right)\geq p \left( \eps, \frac{1}{n}, \frac{1}{k!} \right),
\end{equation}
for some polynomial $p$, where $\sigma \in S_k^n$ is selected uniformly. In particular, we show a lower bound of $\frac{\eps^5}{4 n^7 k^{12} \left(k!\right)^4}$.

An immediate consequence is that
\begin{equation*}
\p \left( \left( \sigma, \sigma' \right) \text{ is a manipulation pair for } f \right)\geq q \left( \eps, \frac{1}{n}, \frac{1}{k!} \right),
\end{equation*}
for some polynomial $q$, where $\sigma \in S_k^n$ is selected uniformly, and $\sigma'$ is obtained from $\sigma$ by uniformly selecting a coordinate $i \in \left\{1, \dots, n \right\}$ and resetting the $i^{\text{th}}$ coordinate to a random preference. In particular, the specific lower bound for $\p \left( \sigma \in M \left( f \right) \right)$ implies that we can take $q \left( \eps, \frac{1}{n}, \frac{1}{k} \right) = \frac{\eps^5}{4 n^8 k^{12} \left(k!\right)^5}$.
\end{theorem}

First we provide an overview of the proof of Theorem~\ref{thm:k_bdd} in Section~\ref{sec:proof_outline_1}. In this overview we use adjectives such as ``big'', and ``not too small'' to describe probabilities---here these are all synonymous with ``has probability at least an inverse polynomial of $n$, $k!$, and $\eps^{-1}$''.

%%%%%%%%%%%%%%%%%%%%%%%%%%%%%%%%%%%%%%%%%%%%%%%%%%%%%%%%%%%%%
\subsection{Overview of proof}\label{sec:proof_outline_1} %%%
%%%%%%%%%%%%%%%%%%%%%%%%%%%%%%%%%%%%%%%%%%%%%%%%%%%%%%%%%%%%%

The tactic in proving Theorem \ref{thm:k_bdd} is roughly the following:
\begin{itemize}
\item By Lemma \ref{lem:boundaries1}, we know that there are at least two boundaries which are big. W.l.o.g.\ we can assume that these are either $B_1^{a,b}$ and $B_2^{a,c}$, or $B_1^{a,b}$ and $B_2^{c,d}$ with $\left\{a,b\right\} \cap \left\{ c,d \right\} = \emptyset$. Our proof works in both cases, but we continue the outline of the proof assuming the former case---this is the more interesting case, since the latter case has been solved already by Isaksson et al.~\cite{isaksson2010geometry}.
\item We partition $B_1^{a,b}$ according to its fibers based on the preferences between $a$ and $b$ of the $n$ voters, just like as described in Section \ref{sec:prelims}. Similarly for $B_2^{a,c}$ and preferences between $a$ and $c$.
\item As in Section \ref{sec:prelims}, we can distinguish small and large fibers for these two boundaries. Now since $B_1^{a,b}$ is big, either the mass of small fibers, or the mass of large fibers is big. Similarly for $B_2^{a,c}$.
\item Suppose first that there is big mass on large fibers in both $B_1^{a,b}$ and $B_2^{a,c}$. In this case the probability of our random ranking $\sigma$ being in $F_1^{a,b}$ is big, and similarly for $F_2^{a,c}$. Being in $F_1^{a,b}$ only depends on the vector $x^{a,b} \left( \sigma \right)$ of preferences between $a$ and $b$, and similarly being in $F_2^{a,c}$ only depends on the vector $x^{a,c}\left( \sigma \right)$ of preferences between $a$ and $c$. We know the correlation between $x^{a,b}\left( \sigma \right)$ and $x^{a,c} \left( \sigma \right)$ and hence we can apply reverse hypercontractivity (Lemma \ref{lem:invhypcontr}), which tells us that the probability that $\sigma$ lies in both $F_1^{a,b}$ and $F_2^{a,c}$ is big as well. If $\sigma \in F_1^{a,b}$, then voter 1 is pivotal between alternatives $a$ and $b$ with big probability, and similarly if $\sigma \in F_2^{a,c}$, then voter 2 is pivotal between alternatives $a$ and $c$ with big probability. So now we have that the probability that both voter 1 is pivotal between $a$ and $b$ and voter 2 is pivotal between $a$ and $c$ is big, and in this case the Gibbard-Satterthwaite theorem tells us that there is a manipulation point which agrees with this ranking profile in all except for perhaps the first two coordinates. So there are many manipulation points.
\item Now suppose that the mass of small fibers in $B_1^{a,b}$ is big. By isoperimetric theory, the size of the boundary of every small fiber is comparable (same order up to $\text{poly}^{-1} \left( \eps^{-1}, n, k! \right)$ factors) to the size of the small fiber. Consequently, the total size of the boundaries of small fibers is comparable to the total size of small fibers, which in this case has to be big.

We then distinguish two cases: either we are on the boundary of a small fiber in the first coordinate, or some other coordinate. If $\sigma$ is on the boundary of a small fiber in some coordinate $j \neq 1$, then the Gibbard-Satterthwaite theorem tells us that there is a manipulation point which agrees with $\sigma$ in all coordinates except perhaps in coordinates 1 and $j$. If our ranking profile $\sigma$ is on the boundary of a small fiber in the first coordinate, then either there exists a manipulation point which agrees with $\sigma$ in all coordinates except perhaps the first, or the SCF on one voter that we obtain from $f$ by fixing the votes of voters 2 through $n$ to be $\sigma_{-1}$ must be a dictator on some subset of the alternatives. So either we get sufficiently many manipulation points this way, or for many votes of voters 2 through $n$, the restricted SCF obtained from $f$ by fixing these votes is a dictator on coordinate 1 on some subset of the alternatives.

Finally, to deal with dictators on the first coordinate, we look at the boundary of the dictators. Since $\Dist \left( f, \overline{\NONMANIP} \right) \geq \eps$, the boundary is big, and we can also show that there is a manipulation point near every boundary point.
\item If the mass of small fibers in $B_2^{a,c}$ is big, then we can do the same thing for this boundary.
\end{itemize}

%%%%%%%%%%%%%%%%%%%%%%%%%%%%%%%%%%%%
\subsection{Division into cases} %%%
%%%%%%%%%%%%%%%%%%%%%%%%%%%%%%%%%%%%

For the remainder of Section \ref{sec:k_bdd}, let us fix the number of voters $n \geq 2$, the number of alternatives $k \geq 3$, and the SCF $f$, which satisfies $\Dist \left( f, \overline{\NONMANIP} \right) \geq \eps$. Accordingly, we typically omit the dependence of various sets (e.g.\ boundaries between two alternatives) on $f$.

Our starting point is Lemma \ref{lem:boundaries1}. W.l.o.g.\ we may assume that the two boundaries that the lemma gives us have $i=1$ and $j=2$, so the lemma tells us that
\begin{equation*}
\p \left( \left( \sigma, \sigma^{\left( 1 \right)} \right) \in B_1^{a,b} \right) \geq \frac{2 \eps}{n k^3},
\end{equation*}
where $\sigma$ is uniform on the ranking profiles, and $\sigma^{\left(1\right)}$ is obtained by rerandomizing the first coordinate. This also means that
\begin{equation*}
\p \left( \sigma \in \cup_{z^{a,b}} B_1 \left( z^{a,b} \right) \right) \geq \frac{2 \eps}{n k^3},
\end{equation*}
and similar inequalities hold for the boundary $B_2^{c,d}$. The following lemma is an immediate corollary.
\begin{lemma}\label{lem:cases_k}
Either
\begin{equation}\label{eq:sm_fbr}
\p \left( \sigma \in \Sm \left( B_1^{a,b} \right) \right) \geq \frac{\eps}{n k^3}
\end{equation}
or
\begin{equation}\label{eq:lg_fbr}
\p \left( \sigma \in \Lg \left( B_1^{a,b} \right) \right) \geq \frac{\eps}{n k^3},
\end{equation}
and the same can be said for the boundary $B_2^{c,d}$.
\end{lemma}
We distinguish cases based upon this: either \eqref{eq:sm_fbr} holds, or \eqref{eq:sm_fbr} holds for the boundary $B_2^{c,d}$, or \eqref{eq:lg_fbr} holds for both boundaries. We only need one boundary for the small fiber case, and we need both boundaries only in the large fiber case. So in the large fiber case we must differentiate between two cases: whether $d\in \left\{a,b\right\}$ or $d \notin \left\{ a, b \right\}$. First of all, in the $d\notin \left\{a,b \right\}$ case the problem of finding a manipulation point with not too small (i.e.\ inverse polynomial in $n$, $k!$ and $\eps^{-1}$) probability has already been solved in \cite{isaksson2010geometry}. But moreover, we will see that if $d \notin \left\{a,b\right\}$ then the large fiber case cannot occur---so this method of proof works as well.

In the rest of the section we first deal with the large fiber case, and then with the small fiber case.

%%%%%%%%%%%%%%%%%%%%%%%%%%%%%%%%%%%%%%%%%%%%%%%%%%%%%%%%%%%%%
\subsection{Big mass on large fibers}\label{sec:k_lg_fbr} %%%
%%%%%%%%%%%%%%%%%%%%%%%%%%%%%%%%%%%%%%%%%%%%%%%%%%%%%%%%%%%%%

We now deal with the case when
\begin{equation}\label{eq:k_lg_fbr_ab}
\p \left( \sigma \in \Lg \left( B_1^{a,b} \right) \right) \geq \frac{\eps}{n k^3}
\end{equation}
and also
\begin{equation}\label{eq:k_lg_fbr_cd}
\p \left( \sigma \in \Lg \left( B_2^{c,d} \right) \right) \geq \frac{\eps}{n k^3}.
\end{equation}
As mentioned before, we must differentiate between two cases: whether $d\in \left\{a,b\right\}$ or $d \notin \left\{ a, b \right\}$.

%%%%%%%%%%%%%%%%%%%%%%%%%%
\subsubsection{Case 1} %%%
%%%%%%%%%%%%%%%%%%%%%%%%%%

Suppose $d\in \left\{a,b\right\}$, in which case we may assume w.l.o.g.\ that $d = a$.
\begin{lemma}\label{lem:lg_fbr_1}
If
\begin{equation}\label{eq:k_lg_fbr_ab2}
\p \left( \sigma \in \Lg \left( B_1^{a,b} \right) \right) \geq \frac{\eps}{n k^3} \qquad \text{ and } \qquad \p \left( \sigma \in \Lg \left( B_2^{a,c} \right) \right) \geq \frac{\eps}{n k^3},
\end{equation}
% and also
% \begin{equation}\label{eq:k_lg_fbr_ac}
% \p \left( \sigma \in \Lg \left( B_2^{a,c} \right) \right) \geq \frac{\eps}{n k^3},
% \end{equation}
then
\begin{equation}\label{eq:k_lg_fbr_manip}
\p \left( \sigma \in M \right) \geq \frac{\eps^{3}}{2 n^{3} k^{9} \left(k!\right)^2}.
\end{equation}
\end{lemma}
\begin{proof}
By \eqref{eq:k_lg_fbr_ab2} we have that %and \eqref{eq:k_lg_fbr_ac} we have that
% \begin{align*}
% \p \left( \sigma \in F_1^{a,b} \right) &\geq \frac{\eps}{n k^3},\\
% \p \left( \sigma \in F_2^{a,c} \right) &\geq \frac{\eps}{n k^3}.
% \end{align*}
\[
 \p \left( \sigma \in F_1^{a,b} \right) \geq \frac{\eps}{n k^3} \qquad \text{ and } \qquad \p \left( \sigma \in F_2^{a,c} \right) \geq \frac{\eps}{n k^3}.
\]
We know that $\left|\E\left( x_i^{a,b} \left( \sigma \right) x_i^{a,c} \left( \sigma \right) \right)\right| = 1/3$, and so by reverse hypercontractivity (Lemma \ref{lem:invhypcontr}) we have that
\begin{equation}\label{eq:inv_hyp_contr_k}
\p \left( \sigma \in F_1^{a,b} \cap F_2^{a,c} \right) \geq \frac{\eps^{3}}{n^3 k^9}.
\end{equation}
Recall that we say that $\sigma$ is \emph{on} the boundary $B_1^{a,b}$ if there exists $\sigma'$ such that $\left( \sigma, \sigma' \right) \in B_1^{a,b}$. If $\sigma \in F_1^{a,b}$, then with big probability $\sigma$ is on the boundary $B_1^{a,b}$, and if $\sigma \in F_2^{a,c}$, then with big probability $\sigma$ is on the boundary $B_2^{a,c}$. Using this and \eqref{eq:inv_hyp_contr_k} we can show that the probability of $\sigma$ lying on both the boundary $B_1^{a,b}$ and the boundary $B_2^{a,c}$ is big. Then we are done, because if $\sigma$ lies on both $B_1^{a,b}$ and $B_2^{a,c}$, then by the Gibbard-Satterthwaite theorem there is a $\hat{\sigma}$ which agrees with $\sigma$ on the last $n-2$ coordinates, and which is a manipulation point, and there can be at most $\left( k! \right)^{2}$ ranking profiles that give the same manipulation point. Let us do the computation:
\begin{multline*}
\p \left( \sigma \text{ on } B_1^{a,b}, \sigma \text{ on } B_2^{a,c} \right) \geq \p \left( \sigma \text{ on } B_1^{a,b}, \sigma \text{ on } B_2^{a,c}, \sigma \in F_1^{a,b} \cap F_2^{a,c} \right)\\
\geq \p \left( \sigma \in F_1^{a,b} \cap F_2^{a,c} \right) - \p \left( \sigma \in F_1^{a,b} \cap F_2^{a,c}, \sigma \text{ not on } B_1^{a,b} \right)- \p \left( \sigma \in F_1^{a,b} \cap F_2^{a,c}, \sigma \text{ not on } B_2^{a,c} \right).
\end{multline*}
The first term is bounded below via \eqref{eq:inv_hyp_contr_k}, while the other two terms can be bounded using \eqref{eq:Aiab}:
\[
\p \left( \sigma \in F_1^{a,b} \cap F_2^{a,c}, \sigma \text{ not on } B_1^{a,b} \right) \leq \p \left( \sigma \in F_1^{a,b}, \sigma \text{ not on } B_1^{a,b} \right) \leq \p \left( \sigma \text{ not on } B_1^{a,b} \, \middle| \, \sigma \in F_1^{a,b} \right)\leq \frac{\eps^{3}}{4 n^3 k^9},
\]
and similarly for the other term. Putting everything together gives us
\begin{equation*}
\p \left( \sigma \text{ on } B_1^{a,b}, \sigma \text{ on } B_2^{a,c} \right) \geq \frac{\eps^{3}}{2 n^3 k^9},
\end{equation*}
which by the discussion above implies \eqref{eq:k_lg_fbr_manip}.
\end{proof}

%%%%%%%%%%%%%%%%%%%%%%%%%%%%%%%%%%%%%%%%%%%%%%%%%%%%%%
\subsubsection{Case 2}\label{sec:lg_fbr_gen_case2} %%%
%%%%%%%%%%%%%%%%%%%%%%%%%%%%%%%%%%%%%%%%%%%%%%%%%%%%%%

\begin{lemma}\label{lem:lg_fbr_2}
If $d\notin \left\{a,b\right\}$, then \eqref{eq:k_lg_fbr_ab} and \eqref{eq:k_lg_fbr_cd} cannot hold simultaneously.
\end{lemma}
\begin{proof}
Suppose on the contrary that \eqref{eq:k_lg_fbr_ab} and \eqref{eq:k_lg_fbr_cd} do both hold. Then
% \begin{align*}
% \p \left( \sigma \in F_1^{a,b} \right) &\geq \frac{\eps}{n k^3},\\
% \p \left( \sigma \in F_2^{c,d} \right) &\geq \frac{\eps}{n k^3}
% \end{align*}
\[
 \p \left( \sigma \in F_1^{a,b} \right) \geq \frac{\eps}{n k^3} \qquad \text{ and } \qquad \p \left( \sigma \in F_2^{c,d} \right) \geq \frac{\eps}{n k^3}
\]
as before. Since $\left\{ a, b \right\} \cap \left\{ c, d \right\} = \emptyset$, $\left\{ \sigma \in F_1^{a,b} \right\}$ and $\left\{ \sigma \in F_2^{c,d} \right\}$ are independent events, and so
\[
\p \left( \sigma \in F_1^{a,b} \cap F_2^{c,d} \right) = \p \left( \sigma \in F_1^{a,b} \right) \p \left( \sigma \in F_2^{c,d} \right) \geq  \frac{\eps^{2}}{n^2 k^6}.
\]
In the same way as before, by the definition of large fibers this implies that
\begin{equation*}
\p \left( \sigma \text{ on } B_1^{a,b}, \sigma \text{ on } B_2^{c,d} \right) \geq \frac{\eps^{2}}{2 n^2 k^6} > 0,
\end{equation*}
but it is clear that
\begin{equation*}
\p \left( \sigma \text{ on } B_1^{a,b}, \sigma \text{ on } B_2^{c,d} \right) = 0,
\end{equation*}
since $\sigma$ on $B_1^{a,b}$ and on $B_2^{c,d}$ requires $f\left( \sigma \right) \in \left\{a,b\right\} \cap \left\{c,d\right\} = \emptyset$. So we have reached a contradiction.
\end{proof}

%%%%%%%%%%%%%%%%%%%%%%%%%%%%%%%%%%%%%%%%%%%%%%%%%%%%%%%%%%%%%
\subsection{Big mass on small fibers}\label{sec:k_sm_fbr} %%%
%%%%%%%%%%%%%%%%%%%%%%%%%%%%%%%%%%%%%%%%%%%%%%%%%%%%%%%%%%%%%

We now deal with the case when \eqref{eq:sm_fbr} holds, i.e.\ when we have a big mass on the small fibers for the boundary $B_1^{a,b}$.
We formalize the ideas of the outline described in Section \ref{sec:proof_outline_1} in a series of statements.

First, we want to formalize that the boundaries of the boundaries are big, when we are on a small fiber.
\begin{lemma}\label{lem:comparable_boundaries_k}
Fix coordinate 1 and the pair of alternatives $a,b$. Let $z^{a,b}$ be such that $B_1 \left( z^{a,b} \right)$ is a small fiber for $B_1^{a,b}$. Then, writing $B \equiv B_1 \left( z^{a,b} \right)$, we have
\[
\left| \partial_e \left( B \right) \right| \geq \frac{\eps^3}{4 n^3 k^9} \left| B \right|
\]
and
\begin{equation}\label{eq:iso_small_fib}
\p \left( \sigma \in \partial \left( B \right) \right) \geq \frac{\eps^3}{2 n^4 k^9 k!} \p \left( \sigma \in B \right),
\end{equation}
where both the edge boundary $\partial_e \left( B \right)$ and the boundary $\partial \left( B \right)$ are with respect to the induced subgraph $F\left( z^{a,b} \right)$, which is isomorphic to $K_{k!/2}^{n}$, the Cartesian product of $n$ complete graphs of size $k!/2$.
\end{lemma}
\begin{proof}
We use Corollary \ref{cor:isoLindsey} with $k$ replaced by $k! / 2$ and the set $A$ being either $B$ or $B^c := F \left( z^{a,b} \right) \setminus B$. Suppose first that $\left| B \right| \leq \left( 1 - \frac{2}{k!} \right)\left( k! / 2 \right)^n$. Then $\left| \partial_e \left( B \right) \right| \geq \left| B \right|$. Suppose now that $\left| B \right| > \left( 1 - \frac{2}{k!} \right)\left( k! / 2 \right)^n$. Since we are in the case of a small fiber, we also know that $\left| B \right| \leq \left( 1 - \frac{\eps^3}{4 n^3 k^9} \right) \left(k! / 2 \right)^n$. Consequently, we get
\[
\left| \partial_e \left( B \right) \right| = \left| \partial_e \left( B^c \right) \right| \geq \left| B^c \right| \geq \frac{\eps^3}{4 n^3 k^9} \left| B \right|,
\]
which proves the first claim.

A ranking profile in $F \left( z^{a,b} \right)$ has $\left( k! / 2 - 1 \right) n \leq n k! / 2$ neighbors in $F \left( z^{a,b} \right)$, which then implies \eqref{eq:iso_small_fib}.
\end{proof}

\begin{corollary}\label{cor:comparable_boundaries_k}
If \eqref{eq:sm_fbr} holds, then
\[
\p\left( \sigma \in \bigcup_{z^{a,b}} \partial \left( B_1 \left( z^{a,b} \right) \right) \right) \geq \frac{\eps^4}{2 n^5 k^{12} k!}.
\]
\end{corollary}

\begin{proof}
Using the previous lemma and \eqref{eq:sm_fbr} we have
\begin{align*}
\p \left( \sigma \in \bigcup_{z^{a,b}} \partial \left( B_1 \left( z^{a,b} \right) \right) \right) &= \sum_{z^{a,b}} \p \left( \sigma \in \partial \left( B_1 \left( z^{a,b} \right) \right) \right)\geq \sum_{z^{a,b} : B_1 \left( z^{a,b} \right) \subseteq \Sm \left( B_1^{a,b} \right)} \p \left( \sigma \in \partial \left( B_1 \left( z^{a,b} \right) \right) \right)\\
&\geq \sum_{z^{a,b} : B_1 \left( z^{a,b} \right) \subseteq \Sm \left( B_1^{a,b} \right)} \frac{\eps^3}{2 n^4 k^9 k!} \p \left( \sigma \in  B_1 \left( z^{a,b} \right) \right) = \frac{\eps^3}{2 n^4 k^9 k!} \p \left( \sigma  \in \Sm \left( B_1^{a,b} \right) \right)\\
&\geq \frac{\eps^4}{2 n^5 k^{12} k!}. \qedhere
\end{align*}
\end{proof}

Next, we want to find manipulation points on the boundaries of boundaries.

\begin{lemma}\label{lem:manip_on_bdry}
Suppose the ranking profile $\sigma$ is on the boundary of a fiber for $B_1^{a,b}$, i.e.\
\[
\sigma \in \bigcup_{z^{a,b}} \partial \left( B_1 \left( z^{a,b} \right) \right).
\]
Then either $\sigma_{-1} \in D_1 \left( a, b \right)$, or there exists a manipulation point $\hat{\sigma}$ which differs from $\sigma$ in at most two coordinates, one of them being the first coordinate.
% \begin{itemize}
% \item either $\sigma_{-1} \in D_1 \left( a, b \right)$,
% \item or there exists a manipulation point $\hat{\sigma}$ which differs from $\sigma$ in at most two coordinates, one of them being the first coordinate.
% \end{itemize}
\end{lemma}

\begin{proof}
First of all, by our assumption that $\sigma$ is on the boundary of a fiber for $B_1^{a,b}$, we know that $\sigma \in B_1 \left( z^{a,b} \right)$ for some $z^{a,b}$, which means that there exists a ranking profile $\sigma' = \left( \sigma'_1, \sigma_{-1} \right)$ such that $\left( \sigma, \sigma' \right) \in B_1^{a,b}$. We may assume $a \stackrel{\sigma_1}{>} b$ and $b \stackrel{\sigma'_1}{>} a$, or else either $\sigma$ or $\sigma'$ is a manipulation point.

Now since $\sigma \in \partial \left( B_1 \left( z^{a,b} \right) \right)$ we also know that there exists a ranking profile $\pi = \left( \pi_j, \sigma_{-j} \right) \in F \left( z^{a,b} \right) \setminus B_1 \left( z^{a,b} \right)$ for some $j \in \left[k \right]$. We distinguish two cases: $j \neq 1$ and $j=1$.

\textbf{Case 1:} $\mathbf{j \neq 1}$\textbf{.} What does it mean for $\pi = \left( \pi_j, \sigma_{-j} \right)$ to be on the same fiber as $\sigma$, but for $\pi$ to not be in $B_1 \left( z^{a,b} \right)$? First of all, being on the same fiber means that $\sigma_j$ and $\pi_j$ both rank $a$ and $b$ in the same order. Now $\pi \notin B_1 \left( z^{a,b} \right)$ means that
\begin{itemize}
\item either $f\left( \pi \right) \neq a$;
\item or $f \left( \pi \right) = a$ and $f\left( \pi'_1, \pi_{-1} \right) \neq b$ for every $\pi'_1 \in S_k$.
\end{itemize}
If $f\left( \pi \right) = b$, then either $\sigma$ or $\pi$ is a manipulation point, since the order of $a$ and $b$ is the same in both $\sigma_j$ and $\pi_j$ (since $\sigma$ and $\pi$ are on the same fiber).

Suppose $f\left( \pi \right) = c \notin \left\{ a, b \right\}$. Then we can define a SCF function on two coordinates by fixing all coordinates except coordinates 1 and $j$ to agree with the respective coordinates of $\sigma$---letting coordinates 1 and $j$ vary we get a SCF function on two coordinates which takes on at least three values ($a$, $b$, and $c$), and does not only depend on one coordinate. Now applying the Gibbard-Satterthwaite theorem we get that this SCF on two coordinates has a manipulation point, which means that our original SCF $f$ has a manipulation point which agrees with $\sigma$ in all coordinates except perhaps in coordinates 1 and $j$.

So the final case is that $f \left( \pi \right) = a$ and $f\left( \pi'_1, \pi_{-1} \right) \neq b$ for every $\pi'_1 \in S_k$. In particular for $\tilde{\pi} := \left( \sigma'_1, \pi_{-1} \right) = \left( \pi_j, \sigma'_{-j} \right)$ we have $f\left( \tilde{\pi} \right) \neq b$. Now if $f\left( \tilde{\pi} \right) = a$ then either $\sigma'$ or $\tilde{\pi}$ is a manipulation point, since the order of $a$ and $b$ is the same in both $\sigma'_j = \sigma_j$ and $\pi_j$. Finally, if $f\left( \tilde{\pi} \right) = c \notin \left\{ a, b \right\}$, then we can apply the Gibbard-Satterthwaite theorem just like in the previous paragraph.

\textbf{Case 2:} $\mathbf{j = 1}$\textbf{.} We can again ask: what does it mean for $\pi = \left( \pi_1, \sigma_{-1} \right)$ to be on the same fiber as $\sigma$, but for $\pi$ to not be in $B_1 \left( z^{a,b} \right)$? First of all, being on the same fiber means that $\sigma_1$ and $\pi_1$ both rank $a$ and $b$ in the same order (namely, as discussed at the beginning, ranking $a$ above $b$, or else we have a manipulation point). Now $\pi \notin B_1 \left( z^{a,b} \right)$ means that
\begin{itemize}
\item either $f\left( \pi \right) \neq a$;
\item or $f \left( \pi \right) = a$ and $f\left( \pi'_1, \pi_{-1} \right) \neq b$ for every $\pi'_1 \in S_k$.
\end{itemize}
However, we know that $f \left( \sigma' \right) = b$ and that $\sigma'$ is of the form $\sigma' = \left( \sigma'_1, \sigma_{-1} \right) = \left( \sigma'_1, \pi_{-1} \right)$, and so the only way we can have $\pi \notin B_1 \left( z^{a,b} \right)$ is if $f\left( \pi \right) \neq a$.

If $f\left( \pi \right) = b$, then $\pi$ is a manipulation point, since $a \stackrel{\pi_1}{>} b$ and $f\left( \sigma \right) = a$.

So the remaining case is if $f\left( \pi \right) = c \notin \left\{a,b\right\}$. This means that $f_{\sigma_{-1}}$ (see Definition \ref{def:induced_SCF}) takes on at least three values. Denote by $H \subseteq \left[ k \right]$ the range of $f_{\sigma_{-1}}$. Now either $\sigma_{-1} \in D_1^{H} \subseteq D_1 \left( a,b \right)$, or there exists a manipulation point $\hat{\sigma}$ which agrees with $\sigma$ in every coordinate except perhaps the first. %This concludes the proof.
\end{proof}

Finally, we need to deal with dictators on the first coordinate.

\begin{lemma}\label{lem:dictators}
Assume that $\Dist \left( f, \overline{\NONMANIP} \right) \geq \eps$. We have that either
\[
\p \left( \sigma_{-1}  \in D_1 \left( a, b \right) \right) \leq \frac{\eps^4}{4 n^5 k^{12} k!},
\]
or
\begin{equation}\label{eq:manip_with_dictators}
\p \left( \sigma \in M \right) \geq \frac{\eps^{5}}{4 n^7 k^{12} \left(k!\right)^4}.
\end{equation}
\end{lemma}

\begin{proof}
Suppose $\p \left( \sigma_{-1} \in D_1 \left( a,b \right) \right) \geq \frac{\eps^{4}}{4 n^5 k^{12} k!}$, which is the same as
\begin{equation}\label{eq:lot_of_dictators}
\sum_{H : \left\{a,b\right\} \subseteq H, \left| H \right| \geq 3} \p \left( \sigma_{-1} \in D_1^{H} \right) \geq  \frac{\eps^{4}}{4 n^5 k^{12} k!}.
\end{equation}
Note that for every $H \subseteq \left[k \right]$ we have
\[
\eps \leq \Dist \left( f, \overline{\NONMANIP} \right) \leq \p \left( f \left( \sigma \right) \neq \tp_H \left( \sigma_1 \right) \right) \leq 1 - \p \left( D_1^H \right), 
\]
and so 
\begin{equation}\label{eq:D1H_not_too_big}
\p \left( D_1^H \right) \leq 1 - \eps.
\end{equation}

The main idea is that \eqref{eq:D1H_not_too_big} implies that the size of the boundary of $D_1^H$ is comparable to the size of $D_1^H$, and if we are on the boundary of $D_1^H$, then there is a manipulation point nearby.

So first let us establish that the size of the boundary of $D_1^H$ is comparable to the size of $D_1^H$. This is done along the same lines as the proof of Lemma \ref{lem:comparable_boundaries_k}.

Notice that $D_1^{H} \subseteq S_k^{n-1}$, where $S_k^{n-1}$ should be thought of as the Cartesian product of $n-1$ copies of the complete graph on $S_k$. We apply Corollary \ref{cor:isoLindsey} with $k$ replaced by $k!$ and with $n-1$ copies, and we see that if $\eps \geq \frac{1}{k!}$, then $\left| \partial_e \left( D_1^H \right) \right| \geq \left| D_1^H \right|$. If $\eps < \frac{1}{k!}$ and $1 - \frac{1}{k!} \leq \p \left( D_1^H \right) \leq 1-\eps$ then
\[
\left| \partial_e \left( D_1^H \right) \right| = \left| \partial_e \left( \left( D_1^H \right)^c \right) \right| \geq \left| \left(D_1^H\right)^c \right| \geq \eps \left| D_1^H \right|.
\]
So in any case we have $\left| \partial_e \left( D_1^H \right) \right| \geq \eps \left| D_1^H \right|$. Since $\sigma_{-1}$ has $\left( n - 1 \right) \left(k! - 1\right) \leq n k!$ neighbors in $S_{k}^{n-1}$, we have that
\[
\p \left( \sigma_{-1} \in \partial \left( D_1^H \right) \right) \geq \frac{\eps}{n k!} \p \left( \sigma_{-1} \in D_1^H \right).
\]
Consequently, by \eqref{eq:lot_of_dictators}, we have
\begin{align*}
\p \left( \sigma_{-1} \in \bigcup_{H : \left\{a,b\right\} \subseteq H, \left|H \right| \geq 3} \partial \left( D_1^H \right) \right) &= \sum_{H : \left\{a,b\right\} \subseteq H, \left|H \right| \geq 3} \p \left( \sigma_{-1} \in \partial \left( D_1^H \right) \right)\\
&\geq \sum_{H : \left\{a,b\right\} \subseteq H, \left|H \right| \geq 3} \frac{\eps}{n k!} \p \left( \sigma_{-1} \in D_1^H \right) \geq \frac{\eps^{5}}{4 n^6 k^{12} \left(k!\right)^2}.
\end{align*}

Next, suppose $\sigma_{-1} \in \partial \left( D_1^H \right)$ for some $H$ such that $\left\{ a, b \right\} \subseteq H, \left| H \right| \geq 3$. We want to show that then there is a manipulation point ``close'' to $\sigma_{-1}$ in some sense. To be more precise: for the manipulation point $\hat{\sigma}$, $\hat{\sigma}_{-1}$ will agree with $\sigma_{-1}$ in all except maybe one coordinate.

If $\sigma_{-1} \in \partial \left( D_1^H \right)$, then there exist  $j \in \left\{2, \dots, n \right\}$ and $\sigma'_j$ such that $\sigma'_{-1} := \left( \sigma'_{j}, \sigma_{-\left\{1, j\right\}} \right) \notin D_1^H$. That is, $f_{\sigma'_{-1}} \left( \cdot \right) \not\equiv \tp_{H} \left( \cdot \right)$. There can be two ways that this can happen---the two cases are outlined below. Denote by $H' \subseteq \left[k\right]$ the range of $f_{\sigma'_{-1}}$.

\textbf{Case 1:} $\mathbf{H' = H}$\textbf{.} In this case we automatically know that there exists a manipulation point $\hat{\sigma}$ such that $\hat{\sigma}_{-1} = \sigma'_{-1}$, and so $\hat{\sigma}_{-1}$ agrees with $\sigma_{-1}$ in all coordinates except coordinate $j$.

\textbf{Case 2:} $\mathbf{H' \neq H}$\textbf{.} W.l.o.g.\ suppose $H' \setminus H \neq \emptyset$, and let $c \in H' \setminus H$. (The other case when $H \setminus H' \neq \emptyset$ works in exactly the same way.) First of all, we may assume that $f_{\sigma'_{-1}} \left( \cdot \right) \equiv \tp_{H'} \left( \cdot \right)$, because otherwise we have a manipulation point just like in Case 1.
 
We can define a SCF on two coordinates by fixing all coordinates except coordinate 1 and $j$ to agree with $\sigma_{-1}$, and varying coordinates 1 and $j$. We know that the outcome takes on at least three different values, since $\sigma_{-1} \in D_1^H$, and $\left| H \right| \geq 3$. 

Now let us show that this SCF is not a function of the first coordinate. Let $\sigma_1$ be a ranking which puts $c$ first, and then $a$. Then $f \left( \sigma_1, \sigma_{-1} \right) = a$, but $f\left( \sigma_1, \sigma'_{-1} \right) = c$, which shows that this SCF is not a function of the first coordinate (since a change in coordinate $j$ can change the outcome).

Consequently, the Gibbard-Satterthwaite theorem tells us that this SCF on two coordinates has a manipulation point, and therefore there exists a manipulation point $\hat{\sigma}$ for $f$ such that $\hat{\sigma}_{-1}$ agrees with $\sigma_{-1}$ in all coordinates except coordinate $j$.

Putting everything together yields \eqref{eq:manip_with_dictators}.
\end{proof}

%%%%%%%%%%%%%%%%%%%%%%%%%%%%%%%%%%%%%%%%%%%%%%%%%%%%%%%%%%%
\subsection{Proof of Theorem \ref{thm:k_bdd} concluded} %%%
%%%%%%%%%%%%%%%%%%%%%%%%%%%%%%%%%%%%%%%%%%%%%%%%%%%%%%%%%%%

\begin{proof}[Proof of Theorem \ref{thm:k_bdd}]
If \eqref{eq:k_lg_fbr_ab} and \eqref{eq:k_lg_fbr_cd} hold, then we are done by Lemmas \ref{lem:lg_fbr_1} and \ref{lem:lg_fbr_2}.

If not, then either \eqref{eq:sm_fbr} holds, or \eqref{eq:sm_fbr} holds for the boundary $B_2^{c,d}$; w.l.o.g.\ assume that \eqref{eq:sm_fbr} holds.

By Corollary \ref{cor:comparable_boundaries_k}, we have
\[
\p\left( \sigma \in \bigcup_{z^{a,b}} \partial \left( B_1 \left( z^{a,b} \right) \right) \right) \geq \frac{\eps^4}{2 n^5 k^{12} k!}.
\]

We may assume that $\p \left( \sigma_{-1} \in D_1 \left( a, b \right) \right) \leq \frac{\eps^4}{4 n^5 k^{12} k!}$, since otherwise we are done by Lemma \ref{lem:dictators}. Consequently, we then have
\[
\p\left( \sigma \in \bigcup_{z^{a,b}} \partial \left( B_1 \left( z^{a,b} \right) \right), \sigma_{-1} \notin D_1 \left( a, b \right) \right) \geq \frac{\eps^4}{4 n^5 k^{12} k!}.
\]
We can then finish our argument using Lemma \ref{lem:manip_on_bdry}:
\[
\p \left( \sigma \in M \right) \geq \frac{1}{n \left( k! \right)^{2}} \p\left( \sigma \in \bigcup_{z^{a,b}} \partial \left( B_1 \left( z^{a,b} \right) \right), \sigma_{-1} \notin D_1 \left( a, b \right) \right) \geq \frac{\eps^4}{4 n^6 k^{12} \left(k!\right)^{3}}. \qedhere
\]
\end{proof}

%%%%%%%%%%%%%%%%%%%%%%%%%%%%%%%%%%%%%%%%%%%%%%%%%%%%%%%%%%%%%%%%%%%%%%%%%%%%
\section{An overview of the refined proof}\label{sec:k_refined_overview} %%%
%%%%%%%%%%%%%%%%%%%%%%%%%%%%%%%%%%%%%%%%%%%%%%%%%%%%%%%%%%%%%%%%%%%%%%%%%%%%

In order to improve on the result of Theorem \ref{thm:k_bdd}---in particular to get rid of the factor of $\frac{1}{\left(k!\right)^4}$---we need to refine the methods used in the previous section. We continue the approach of Isaksson, Kindler and Mossel \cite{isaksson2010geometry}, where the authors first proved a quantitative Gibbard-Satterthwaite theorem for neutral SCFs with a bound involving factors of $\frac{1}{k!}$, and then with a refined method were able to remove these factors.

%Our goal is now to find not just many manipulation points, but many $r$-manipulation points for some small (constant) $r$.

The key to the refined method is to consider the so-called \emph{refined rankings graph} instead of the general rankings graph studied in Section \ref{sec:k_bdd}. The vertices of this graph are again ranking profiles (elements of $S_k^n$), and two vertices are connected by an edge if they differ in exactly one coordinate, and by an adjacent transposition in that coordinate. Again, the SCF $f$ naturally partitions the vertices of this graph into $k$ subsets, depending on the value of $f$ at a given vertex. Clearly a 2-manipulation point can only be on the edge boundary of such a subset in the refined rankings graph, and so it is important to study these boundaries.

One of the important steps of the proof in Section \ref{sec:k_bdd} is creating a configuration where we fix all but two coordinates, and the SCF $f$ takes on at least three values when we vary these two coordinates---then we can define another SCF on two voters and $k$ alternatives which must have a manipulation point by the Gibbard-Satterthwaite theorem. The advantage of the refined rankings graph is that we can create a configuration where we fix all but two coordinates, and in these two coordinates we also fix all but constantly many adjacent alternatives, and the SCF takes on at least three values when we vary these constantly many adjacent alternatives in the two coordinates. Then we can define another SCF on two voters and $r$ alternatives, where $r$ is a small constant, which must have a manipulation point by the Gibbard-Satterthwaite theorem. Since $r$ is a constant, we only lose a constant factor in our estimates, not factors of $\frac{1}{k!}$.

We state the refined result in Theorem~\ref{thm:k_refined}, which we also prove in Section~\ref{sec:manip_ref}. The proof of Theorem \ref{thm:k_refined} follows the outline of the proof of Theorem \ref{thm:k_bdd}: we know that there are at least two refined boundaries which are big (by Isaksson et al.~\cite{isaksson2010geometry}); we partition them according to their fibers; we distinguish small and large fibers; and we consider two cases: the small fiber case and the large fiber case. The ideas in both cases are roughly the same as in Section~\ref{sec:k_bdd}, except the proofs are more involved. There is, however, one major difference in the small fiber case, which is the following.

The difficulty is dealing with the case when we are on the boundary of a small fiber in the first coordinate. Suppose $\sigma = \left( \sigma_1, \sigma_{-1} \right)$ is on such a boundary. We know that there are $k!$ ranking profiles which agree with $\sigma$ in coordinates 2 through $n$. The difficulty comes from the fact that---in order to obtain a polynomial bound in $k$---we are only allowed to look at a polynomial number (in $k$) of these ranking profiles when searching for a manipulation point. If there is an $r$-manipulation point among them for some small constant $r$, then we are done. If this is not the case then $\sigma$ is what we call a \emph{local dictator} on some subset of the alternatives in coordinate 1. We say that $\sigma$ is a local dictator on some subset $H \subseteq \left[k\right]$ of the alternatives in coordinate 1 if the alternatives in $H$ are adjacent in $\sigma_1$, and permuting the alternatives in $H$ in every possible way in the first coordinate, the outcome of the SCF $f$ is always the top-ranked alternative in $H$.

So instead of dealing with dictators on some subset in coordinate 1, as in Section \ref{sec:k_bdd}, we have to deal with \emph{local dictators} on some subset in coordinate 1. This analysis involves essentially only the first coordinate, in essence proving a quantitative Gibbard-Satterthwaite theorem for one voter. As discussed in Section \ref{sec:discussion}, this has not been studied in the literature before, and, moreover, we were not able to utilize previous quantitative Gibbard-Satterthwaite theorems to solve this problem easily. Hence we separate this argument from the rest of the proof of Theorem \ref{thm:k_refined} and formulate a quantitative Gibbard-Satterthwaite theorem for one voter, Theorem \ref{thm:quant_GS_1voter}, which is proven in Section \ref{sec:1voter}. This proof forms the backbone for the proof of Theorem \ref{thm:k_refined}, which is then proven in Section \ref{sec:manip_ref}. % The proof of Theorem \ref{thm:k_refined} mirrors that of Theorem \ref{thm:quant_GS_1voter}, modifying it when necessary to deal with the rest of the coordinates.

%%%%%%%%%%%%%%%%%%%%%%%%%%%%%%%%%%%%%%%%%%%%%%%%%%%%%%%%%%%%%%%%%%%%%%%%%%%%%%%%%%%%%%%%%%%%%%%%%%
\section{Refined rankings graph---introduction and preliminaries}\label{sec:k_refined_prelims} %%%
%%%%%%%%%%%%%%%%%%%%%%%%%%%%%%%%%%%%%%%%%%%%%%%%%%%%%%%%%%%%%%%%%%%%%%%%%%%%%%%%%%%%%%%%%%%%%%%%%%

This section presents the definitions and results that are needed for the rest of the paper, in the proofs of Theorems \ref{thm:k_refined}, \ref{thm:quant_GS_1voter} and \ref{thm:TRUENONMANIP}. 

%%%%%%%%%%%%%%%%%%%%%%%%%%%%%%%%%%%%%%%%%%%%%%%%%%%%%%%%%%%
\subsection{Transpositions, boundaries, and influences} %%%
%%%%%%%%%%%%%%%%%%%%%%%%%%%%%%%%%%%%%%%%%%%%%%%%%%%%%%%%%%%

\begin{definition}[Adjacent transpositions]
Given two elements $a,b \in \left[k \right]$, the adjacent transposition $\left[a:b\right]$ between them is defined as follows. If $\sigma \in S_k$ has $a$ and $b$ adjacent, then $\left[a:b\right] \sigma$ is obtained from $\sigma$ by exchanging $a$ and $b$. Otherwise $\left[a:b\right] \sigma = \sigma$.

We let $T$ denote the set of all $k \left( k - 1 \right) / 2$ adjacent transpositions.

For $\sigma \in S_k^n$, we let $\left[a:b\right]_i \sigma$ denote the ranking profile obtained by applying $\left[a:b\right]$ on the $i^{\text{th}}$ coordinate of $\sigma$ while leaving all other coordinates unchanged.
\end{definition}

\begin{definition}[Boundaries]
For a given SCF $f$ and a given alternative $a \in \left[k \right]$, we define
\[
H^{a} \left( f \right) = \left\{ \sigma \in S_k^n : f\left( \sigma \right) = a \right\},
\]
the set of ranking profiles where the outcome of the vote is $a$. The edge boundary of this set (with respect to the underlying refined rankings graph) is denoted by $B^{a;T} \left( f \right)$ : $B^{a;T} \left( f \right) = \partial_e \left( H^a \left( f \right) \right)$. This boundary can be partitioned: we say that the edge boundary of $H^a \left( f \right)$ in the direction of the $i^{\text{th}}$ coordinate is
\[
B_i^{a;T} \left( f \right) = \left\{ \left( \sigma, \sigma' \right) \in B^{a;T} \left( f \right):  \sigma_i \neq \sigma'_i \right\}.
\]
The boundary $B^a \left( f \right)$ can be therefore written as $B^{a;T} \left( f \right)= \cup_{i=1}^{n} B_i^{a;T} \left( f \right)$. We can also define the boundary between two alternatives $a$ and $b$ in the direction of the $i^{\text{th}}$ coordinate:
\[
B_i^{a,b;T} \left( f \right) = \left\{ \left( \sigma, \sigma' \right) \in B_i^{a;T} \left( f \right) : f\left( \sigma' \right) = b \right\}.
\]
Moreover, we can define the boundary between two alternatives $a$ and $b$ in the direction of the $i^{\text{th}}$ coordinate with respect to the adjacent transposition $z \in T$:
\[
B_i^{a,b;z} \left( f \right) = \left\{ \left( \sigma, \sigma' \right) \in B_i^{a;T} \left( f \right) : \sigma' = z_i \sigma, f\left( \sigma' \right) = b \right\}.
\]
We also say that $\sigma$ is \emph{on} the boundary $B_i^{a,b;z} \left( f \right)$ if $\left( \sigma, z_i \sigma \right) \in B_i^{a,b;z} \left( f \right)$.
Clearly we have
\[
B_i^{a,b;T} \left( f \right) = \bigcup_{z \in T} B_i^{a,b;z} \left( f \right).
\]
\end{definition}

\begin{definition}[Influences]
Given $z \in T$, we define
\begin{align*}
\Inf_i^{a,b;z} \left( f \right) &= \p \left( f \left( \sigma \right) = a, f \left( \sigma^{\left( i \right)} \right) = b \right)\\
\Inf_i^{a;z} \left( f \right) &= \p \left( f \left( \sigma \right) = a, f \left( \sigma^{\left( i \right)} \right) \neq a \right)\\
\Inf_i^{a,b;T} \left( f \right) &= \sum_{z \in T} \Inf_i^{a,b;z} \left( f \right),
\end{align*}
where $\sigma$ is uniformly distributed in $S_k^n$ and $\sigma^{\left(i\right)}$ is obtained from $\sigma$ by rerandomizing the $i^{\text{th}}$ coordinate $\sigma_i$ in the following way: with probability 1/2 we keep it as $\sigma_i$, and otherwise we replace it by $z \sigma_i$.
\end{definition}
Note that for $a \neq b$,
\[
\Inf_i^{a,b;z} \left( f \right) = \frac{1}{2} \p \left( f\left( \sigma \right) = a, f\left( z_i \sigma \right) = b \right) = \frac{1}{2} \frac{\left| B_i^{a,b;z} \left( f \right) \right|}{\left( k! \right)^{n}}.
\]

Again, most of the time the specific SCF $f$ will be clear from the context, in which case we omit the dependence on $f$.

%%%%%%%%%%%%%%%%%%%%%%%%%%%%%%%%%%%%%%%%%%%%%%%%%%%%%%%%%%
\subsection{Manipulation points on refined boundaries} %%%
%%%%%%%%%%%%%%%%%%%%%%%%%%%%%%%%%%%%%%%%%%%%%%%%%%%%%%%%%%

The following two lemmas from Isaksson, Kindler and Mossel \cite{isaksson2010geometry} identify manipulation points on (or close to) these refined boundaries.

\begin{lemma}
  \label{lem:nonManipBoundary}\cite[Lemma 7.1.]{isaksson2010geometry}
  Fix $f : S_k^n \to [k]$,
  distinct $a,b \in [k]$
  and $(\sigma, \pi) \in B_i^{a,b ; T}$.
  Then either $\sigma_i=\adj{a}{b} \pi_i$,
  or
  one of $\sigma$ and $\pi$ is a $2$-manipulation point for $f$.
\end{lemma}

\begin{lemma}
  \label{lem:nonManipTriple}\cite[Lemma 7.2.]{isaksson2010geometry}
  Fix $f : S_k^n \to [k]$
  and points
  $\sigma, \pi, \mu \in S_k^n$ such that
  $(\sigma, \pi) \in B_i^{a,b; T}$,
  $(\mu, \pi) \in B_j^{c,b; T}$
  % $x = \adj{?}{?}_i y$, $z = \adj{?}{?}_j y$ where $f(x)=a,f(y)=b,f(z)=c$ are
  % distinct.
  where $a,b,c$ are distinct and $i \neq j$.
  Then there exists a $3$-manipulation point $\nu \in S_k^n$ for $f$
  such that $\nu_k=\pi_k$ for $k \notin \{i,j\}$
  and $\nu_i$ is equal to $\sigma_i$ or $\pi_i$ except that the position of $c$ may be
  shifted arbitrarily
  and $\nu_j$ is equal to $\mu_j$ or $\pi_j$ except that the position of $a$ may be
  shifted arbitrarily.
\end{lemma}

%%%%%%%%%%%%%%%%%%%%%%%%%%%%%%%%%%%%%%%%%
\subsection{Large refined boundaries} %%%
%%%%%%%%%%%%%%%%%%%%%%%%%%%%%%%%%%%%%%%%%

An essential result that will be our starting point in Section \ref{sec:manip_ref} is the following lemma, again from Isaksson, Kindler and Mossel \cite{isaksson2010geometry}, which shows that there are large refined boundaries (or else we have a lot of 2-manipulation points automatically).
\begin{lemma} \label{lem:boundaries2}\cite[Lemma 7.3.]{isaksson2010geometry}
  Fix $k \ge 3$ and
  $f : S_k^n \to \left[k\right]$
  satisfying $\Dist(f, \overline{\NONMANIP}) \ge \eps$.
  Let $\sigma$ be uniformly selected from $S_k^n$.
  Then either
  \begin{equation}
    \label{eq:neutralPairs3ManipProb}
    \p\left( \sigma \in M_2 \left( f \right) \right) \ge \frac{4\eps}{n k^7},
  \end{equation}
  or there exist distinct $i,j \in [n]$
  and $\{a,b\},\{c,d\} \subseteq [k]$ such that $c \notin \{a,b\}$ and
  \begin{equation}\label{eq:k_inf_ref}
   \Inf_i^{a,b;[a:b]} (f) \ge \frac{2\eps}{n k^7} \quad
    \text{ and } \quad
   \Inf_j^{c,d;[c:d]} (f) \ge \frac{2\eps}{n k^7}.
  \end{equation}
\end{lemma}

%%%%%%%%%%%%%%%%%%%%%%%
\subsection{Fibers} %%%
%%%%%%%%%%%%%%%%%%%%%%%

We again use fibers $F\left( z^{a,b} \right)$ as defined in Definition \ref{def:fibers}. However, we need more than this. We note that the following definitions only apply in Section \ref{sec:manip_ref}, i.e.\ when when we have at least two voters; in Section \ref{sec:1voter}, when we have only one voter, things are simpler.

Given the result of Lemma \ref{lem:boundaries2}, our primary interest is in the boundary $B_i^{a,b;\left[a:b\right]}$. For ranking profiles on this boundary, we know that the alternatives $a$ and $b$ are adjacent in coordinate $i$---so we know more than just the preference between $a$ and $b$ in coordinate $i$. Consequently we would like to divide the set of ranking profiles with $a$ and $b$ adjacent in coordinate $i$ according to the preferences between $a$ and $b$ in all coordinates except coordinate $i$. The following definitions make this precise.

As done in Section \ref{sec:dict_misc_def} for ranking profiles, we can write $x_{-i}^{a,b} \equiv x_{-i}^{a,b} \left( \sigma \right)$ for the vector of preferences between $a$ and $b$ for all coordinates except coordinate $i$, i.e.\ the whole vector of preferences between $a$ and $b$ is $x^{a,b} \left( \sigma \right) = \left( x_{i}^{a,b} \left( \sigma \right), x_{-i}^{a,b} \left( \sigma \right) \right)$.

We can define $F\left( z_{-i}^{a,b} \right)$ analogously to $F\left( z^{a,b} \right)$:
\[
F\left( z_{-i}^{a,b} \right) := \left\{ \sigma : x_{-i}^{a,b} \left( \sigma \right) = z_{-i}^{a,b} \right\}.
\]
We also define the subset of $F\left( z_{-i}^{a,b} \right)$ where $a$ and $b$ are adjacent in coordinate $i$, with $a$ above $b$:
\[
\bar{F}\left( z_{-i}^{a,b} \right) := \left\{ \sigma \in F\left( z_{-i}^{a,b} \right) : a \text{ and } b \text{ are adjacent in coordinate } i, \text{ with } a \text{ above } b \right\}.
\]

Given a SCF $f$, for any pair of alternatives $a, b \in \left[k \right]$ and coordinate $i \in \left[n\right]$, we can also partition the boundary $B_i^{a,b} \left( f \right)$ according to its fibers. There are multiple, slightly different ways of doing this, but for our purposes the following definition is most useful.

Define
\[
B_i \left( z_{-i}^{a,b} \right) := \left\{ \sigma \in \bar{F}\left( z_{-i}^{a,b} \right) : f\left( \sigma \right) = a, f\left( \left[a:b\right]_{i} \sigma \right) = b \right\},
\]
where we omit the dependence of $B_i \left( z_{-i}^{a,b} \right)$ on $f$. We call sets of the form $B_i \left( z_{-i}^{a,b} \right) \subseteq \bar{F} \left( z_{-i}^{a,b} \right)$ \emph{fibers for the boundary $B_i^{a,b;\left[a:b\right]}$}.

We now distinguish between small and large fibers for the boundary $B_i^{a,b;\left[a:b\right]}$.
\begin{definition}[Small and large fibers]
We say that the fiber $B_i \left( z_{-i}^{a,b} \right) \subseteq \bar{F} \left( z_{-i}^{a,b} \right)$ is \emph{large} if
\[
\p \left( \sigma \in B_i \left( z_{-i}^{a,b} \right) \, \middle| \, \sigma \in \bar{F} \left( z_{-i}^{a,b} \right) \right) \geq 1 - \gamma,
\]
where $\gamma = \frac{\eps^3}{10^3 n^3 k^{24}}$, and \emph{small} otherwise.

As before, we denote by $\Lg \left( B_i^{a,b;\left[a:b\right]} \right)$ the union of large fibers for the boundary $B_i^{a,b;\left[a:b\right]}$, i.e.\
\[
\Lg \left( B_i^{a,b;\left[a:b\right]} \right) := \bigcup_{B_i \left( z_{-i}^{a,b} \right) \text{ is a large fiber}} B_i \left( z_{-i}^{a,b} \right),
\]
and similarly, we denote by $\Sm \left( B_i^{a,b;\left[a:b\right]} \right)$ the union of small fibers.
\end{definition}
As in Definition \ref{def:lg_fibers}, we remark that what is important is that $\gamma$ is a polynomial of $\frac{1}{n}$, $\frac{1}{k}$ and $\eps$---the specific polynomial in this definition is the end result of the computation in the proof.

The following definition is used in Section \ref{sec:lg_fbr_ref} in dealing with the large fiber case in the refined setting.
\begin{definition}
For a coordinate $i$ and a pair of alternatives $a$ and $b$, define $F_i^{a,b}$ to be the set of ranking profiles $\sigma$ such that $x^{a,b} \left( \sigma \right)$ satisfies
\[
\p \left( f \left( \tilde{\sigma} \right) = \tp_{\left\{a,b\right\}} \left( \tilde{\sigma}_i \right) \, \middle| \, \tilde{\sigma} \in F \left( x_{-i}^{a,b} \left( \sigma \right) \right) \right) \geq 1 - 2k \gamma.
\]
\end{definition}
Clearly $F_i^{a,b}$ is the union of fibers of the form $F \left( z^{a,b} \right)$, and also $F \left( \left(1, x_{-i}^{a,b} \right) \right) \subseteq F_i^{a,b}$ if and only if $F \left( \left(-1, x_{-i}^{a,b} \right) \right) \subseteq F_i^{a,b}$.

%%%%%%%%%%%%%%%%%%%%%%%%%%%%%%%%%%%%%%%%%%%%%%%%%%%%%%%%%%%%%%%%%
\subsection{Boundaries of boundaries}\label{sec:bdry_of_bdry} %%%
%%%%%%%%%%%%%%%%%%%%%%%%%%%%%%%%%%%%%%%%%%%%%%%%%%%%%%%%%%%%%%%%%

In the refined graph setting, just like in the general rankings graph setting, we also look at boundaries of boundaries.

For a given vector $z_{-i}^{a,b}$ of preferences between $a$ and $b$, we can think of $\bar{F} \left( z_{-i}^{a,b} \right)$ as a subgraph of the original refined rankings graph $S_k^n$, i.e.\ two ranking profiles in $\bar{F} \left( z_{-i}^{a,b} \right)$ are adjacent if they differ by one adjacent transposition in exactly one coordinate. Since both of the ranking profiles are in $\bar{F} \left( z_{-i}^{a,b} \right)$, this adjacent transposition keeps the order of $a$ and $b$ in all coordinates, and moreover it keeps $a$ and $b$ adjacent in coordinate $i$.

We choose to slightly modify this graph: the vertex set is still $\bar{F} \left( z_{-i}^{a,b} \right)$, but we modify the edge set by adding new edges. Suppose $\sigma \in \bar{F} \left( z_{-i}^{a,b} \right)$ and
\[
\sigma_i = 
\begin{pmatrix}
\vdots\\
c\\
a\\
b\\
d\\
\vdots
\end{pmatrix};
\qquad\quad
\sigma'_i =
\begin{pmatrix}
\vdots\\
a\\
b\\
c\\
d\\
\vdots
\end{pmatrix};
\qquad\quad
\sigma''_i =
\begin{pmatrix}
\vdots\\
c\\
d\\
a\\
b\\
\vdots
\end{pmatrix}.
\]
Define in this way $\sigma' = \left( \sigma'_i, \sigma_{-i} \right)$ and $\sigma'' = \left( \sigma''_i, \sigma_{-i} \right)$, and add $\left( \sigma, \sigma' \right)$ and $\left( \sigma, \sigma'' \right)$ to the edge set. So basically, we consider the block of $a$ and $b$ in coordinate $i$ as a single element, and connect two ranking profiles in $\bar{F} \left( z_{-i}^{a,b} \right)$ if they differ in an adjacent transposition in a single coordinate, allowing this transposition to move the block of $a$ and $b$ in coordinate $i$. We call this graph $G\left( z_{-i}^{a,b} \right) = \left( \bar{F} \left( z_{-i}^{a,b} \right), E \left( z_{-i}^{a,b} \right) \right)$, where $E \left( z_{-i}^{a,b} \right)$ is the edge set.

When we write $\partial_e \left( B_i \left( z_{-i}^{a,b} \right) \right)$, we mean the edge boundary of $B_i \left( z_{-i}^{a,b} \right)$ in the graph $G \left( z_{-i}^{a,b} \right)$, and similarly when we write $\partial \left( B_i \left( z_{-i}^{a,b} \right) \right)$, we mean the vertex boundary of $B_i \left( z_{-i}^{a,b} \right)$ in the graph $G \left( z_{-i}^{a,b} \right)$.

%%%%%%%%%%%%%%%%%%%%%%%%%%%%%%%%%%%%%%%%%%%
\subsection{Reverse hypercontractivity} %%%
%%%%%%%%%%%%%%%%%%%%%%%%%%%%%%%%%%%%%%%%%%%

We again use Lemma \ref{lem:invhypcontr} about reverse hypercontractivity.

%%%%%%%%%%%%%%%%%%%%%%%%%%%%%%%%%%%%%%%%%%%%%%%%%%%%%%%%%%%%%%%%%%%%%%%%%%%%
\subsection{Local dictators, conditioning and miscellaneous definitions} %%%
%%%%%%%%%%%%%%%%%%%%%%%%%%%%%%%%%%%%%%%%%%%%%%%%%%%%%%%%%%%%%%%%%%%%%%%%%%%%

In the general rankings graph setting we defined a dictator on a subset of the alternatives, but in the refined rankings graph setting we need to define so-called \emph{local dictators}.

\begin{definition}[Local dictators]
For a coordinate $i$ and a subset of alternatives $H \subseteq \left[k\right]$, define $\LD_i^H$ to be the set of ranking profiles $\sigma$ such that the alternatives in $H$ form an adjacent block in $\sigma_i$, and permuting them among themselves in any order, the outcome of the SCF $f$ is always the top ranked alternative among those in $H$.
If $\sigma \in \LD_i^{H}$, then we call $\sigma$ a local dictator on $H$ in coordinate $i$.

Also, for a pair of alternatives $a$ and $b$, define
\[
\LD_i \left( a, b\right) := \bigcup_{c \notin \left\{ a, b \right\}} \LD_i^{\left\{a,b,c\right\}},
\]
the set of local dictators on three alternatives, two of which are $a$ and $b$, in coordinate $i$.
\end{definition}

In dealing with local dictators, we will condition on the top of a particular coordinate being fixed. We therefore introduce the following notation.
\begin{definition}[Conditioning]\label{def:cond}
For any coordinate $i \in \left[n\right]$ and any vector $\mathbf{v}$ of alternatives we define
\[
\p_i^{\mathbf{v}} \left( \cdot \right) := \p \left( \cdot \, \middle| \, \left( \sigma_i \left( 1 \right), \dots, \sigma_i \left(\left|\mathbf{v} \right| \right) \right) = \mathbf{v} \right),
\]
where $\left| \mathbf{v} \right|$ denotes the length of the vector $\mathbf{v}$. E.g.\ $\p_1^{\left(a \right)} \left( \cdot \right) = \p \left( \cdot \, \middle| \, \sigma_1 \left( 1 \right) = a \right)$ and
\[
 \p_1^{\left( a, b, c \right)} = \p \left( \cdot \, \middle| \, \left( \sigma_1 \left( 1 \right), \sigma_1 \left( 2 \right), \sigma_1 \left( 3 \right) \right) = \left( a, b, c \right) \right).
\]
\end{definition}

We use the following notation in the proof of Theorem \ref{thm:TRUENONMANIP}.
\begin{definition}[Majority function]\label{def:maj}
For a function $f$ whose domain $X$ is finite and whose range is the set $\left\{a,b\right\}$, define $\Maj \left( f \right)$ by
\begin{equation*}
\Maj \left( f \right) =
\begin{cases}
a & \text{ if } \quad \# \left\{ x \in X : f\left( x \right) = a \right\} \geq \# \left\{ x \in X : f\left( x \right) = b \right\},\\
b & \text{ if } \quad \# \left\{ x \in X : f\left( x \right) = a \right\} < \# \left\{ x \in X : f\left( x \right) = b \right\}.
\end{cases}
\end{equation*}
\end{definition}

%%%%%%%%%%%%%%%%%%%%%%%%%%%%%%%%%%%%%%%%%%%%%%%%%%%%%%%%%%%%%%%%%%%%%%%%%%%%%%%%%%%%%%%%
\section{Quantitative Gibbard-Satterthwaite theorem for one voter}\label{sec:1voter} %%%
%%%%%%%%%%%%%%%%%%%%%%%%%%%%%%%%%%%%%%%%%%%%%%%%%%%%%%%%%%%%%%%%%%%%%%%%%%%%%%%%%%%%%%%%

In this section we prove our quantitative Gibbard-Satterthwaite theorem for one voter, Theorem \ref{thm:quant_GS_1voter}. As mentioned before, we present this proof before proving Theorem \ref{thm:k_refined}, because the proof of Theorem \ref{thm:k_refined} follows the lines of this proof, with slight modifications needed to deal with having $n > 1$ coordinates.

For the remainder of this section, let us fix the number of voters to be 1, the number of alternatives $k\geq 3$, and the SCF $f$, which satisfies $\Dist \left( f, \NONMANIP \right) \geq \eps$. Accordingly, we typically omit the dependence of various sets (e.g.\ boundaries between two alternatives) on $f$. 

An additional notational remark: since our SCF is on one voter only, we omit the subscripts that denote the coordinate we are on. E.g.\ we write simply $\Inf^{a,b}$ instead of $\Inf_1^{a,b}$, etc.\

We present the proof in several steps.

%%%%%%%%%%%%%%%%%%%%%%%%%%%%%%%%%%%%%%%%%%%%%%%%%%%%%%%%
\subsection{Large boundary between two alternatives} %%%
%%%%%%%%%%%%%%%%%%%%%%%%%%%%%%%%%%%%%%%%%%%%%%%%%%%%%%%%

The first thing we have to establish is a large boundary between two alternatives. This can be done just like in Lemma~\ref{lem:boundaries2}, except there are two small differences. On the one hand, the assumption of the lemma, namely that $\Dist \left( f, \NONMANIP \right) \geq \eps$, is weaker than that of the original lemma. On the other hand, here we only need one big boundary, unlike in Lemma~\ref{lem:boundaries2}, where Isaksson et al.~\cite{isaksson2010geometry} showed that there are two big boundaries in two different coordinates. The following lemma formulates what we need.

\begin{lemma}\label{lem:big_ab_bdry}
Recall that $f$ is a SCF on 1 voter and $k \geq 3$ alternatives which satisfies $\Dist \left( f, \NONMANIP \right) \geq \eps$. Let $\sigma \in S_k$ be selected uniformly. Then either
\begin{equation}\label{eq:2manip_pt}
\p \left( \sigma \in M_2 \right) \geq \frac{4\eps}{k^6}
\end{equation}
or there exist alternatives $a,b \in \left[ k \right]$, $a\neq b$ such that
\begin{equation}\label{eq:big_ab_bdry}
\Inf^{a,b;\left[a:b\right]} \geq \frac{2\eps}{k^6}.
\end{equation}
\end{lemma}

\begin{proof}
The proof is just like the proof of Lemma \ref{lem:boundaries2}. First, suppose that $\Inf^{a,b;z} \geq \frac{2\eps}{k^6}$ for some pair of alternatives $a \neq b$, and transposition $z \neq \left[a:b \right]$. Then by Lemma \ref{lem:nonManipBoundary}, for any point $\left( \sigma, \sigma' \right) \in B^{a,b;z}$, at least one of $\sigma$ or $\sigma' = z \sigma$ is a 2-manipulation point. Then
\[
\left| M_2 \right| \geq \left| B^{a,b;z} \right| = 2 \cdot k! \cdot \Inf^{a,b;z} \geq \frac{4\eps}{k^6} k!,
\]
and dividing with $k!$ gives \eqref{eq:2manip_pt}. So for the remainder of the proof we may assume that $\Inf^{a,b;z} < \frac{2\eps}{k^6}$ for every $a \neq b$ and $z \neq \left[a:b \right]$.

For every $a \in \left[k\right]$, $\Dist \left( f, \tp_{\left\{a\right\}} \right) \geq \eps$, so $\p \left( f \left( \sigma \right) = a \right) \leq 1 - \eps$. On the other hand, there exists an alternative, say $a \in \left[ k \right]$, such that $\p \left( f\left( \sigma \right) = a \right) \geq \frac{1}{k}$. So for this alternative we have
\[
\Var \left( \mathbf{1} \left[ f\left( \sigma \right) = a \right] \right) \geq \frac{\eps}{k},
\]
and consequently using \cite[Corollary 6.5]{isaksson2010geometry} and \cite[Proposition 2.3]{isaksson2010geometry} we have
\[
\sum_{w \in T} \sum_{b \neq a} \Inf^{a,b;w} = \sum_{w \in T} \Inf^{a;w} \geq \frac{1}{k^2} \Var \left( \mathbf{1} \left[ f\left( \sigma \right) = a \right] \right) \geq \frac{\eps}{k^3}.
\]
Hence there must exist some $w \in T$ and $b \neq a$ such that $\Inf^{a,b;w} \geq \frac{2\eps}{k^6}$, but by our assumption we must have $w = \left[a:b\right]$.
\end{proof}

If \eqref{eq:2manip_pt} holds, then we are done, so in the following we assume that \eqref{eq:big_ab_bdry} holds.

We know that $\sigma$ is on $B^{a,b;\left[a:b\right]}$ if $f\left( \sigma \right) = a$ and $f \left( \left[a:b\right] \sigma \right) = b$. We know that if $b \stackrel{\sigma}{>} a$, then $\sigma$ is a 2-manipulation point, so if this happens in more than half of the cases when $\sigma$ is on $B^{a,b;\left[a:b\right]}$, then we have 
\[
\p \left( \sigma \in M_2 \right) \geq \frac{2\eps}{k^6},
\]
in which case we are again done. So we may assume in the following that
\begin{equation}\label{eq:B_1vot}
\p \left( \sigma \in B \right) \geq \frac{2\eps}{k^6},
\end{equation}
where
\[
B := \left\{ \sigma : f\left( \sigma \right) = a, f\left( \left[a:b \right] \sigma \right) = b, a \stackrel{\sigma}{>} b \right\}.
\]

%%%%%%%%%%%%%%%%%%%%%%%%%%%%%%%%%%%%%%%%%%%%%%%%%%%%%%%%%%
\subsection{Division into cases}\label{sec:cases_1vot} %%%
%%%%%%%%%%%%%%%%%%%%%%%%%%%%%%%%%%%%%%%%%%%%%%%%%%%%%%%%%%

We again divide into two cases. %The difference now is that we only have one coordinate, so we only have one fiber---which is either small or large.

We introduce the set $\bar{F}$ of permutations where $a$ is directly above $b$:
\[
\bar{F} := \left\{ \sigma \in S_k : a \stackrel{\sigma}{>} b \text{ and } b \stackrel{\sigma'}{>} a, \text{ where } \sigma' = \left[a:b\right] \sigma \right\}.
\]
One of the following two cases must hold.

\textbf{Case 1: Small fiber case.} We have
\begin{equation}\label{eq:sm_fbr_case}
\p \left( \sigma \in B \, \middle| \, \sigma \in \bar{F} \right) \leq 1 - \frac{\eps}{4k}.
\end{equation}

\textbf{Case 2: Large fiber case.} We have
\begin{equation}\label{eq:lg_fbr_case}
\p \left( \sigma \in B \, \middle| \, \sigma \in \bar{F} \right) > 1 - \frac{\eps}{4k}.
\end{equation}

%%%%%%%%%%%%%%%%%%%%%%%%%%%%%%%%%%%%%%%%%%%%%%%%%%%%%%%%
\subsection{Small fiber case}\label{sec:sm_fbr_1vot} %%%
%%%%%%%%%%%%%%%%%%%%%%%%%%%%%%%%%%%%%%%%%%%%%%%%%%%%%%%%

In this section we assume that \eqref{eq:sm_fbr_case} holds.

We first formalize that the boundary $\partial \left( B \right)$ of $B$ is big. (Recall the definition of $\partial\left( B \right)$ from Section \ref{sec:bdry_of_bdry}, where we first introduced the graph $G$ on the vertex set $\bar{F}$ where two ranking profiles are adjacent if they differ in a single transposition, but treating the block formed by adjacent alternatives $a$ and $b$ as a single element.) The proof uses the canonical path method, as successfully adapted to this setting by Isaksson, Kindler and Mossel \cite{isaksson2010geometry}.

\begin{lemma}\label{lem:comparable_bdries_1vot}
If \eqref{eq:sm_fbr_case} holds, then
\begin{equation}\label{eq:comp_bdry_of_B}
\p \left( \sigma \in \partial \left( B \right) \right) \geq \frac{\eps}{2 k^4} \p \left( \sigma \in B \right).
\end{equation}
\end{lemma}

\begin{proof}
As mentioned before, the proof is a canonical path argument. Let $B^c = \bar{F} \setminus B$. For every $\left( \sigma, \sigma' \right) \in B \times B^c$, we define a canonical path from $\sigma$ to $\sigma'$, which has to pass through at least one edge in $\partial_e \left( B \right)$. Then if we show that every edge in $\partial_e \left( B \right)$ lies on at most $r$ canonical paths, then it follows that $\left| \partial_e \left( B \right) \right| \geq \left|B \right| \left|B^c \right| / r$.

So let $\left( \sigma, \sigma' \right) \in B \times B^c$. We apply the path construction of \cite[Proposition 6.4.]{isaksson2010geometry}, but considering the block formed by $a$ and $b$ as a single element. Since this path goes from $\sigma$ (which is in $B$) to $\sigma'$ (which is in $B^c$), it must pass through at least one edge in $\partial_e \left( B \right)$.

For a given edge $\left( \pi, \pi' \right) \in \partial_e \left( B \right)$, at most how many possible $\left( \sigma, \sigma' \right) \in B \times B^c$ pairs are there such that the canonical path between $\sigma$ and $\sigma'$ defined above passes through $\left( \pi, \pi' \right)$? We learn from \cite[Proposition 6.4.]{isaksson2010geometry} that there are at most $\left( k - 1\right)^2 \left( k - 1 \right)! / 2 < k^2 \left( k - 1 \right)! / 2$ possibilities for the pair $\left( \sigma, \sigma' \right)$.

Recall that $\left| \bar{F} \right| = \left( k - 1 \right)!$. By our assumption we have $\left| B \right| \leq \left( 1 - \frac{\eps}{4k} \right) \left( k - 1 \right)!$, and so $\left| B^c \right| \geq \frac{\eps}{4k} \left( k - 1 \right)!$. Therefore
\[
\left| \partial_e \left( B \right) \right| \geq \frac{\left|B \right| \left| B^c \right|}{\frac{k^2}{2} \left(k-1\right)!} \geq \frac{\eps}{2 k^3} \left| B \right|.
\]
Now in $G$ every ranking profile has $k - 2 < k$ neighbors, which implies \eqref{eq:comp_bdry_of_B}.
\end{proof}

\begin{corollary}\label{cor:big_bdry_of_B}
If \eqref{eq:sm_fbr_case} holds, then
\begin{equation}\label{eq:big_bdry_of_B}
\p \left( \sigma \in \partial \left( B \right) \right) \geq \frac{\eps^2}{k^{10}}.
\end{equation}
\end{corollary}

\begin{proof}
Combine Lemma \ref{lem:comparable_bdries_1vot} and \eqref{eq:B_1vot}.
\end{proof}

Next we want to find manipulation points on the boundary $\partial \left( B \right)$. The next lemma tells us that if we are on the boundary $\partial \left( B \right)$, then either we can find manipulation points easily, or we are at a local dictator on three alternatives.

\begin{lemma}\label{lem:bdry_of_B}
Suppose that $\sigma \in \partial \left( B \right)$. Then
\begin{itemize}
\item either $\sigma \in \LD \left( a,b \right)$,
\item or there exists $\hat{\sigma} \in M_3$ such that $\hat{\sigma}$ is equal to $\sigma$ or $\left[a:b\right] \sigma$ except that the position of a third alternative $c$ might be shifted arbitrarily.
\end{itemize}
\end{lemma}

\begin{proof}
Since $\sigma \in \partial \left( B \right) \subseteq B$, we know that $f \left( \sigma \right) = a$, and if $\sigma' = \left[a:b\right] \sigma$, then $f\left( \sigma' \right) = b$. Let $\pi \in B^c$ denote the ranking profile such that $\left( \sigma, \pi \right) \in \partial_e \left( B \right)$, and let $\pi' = \left[a:b \right] \pi$. Since $\pi \notin B$, $\left( f\left( \pi \right), f \left( \pi' \right) \right) \neq \left( a, b \right)$. Then, by Lemma~\ref{lem:nonManipBoundary}, if $f\left( \pi \right) \neq f\left( \pi' \right)$, then one of $\pi$ and $\pi'$ is a 2-manipulation point. So assume $f\left( \pi \right) = f \left( \pi' \right)$.

There are two cases to consider: either $\sigma$ and $\pi$ differ by an adjacent transposition not involving the block of $a$ and $b$, or they differ by an adjacent transposition that moves the block of $a$ and $b$.

In the former case, it is not hard to see that one of $\sigma$, $\sigma'$, $\pi$, $\pi'$ is a 2-manipulation point, by Lemma \ref{lem:nonManipBoundary}.

If $\sigma$ and $\pi$ differ by an adjacent transposition that involves the block of $a$ and $b$, then there are again two cases to consider: either this transposition moves the block of $a$ and $b$ up in the ranking, or it moves it down.

If the block of $a$ and $b$ is moved up to get from $\sigma$ to $\pi$, then we must have $f\left( \pi \right) = a$, or else $\sigma$ or $\pi$ is a 3-manipulation point. Then we must have $f \left( \pi' \right) = f\left( \pi \right) = a$, in which case $\pi'$ is a 3-manipulation point, since $f\left( \sigma' \right) = b$.

The final case is when the block of $a$ and $b$ is moved down to get from $\sigma$ to $\pi$, and a third alternative, call it $c$, is moved up, directly above the block of $a$ and $b$. Now if $f\left( \pi \right) = d \notin \left\{a, b, c \right\}$, then $\sigma$ or $\pi$ is a 3-manipulation point. If $f\left( \pi \right) = f \left( \pi' \right) = a$, then $\pi'$ is a 3-manipulation point, whereas if $f\left( \pi \right) = f \left( \pi' \right) = b$, then $\pi$ is a 3-manipulation point. The remaining case is when $f\left( \pi \right) = f \left( \pi' \right) = c$. Now if $f\left( \left[b:c\right] \sigma \right) \neq a$ or $f\left( \left[a:c\right] \sigma' \right) \neq b$, then we again have a 3-manipulation point close to $\sigma$. Otherwise $\sigma \in \LD \left( a, b \right)$.
\end{proof}

The following corollary then tells us that either we have found many 3-manipulation points, or we have many local dictators on three alternatives.

\begin{corollary}\label{cor:manip_or_loc_dict}
If \eqref{eq:sm_fbr_case} holds, then either
\begin{equation}\label{eq:lot_of_loc_dict}
\sum_{c \notin \left\{ a,b \right\}} \p \left( \sigma \in \LD^{\left\{a,b,c\right\}} \right) = \p \left( \sigma \in \LD \left( a, b \right) \right) \geq \frac{\eps^2}{2 k^{10}}
\end{equation}
or
\[
\p \left( \sigma \in M_3 \right) \geq \frac{\eps^2}{4 k^{12}}.
\]
\end{corollary}

%%%%%%%%%%%%%%%%%%%%%%%%%%%%%%%%%%%%%%%%%%%%%%%%%%%%%%%%%%%%%%%%%%%%%%%%%
\subsubsection{Dealing with local dictators}\label{sec:loc_dict_1vot} %%%
%%%%%%%%%%%%%%%%%%%%%%%%%%%%%%%%%%%%%%%%%%%%%%%%%%%%%%%%%%%%%%%%%%%%%%%%%

So the remaining case we have to deal with in this small fiber case is when \eqref{eq:lot_of_loc_dict} holds, i.e.\ we have many local dictators on three alternatives.

\begin{lemma}\label{lem:loc_dict_abc_to_top_1vot}
Suppose $\sigma \in \LD^{\left\{a,b,c\right\}}$ for some alternative $c\notin \left\{a,b\right\}$. Let $\sigma'$ be equal to $\sigma$ except that the block of $a$, $b$ and $c$ is moved to the top of the coordinate. Then
\begin{itemize}
\item either $\sigma' \in \LD^{\left\{a,b,c\right\}}$,
\item or there exists a 3-manipulation point $\hat{\sigma}$ which agrees with $\sigma$ except that the positions of $a$, $b$ and $c$ might be shifted arbitrarily.
\end{itemize}
\end{lemma}

\begin{proof}
W.l.o.g.\ we may assume that in $\sigma$ alternative $a$ is ranked above $b$, which is ranked above $c$. Now bubble up $a$ to the top, using adjacent transpositions. If at any point during this the outcome of $f$ is not $a$, then we have found a 2-manipulation point. Now bubble up $b$ to right below $a$, and then bubble up $c$ to be right below $b$. Again, if at any point during this the outcome of $f$ is not $a$, then there is a 2-manipulation point. Otherwise we now have $a$, $b$ and $c$ at the top (in this order), with the outcome of $f$ being $a$. Now permuting alternatives $a$, $b$ and $c$ at the top, we either have a 3-manipulation point, or $\sigma' \in \LD^{\left\{a,b,c\right\}}$.
\end{proof}

\begin{corollary}\label{cor:many_dict_at_top_1vot}
If \eqref{eq:lot_of_loc_dict} holds, then either
\begin{equation}\label{eq:loc_dict_abc_on_top_1vot}
\sum_{c \notin \left\{ a, b \right\}} \p \left( \sigma \in \LD^{\left\{a,b,c\right\}},  \left\{ \sigma \left( 1 \right), \sigma \left( 2 \right), \sigma \left( 3 \right) \right\} = \left\{a,b,c\right\} \right) \geq \frac{\eps^2}{4 k^{11}}
\end{equation}
or
\[
\p \left( \sigma \in M_3 \right) \geq \frac{\eps^2}{4 k^{13}}.
\]
\end{corollary}

\begin{proof}
Lemma \ref{lem:loc_dict_abc_to_top_1vot} tells us that when we move the block of $a$, $b$, and $c$ up to the top, we either encounter a 3-manipulation point, or we get a local dictator on $\left\{ a, b, c \right\}$ at the top.

If we get a 3-manipulation point, by the describtion of this manipulation point in the lemma, there can be at most $k^3$ ranking profiles that give the same manipulation point.

If we arrive at a local dictator at the top, then there could have been at most $k$ different places where the block of $a$, $b$ and $c$ could have come from.
\end{proof}

Now \eqref{eq:loc_dict_abc_on_top_1vot} is equivalent to
\begin{equation}\label{eq:loc_dict_abc_on_top2_1vot}
\sum_{c \notin \left\{ a, b \right\}} \p \left( \sigma \in \LD^{\left\{a,b,c\right\}}, \left( \sigma \left( 1 \right), \sigma \left( 2 \right), \sigma \left( 3 \right) \right) = \left(a,b,c\right) \right) \geq \frac{\eps^2}{24 k^{11}}.
\end{equation}
We know that
\[
\p \left( \left( \sigma \left( 1 \right), \sigma \left( 2 \right), \sigma \left( 3 \right) \right) = \left(a,b,c\right) \right) = \frac{1}{k \left( k-1 \right) \left( k-2 \right)} \leq \frac{6}{k^3},
\]
and so \eqref{eq:loc_dict_abc_on_top2_1vot} implies (recall Definition \ref{def:cond})
\begin{equation}\label{eq:loc_dict_abc_on_top3_1vot}
\sum_{c \notin \left\{ a, b \right\}} \p^{\left(a,b,c\right)} \left( \sigma \in \LD^{\left\{a,b,c\right\}}  \right) \geq \frac{\eps^2}{144 k^{8}}.
\end{equation}

Now fix an alternative $c \notin \left\{a,b\right\}$ and define the graph $G_{\left( a, b, c \right)} = \left( V_{\left( a, b, c \right)}, E_{\left( a, b, c \right)} \right)$ to have vertex set
\[
V_{\left( a, b, c \right)} := \left\{ \sigma \in S_k : \left( \sigma \left( 1 \right), \sigma \left( 2 \right), \sigma \left( 3 \right) \right) = \left(a,b,c\right) \right\}
\]
and for $\sigma, \pi \in V_{\left( a, b, c \right)}$ let $\left( \sigma, \pi \right) \in E_{\left( a, b, c \right)}$ if and only if $\sigma$ and $\pi$ differ by an adjacent transposition. So $G_{\left( a, b, c \right)}$ is the subgraph of the refined rankings graph induced by the vertex set $V_{\left( a, b, c \right)}$. (If $k = 3$ or $k=4$, then this graph consists of only one vertex, and no edges.)

Let
\[
T \left( a, b, c \right) := V_{\left( a, b, c \right)} \cap \LD^{\left\{a,b,c\right\}},
\]
and let $\partial_e \left( T \left( a, b, c \right) \right)$ and $\partial \left( T \left( a, b, c \right) \right)$ denote the edge- and vertex-boundary of $T \left( a, b, c \right)$ in $G_{\left( a, b, c \right)}$, respectively.

The next lemma shows that unless $T \left( a, b, c \right)$ is almost all of $V_{\left( a, b, c \right)}$, the size of the boundary $\partial \left( T \left( a, b, c \right) \right)$ is comparable to the size of $T \left( a, b, c \right)$. The proof uses a canonical path argument, just like in Lemma \ref{lem:comparable_bdries_1vot}.

\begin{lemma}\label{lem:lg_bdry_for_loc_dict_1vot}
Let $c \notin \left\{a,b \right\}$ be arbitrary. Write $T \equiv T \left( a, b, c \right)$ for simplicity. If $\p^{\left(a,b,c\right)} \left( \sigma \in T \right) \leq 1 - \delta$, then
\begin{equation}\label{eq:lg_bdry_for_loc_dict_1vot}
\p^{\left( a, b, c \right)} \left( \sigma \in \partial \left( T \right) \right) \geq \frac{\delta}{k^3} \p^{\left( a, b, c \right)} \left( \sigma \in T \right).
\end{equation}
\end{lemma}

\begin{proof}
Let $T^c = V_{\left( a, b, c \right)} \setminus T \left( a, b, c \right)$. For every $\left( \sigma, \sigma' \right) \in T \times T^c$, we define a canonical path from $\sigma$ to $\sigma'$, which has to pass through at least one edge in $\partial_e \left( T \right)$. Then if we show that every edge in $\partial_e \left( T \right)$ lies on at most $r$ canonical paths, then it follows that $\left| \partial_e \left( T \right) \right| \geq \left| T \right| \left| T^c \right| / r$.

So let $\left( \sigma, \sigma' \right) \in T \times T^c$. We apply the path construction of \cite[Proposition 6.4.]{isaksson2010geometry}, but only to alternatives $\left[k \right] \setminus \left\{ a, b, c \right\}$.

The analysis of this construction is done in exactly the same way as in Lemma \ref{lem:comparable_bdries_1vot}; in the end we get that there are at most $k^2 \left( k - 3 \right)!$ paths that pass through a given edge in $\partial_e \left( T \right)$.

Recall that $\left| V_{\left(a,b,c\right)} \right| = \left( k - 3 \right)!$ and that by our assumption $\left| T \right| \leq \left( 1 - \delta \right) \left( k - 3 \right)!$, so $\left| T^c \right| \geq \delta \left( k - 3 \right)!$. Therefore
\[
\left| \partial_e \left( T \right) \right| \geq \frac{\left| T \right| \left| T^c \right|}{k^2  \left( k - 3 \right)!} \geq \frac{\delta}{k^2} \left| T \right|.
\]
Now every vertex in $V_{\left(a,b,c\right)}$ has $k - 4 < k$ neighbors, which implies \eqref{eq:lg_bdry_for_loc_dict_1vot}.
\end{proof}

The next lemma tells us that if $\sigma$ is on the boundary of a set of local dictators on $\left\{a,b,c \right\}$ for some alternative $c \notin \left\{a,b\right\}$, then there is a 2-manipulation point $\hat{\sigma}$ which is close to $\sigma$.

\begin{lemma}\label{lem:manip_pts_on_bdry_of_loc_dict_1vot}
Suppose $\sigma \in \partial \left( T \left( a,b, c \right) \right)$ for some $c\notin \left\{a,b\right\}$. Then there exists $\hat{\sigma} \in M_2$ which equals $z \sigma$ for some adjacent transposition $z$ that does not involve $a$, $b$ or $c$, except that the order of the block of $a$, $b$ and $c$ might be rearranged.
\end{lemma}

\begin{proof}
Let $\pi$ be the ranking profile such that $\left( \sigma, \pi \right) \in \partial_e \left( T \left( a, b , c \right) \right)$, and let $z$ be the adjacent transposition in which they differ, i.e.\ $\pi = z \sigma$. Since $\pi \notin T \left( a, b, c \right)$, there exists a reordering of the block of $a$, $b$, and $c$ at the top of $\pi$ such that the outcome of $f$ is not the top ranked alternative. Call the resulting vector $\pi'$. W.l.o.g.\ let us assume that $\pi' \left( 1 \right) = a$. Let us also define $\sigma' := z \pi'$. Now $\pi'$ is a 2-manipulation point, since $f\left( \sigma' \right) = a$.
\end{proof}

The next corollary puts together Corollary \ref{cor:many_dict_at_top_1vot} and Lemmas \ref{lem:lg_bdry_for_loc_dict_1vot} and \ref{lem:manip_pts_on_bdry_of_loc_dict_1vot}.

\begin{corollary}\label{cor:manip_by_bdry_loc_dict_1vot}
Suppose \eqref{eq:loc_dict_abc_on_top_1vot} holds. Then if for every $c \notin \left\{a,b\right\}$ we have $\p^{\left(a,b,c\right)} \left( \sigma \in T \left( a,b, c \right) \right) \leq 1 - \frac{\eps}{100 k}$, then
\[
\p \left( \sigma \in M_2 \right) \geq \frac{\eps^3}{10^5 k^{16}}.
\]
\end{corollary}

\begin{proof}
We know that \eqref{eq:loc_dict_abc_on_top_1vot} implies
\[
\sum_{c \notin \left\{a,b \right\}} \p^{a,b,c} \left( \sigma \in T \left( a,b,c\right) \right) \geq \frac{\eps^2}{144 k^8}.
\]
Now using the assumptions, Lemma \ref{lem:lg_bdry_for_loc_dict_1vot} with $\delta = \frac{\eps}{100k}$, and Lemma \ref{lem:manip_pts_on_bdry_of_loc_dict_1vot}, we have
\begin{align*}
\p \left( \sigma \in M_2 \right) &\geq \sum_{c \neq \left\{a,b \right\}} \frac{1}{k^3} \p^{\left(a,b,c \right)} \left( \sigma \in M_2 \right) \geq \sum_{c\notin \left\{a,b\right\}} \frac{1}{6 k^4} \p^{\left( a, b, c \right)} \left( \sigma \in \partial \left( T \left( a,b,c \right) \right) \right)\\
&\geq \sum_{c \notin \left\{a,b\right\}} \frac{\eps}{600 k^{8}} \p^{\left( a, b , c \right)} \left( \sigma \in T \left( a, b, c \right) \right)\geq \frac{\eps^3}{86400 k^{16}} \geq \frac{\eps^3}{10^5 k^{16}}. \qedhere
\end{align*}
\end{proof}

So again we are left with one case to deal with: if there exists an alternative $c \notin \left\{ a, b \right\}$ such that $\p^{\left(a,b,c\right)} \left( \sigma \in T \left( a,b, c \right) \right) > 1 - \frac{\eps}{100 k}$. Define a subset of alternatives $K \subseteq \left[k \right]$ in the following way:
\[
K := \left\{a,b \right\} \cup \left\{ c \in \left[k \right] \setminus \left\{a,b \right\} : \p^{\left(a,b,c\right)} \left( \sigma \in T \left( a,b, c \right) \right) > 1 - \frac{\eps}{100 k} \right\}.
\]
In addition to $a$ and $b$, $K$ contains those alternatives that whenever they are at the top with $a$ and $b$, they form a local dictator with high probability.

So our assumption now is that $\left| K \right| \geq 3$.

Our next step is to show that unless we have many manipulation points, for any alternative $c \in K$, conditioned on $c$ being at the top, the outcome of $f$ is $c$ with probability close to 1.

\begin{lemma}\label{lem:cond_on_top_1vot}
Let $c \in K$. Then either
\begin{equation}\label{eq:cond_on_top1_1vot}
\p^{\left( c \right)} \left( f \left( \sigma \right) = c \right) \geq 1 - \frac{\eps}{50 k},
\end{equation}
or
\begin{equation}\label{eq:else_2manip_1vot}
\p \left( \sigma \in M_2 \right) \geq \frac{\eps}{100 k^4}.
\end{equation}
\end{lemma}

\begin{proof}
First assume that $c \notin \left\{ a, b \right\}$.

Let $\sigma$ be uniform according to $\p^{\left( c \right)}$, i.e.\ uniform on $S_k$ conditioned on $\sigma \left( 1 \right) = c$. Define $\sigma'$, where $\sigma'$ is constructed from $\sigma$ by first bubbling up alternative $a$ to just below $c$, using adjacent transpositions, and then bubbling up $b$ to just below $a$. Clearly $\sigma'$ is distributed according to $\p^{\left(c,a,b\right)}$, i.e.\ it is uniform on $S_k$ conditioned on $\left( \sigma \left( 1 \right), \sigma \left( 2 \right), \sigma \left( 3 \right) \right) = \left( c, a, b \right)$.

Since $c \in K$, we know that $\p^{\left(c,a,b \right)} \left( \sigma \in  \LD^{\left\{a,b,c\right\}} \right) > 1 - \frac{\eps}{100k}$. This also means that
\[
\p^{\left(c \right)} \left( \sigma' \in \LD^{\left\{a,b,c\right\}} \right) > 1 - \frac{\eps}{100k}.
\]
Now we can partition the ranking profiles into three parts, based on the outcome of the SCF $f$ at $\sigma$ and $\sigma'$:
\begin{align*}
I_1 &= \left\{ \sigma : f\left( \sigma \right) = c, f \left( \sigma' \right) = c \right\}\\
I_2 &= \left\{ \sigma : f\left( \sigma \right) \neq c, f \left( \sigma' \right) = c \right\}\\
I_3 &= \left\{ \sigma : f \left( \sigma' \right) \neq c \right\}.
\end{align*}
If $\p^{\left( c \right)} \left( I_1 \right) \geq 1 - \frac{\eps}{50k}$, then \eqref{eq:cond_on_top1_1vot} holds. Otherwise we have $\p^{\left(c\right)} \left( I_2 \cup I_3 \right) \geq \frac{\eps}{50k}$, and since $\p^{\left( c \right)} \left( I_3 \right) \leq \frac{\eps}{100 k}$, we have $\p^{\left( c \right)} \left( I_2 \right) \geq \frac{\eps}{100k}$.

Now if $\sigma \in I_2$, then we know that there is a 2-manipulation point along the way as we go from $\sigma$ to $\sigma'$. I.e.\ to every $\sigma \in I_2$ there exists $\hat{\sigma} \in M_2$ such that $\hat{\sigma}$ is equal to $\sigma$ except perhaps $a$ and $b$ are shifted arbitrarily. So there can be at most $k^2$ ranking profiles $\sigma$ giving the same 2-manipulation point $\hat{\sigma}$, and so we have
\[
\p\left( \sigma \in M_2 \right) \geq \frac{1}{k} \p^{\left(c\right)} \left( \sigma \in M_2 \right) \geq \frac{1}{k^3} \p^{\left( c \right)} \left( I_2 \right) \geq \frac{\eps}{100 k^4},
\]
showing \eqref{eq:else_2manip_1vot}.

Now suppose $c \in \left\{a,b\right\}$, w.l.o.g.\ assume $c = a$. We know that $\left|K \right| \geq 3$ and so there exists an alternative $d\in K \setminus \left\{a,b \right\}$. We can then do the same thing as above, but we now bubble up $b$ and $d$.
\end{proof}

We now deal with alternatives that are not in $K$: either we have many manipulation points, or for any alternative $d \notin K$, the outcome of $f$ is \emph{not} $d$ with probability close to 1.

\begin{lemma}\label{lem:d_notin_K_1vot}
Let $d \notin K$. If $\p \left( f \left( \sigma \right) = d \right) \geq \frac{\eps}{4k}$, then
\[
\p \left( \sigma \in M_2 \right) \geq \frac{\eps^2}{10^6 k^{9}}.
\]
\end{lemma}

\begin{proof}
Let $\sigma$ be such that $f \left( \sigma \right) = d$. Bubble up $d$ to the top, and call this ranking profile $\sigma'$. Now if $f\left( \sigma' \right) \neq d$, then we know that there exists a 2-manipulation point $\hat{\sigma}$ along the way, i.e.\ a 2-manipulation $\hat{\sigma}$ which agrees with $\sigma$ except perhaps $d$ is shifted arbitrarily. Consequently, either
\[
\p \left( \sigma \in M_2 \right) \geq \frac{\eps}{8k^2},
\]
in which case we are done, or 
\[
\p \left( \sigma: f\left( \sigma \right) = f\left( \sigma' \right) = d \right) \geq \frac{\eps}{8k}.
\]

Next, let us bubble up $a$ to just below $d$, and then bubble up $b$ to just below $d$. Denote this ranking profile by $\sigma^{\left(d,b,a\right)}$, and analogously define $\sigma^{\left(d,a,b\right)}, \sigma^{\left(a,b,d\right)}, \sigma^{\left(a,d,b\right)}, \sigma^{\left(b,a,d\right)}$, and $\sigma^{\left(b,d,a\right)}$. Either we encounter a 2-manipulation point $\hat{\sigma}$ along the way of bubbling up to $\sigma^{\left(d,b,a\right)}$ ($\hat{\sigma}$ agrees with $\sigma$ except $d$ is at the top, and $a$ and $b$ might be arbitrarily shifted), or the outcome of the SCF $f$ is $d$ all along. So we have that either
\[
\p \left( \sigma \in M_2 \right) \geq \frac{\eps}{16 k^3},
\]
in which case we are done, or
\[
\p \left( \sigma : f\left( \sigma \right) = f\left( \sigma' \right) = f\left( \sigma^{\left(d,b,a\right)} \right) = f\left( \sigma^{\left(d,a,b\right)} \right) = d \right) \geq \frac{\eps}{16 k}.
\]

Now start from $\sigma^{\left(d,a,b\right)}$. First swap $a$ and $d$ to get $\sigma^{\left(a,d,b\right)}$, then swap $d$ and $b$ to get $\sigma^{\left(a,b,d\right)}$, and finally bubble $d$ and $b$ down to their original positions in $\sigma$, except for the fact that $a$ is now at the top of the coordinate. Call this profile $\bar{\sigma}$. Since $\sigma$ is uniformly distributed, $\bar{\sigma}$ is distributed according to $\p_1^{\left( a \right)}$, i.e.\ uniformly conditional on $\bar{\sigma} \left( 1 \right) = a$. Now note that one of the following three events has to happen. (These events are not mutually exclusive.)
\begin{align*}
I_1 &= \left\{ f\left( \sigma^{\left(a,d,b\right)} \right) = f\left( \sigma^{\left(a,b,d\right)} \right) = a \right\}\\
I_2 &= \left\{ f \left( \bar{\sigma} \right) \neq a \right\}\\
I_3 &= \{ \sigma: \exists\ \hat{\sigma} \in M_2 \text{ which is equal to } \sigma \text{ except } a \text{ is shifted}\\
&\qquad \ \text{ to the top, and }b \text{ and } d \text{ may be shifted arbitrarily}\}.
\end{align*}
Since $a\in K$, we know by Lemma \ref{lem:cond_on_top_1vot} that (unless we already have enough manipulation points by the lemma) we must have
\[
\p \left( f\left( \bar{\sigma} \right) \neq a \right) = \p^{\left( a \right)} \left( f\left( \bar{\sigma} \right) \neq a \right) \leq \frac{\eps}{50k}.
\]
Consequently 
\[
\p \left( I_1 \cup I_3, f\left( \sigma \right) = f\left( \sigma' \right) = f\left( \sigma^{\left(d,b,a\right)} \right) = f\left( \sigma^{\left(d,a,b\right)} \right) = d  \right) \geq \frac{\eps}{16 k} - \frac{\eps}{50 k} = \frac{17 \eps}{400 k},
\]
and so either
\[
\p \left( \sigma \in M_2 \right) \geq \frac{17 \eps}{800 k^3},
\]
in which case we are done, or
\[
\p \left( \sigma: f\left( \sigma^{\left(d,b,a\right)} \right) = f\left( \sigma^{\left(d,a,b\right)} \right) = d, f\left( \sigma^{\left(a,b,d\right)} \right) = f\left( \sigma^{\left(a,d,b\right)} \right) = a \right) \geq \frac{17 \eps}{800 k}.
\]

Next, we can do the same thing with $b$ on top, and we ultimately get that either
\[
\p \left( \sigma \in M_2 \right) \geq \frac{\eps}{1600 k^3},
\]
in which case we are done, or
\begin{equation}\label{eq:loc_dict_abd_1vot}
\p^{\left( a, b, d \right)} \left( \sigma^{\left(a,b,d\right)} \in \LD^{\left\{a,b,d\right\}} \right) = \p \left( \sigma: \sigma^{\left( a, b, d \right)} \in \LD^{\left\{a,b,d\right\}} \right) \geq \frac{\eps}{1600 k}.
\end{equation}
Define $G_{\left(a,b,d\right)}$ and $T_{\left( a, b, d \right)}$ analogously to $G_{\left( a, b, c \right)}$ and $T_{\left( a, b, c \right)}$, respectively.

Suppose that \eqref{eq:loc_dict_abd_1vot} holds. We also know that $d\notin K$, so Lemma \ref{lem:lg_bdry_for_loc_dict_1vot} applies, and then Lemma \ref{lem:manip_pts_on_bdry_of_loc_dict_1vot} shows us how to find manipulation points. We can put these arguments together, just like in the proof of Corollary \ref{cor:manip_by_bdry_loc_dict_1vot}, to show what we need:
\begin{align*}
\p \left( \sigma \in M_2 \right) &\geq \frac{1}{k^3} \p^{\left(a,b,d \right)} \left( \sigma \in M_2 \right) \geq \frac{1}{6 k^4} \p^{\left( a, b, d \right)} \left( \sigma \in \partial \left( T \left( a,b,d \right) \right) \right)\\
&\geq \frac{\eps}{600 k^{8}} \p^{\left( a, b , d \right)} \left( \sigma \in T \left( a, b, d \right) \right) \geq \frac{\eps^2}{10^6 k^{9}}. \qedhere
\end{align*}
\end{proof}

Putting together the results of the previous lemmas, there is only one case to be covered, which is covered by the following final lemma. Basically, this lemma says that unless there are enough manipulation points, our function is close to a dictator on the subset of alternatives $K$.
\begin{lemma}\label{lem:final_loc_dict_1vot}
Recall that we assume that $\Dist \left( f, \NONMANIP \right) \geq \eps$. Furthermore assume that $\left| K \right| \geq 3$, for every $c \in K$ we have
\begin{equation}\label{eq:dict_cond_on_top_1vot}
\p^{\left( c \right)} \left( f \left( \sigma \right) = c \right) \geq 1 - \frac{\eps}{50k},
\end{equation}
and for every $d \notin K$ we have
\[
\p \left( f\left( \sigma \right) = d \right) \leq \frac{\eps}{4k}.
\]
Then
\begin{equation}\label{eq:final_manip_1vot}
\p \left( \sigma \in M_2 \right) \geq \frac{\eps}{4k^2}.
\end{equation}
\end{lemma}

\begin{proof}
First note that
\[
\p \left( f \left( \sigma \right) \neq \tp_K \left( \sigma \right) \right) \leq \p \left( f \left( \sigma \right) \notin K \right) + \p \left( f\left( \sigma \right) \neq \tp_K \left( \sigma \right), f\left( \sigma \right) \in K \right).
\]
We know that
\[
\eps \leq \Dist \left( f, \NONMANIP \right) \leq \p \left( f \left( \sigma \right) \neq \tp_K \left( \sigma \right) \right)
\]
and also that
\[
\p\left( f\left( \sigma \right) \notin K \right) \leq \frac{\left| K \right| \eps}{4k} \leq \frac{\eps}{2},
\]
which together imply that necessarily
\[
\p \left( f\left( \sigma \right) \neq \tp_K \left( \sigma \right), f\left( \sigma \right) \in K \right) \geq \frac{\eps}{2}.
\]
Let $\sigma$ be such that $ f\left( \sigma \right) \neq \tp_K \left( \sigma \right)$ and $f\left( \sigma \right) \in K$. Now bubble $\tp_K \left( \sigma \right)$ up to the top in $\sigma$, call this ranking profile $\bar{\sigma}$. Clearly then $\tp_K \left( \bar{\sigma} \right) = \tp_K \left( \sigma \right)$.

There are two cases to consider. If $f\left( \sigma \right) \neq f\left( \bar{\sigma} \right)$, then there is a 2-manipulation point along the way from $\sigma$ to $\bar{\sigma}$, i.e.\ a 2-manipulation point $\hat{\sigma}$ such that $\hat{\sigma}$ agrees with $\sigma$ except perhaps some alternative $c$ is arbitrarily shifted. Otherwise $f \left( \sigma \right) = f\left( \bar{\sigma} \right)$, and so $f \left( \bar{ \sigma} \right) \neq \tp_K \left( \bar{\sigma} \right)$.

Consequently we have that either \eqref{eq:final_manip_1vot} holds, or that
\begin{equation}\label{eq:last_eq_1vot}
\p \left( \sigma : f \left( \bar{ \sigma} \right) \neq \tp_K \left( \bar{\sigma} \right) \right) \geq \frac{\eps}{4}.
\end{equation}
By the construction of $\bar{\sigma}$, we know that $\bar{\sigma}$ is uniformly distributed conditional on $\bar{\sigma} \left( 1 \right) \in K$. Consequently, by \eqref{eq:dict_cond_on_top_1vot}, we have that
\[
\p \left( \sigma : f \left( \bar{ \sigma} \right) \neq \tp_K \left( \bar{\sigma} \right) \right) \leq \frac{\eps}{50 k},
\]
which contradicts with \eqref{eq:last_eq_1vot} since $\frac{\eps}{50k} < \frac{\eps}{4}$.
\end{proof}
This concludes the proof of the small fiber case.

%%%%%%%%%%%%%%%%%%%%%%%%%%%%%%%%%%%%%%%%%%%%%%%%%%%%%%%%%%%%%
\subsection{Large fiber case}\label{sec:lg_fbr_case_1vot} %%%
%%%%%%%%%%%%%%%%%%%%%%%%%%%%%%%%%%%%%%%%%%%%%%%%%%%%%%%%%%%%%

In this section we assume that \eqref{eq:lg_fbr_case} holds. We show that we either have a lot 2-manipulation points or we have a lot of local dictators on three alternatives.

Our first step towards this is the following lemma.

\begin{lemma}\label{lem:cond_ab_top_still_lg_fbr_1vot}
Suppose \eqref{eq:lg_fbr_case} holds. Then
\begin{equation}\label{eq:cond_ab_top_still_lg_fbr_1vot}
\p^{\left(a,b\right)} \left( \sigma \in B \right) \geq 1 - \frac{\eps}{4}.
\end{equation}
\end{lemma}

\begin{proof}
Let $B^c = \bar{F} \setminus B$. Our assumption \eqref{eq:lg_fbr_case} implies that $\p \left( \sigma \in B^c \, \middle| \, \sigma \in \bar{F} \right) \leq \frac{\eps}{4k}$, which means that $\left| B^c \right| \leq \frac{\eps \left( k - 1 \right)!}{4k}$, and so
\[
\p^{\left( a, b \right)} \left( \sigma \notin B \right) \leq \frac{\eps \left( k - 1 \right)!}{ 4k \left( k - 2 \right)!} < \frac{\eps}{4},
\]
which is equivalent to \eqref{eq:cond_ab_top_still_lg_fbr_1vot}.
\end{proof}

The next lemma (together with Section \ref{sec:loc_dict_1vot}) concludes the proof in the large fiber case.

\begin{lemma}\label{lem:lg_fbr_final_1vot}
Suppose \eqref{eq:lg_fbr_case} holds and recall that our SCF $f$ satisfies $\Dist \left( f, \NONMANIP \right) \geq \eps$. Then either
\begin{equation}\label{eq:2manip_again_1vot}
\p \left( \sigma \in M_2 \right) \geq \frac{\eps}{4k^2}
\end{equation}
or
\begin{equation}\label{eq:loc_dict_again_1vot}
\p \left( \sigma \in \LD \left(a,b\right) \right) \geq \frac{\eps}{4k^2}.
\end{equation}
\end{lemma}

\begin{proof}
By Lemma \ref{lem:cond_ab_top_still_lg_fbr_1vot} we know that \eqref{eq:cond_ab_top_still_lg_fbr_1vot} holds.

Let $\sigma \in S_k$ be uniform. Define $\sigma'$ by being the same as $\sigma$ except alternatives $a$ and $b$ are moved to the top of the coordinate: $\sigma' \left( 1 \right) = a$ and $\sigma' \left( 2 \right) = b$. Clearly $\sigma'$ is distributed according to $\p^{\left(a,b\right)} \left( \cdot \right)$. Also define $\sigma'' = \left[a:b \right] \sigma'$.

We partition the set of ranking profiles $S_k$ into three parts:
\begin{align*}
I_1 &:= \left\{ \sigma \in S_k : f\left( \sigma \right) = \tp_{\left\{a,b \right\}} \left( \sigma \right), \left( f \left( \sigma' \right), f \left( \sigma'' \right) \right) = \left( a, b \right) \right\}\\
I_2 &:= \left\{ \sigma \in S_k : f\left( \sigma \right) \neq \tp_{\left\{a,b \right\}} \left( \sigma \right), \left( f \left( \sigma' \right), f \left( \sigma'' \right) \right) = \left( a, b \right) \right\}\\
I_3 &:= \left\{ \sigma \in S_k : \left( f \left( \sigma' \right), f \left( \sigma'' \right) \right) \neq \left( a, b \right) \right\}.
\end{align*}

By \eqref{eq:cond_ab_top_still_lg_fbr_1vot} we know that $\p \left( \sigma \in I_3 \right) \leq \frac{\eps}{4}$. We also know that $\p \left( \sigma \in I_1 \right) \leq 1 - \eps$, since $\Dist \left( f, \NONMANIP \right) \geq \eps$. Therefore we must have
\[
\p \left( \sigma \in I_2 \right) \geq \frac{3\eps}{4} > \frac{\eps}{2}.
\]

Let us partition $I_2$ further, and write it as $I_2 = I_2' \cup \left( \cup_{c \notin \left\{a,b \right\}} I_{2,c} \right)$, where 
\[
I_2' := \left\{ \sigma \in I_2 : f \left( \sigma \right) \neq \tp_{\left\{a,b\right\}} \left( \sigma \right), f\left( \sigma \right) \in \left\{a,b\right\} \right\}
\]
and for any $c \notin \left\{a,b \right\}$,
\[
I_{2,c} := \left\{ \sigma \in I_2 : f \left( \sigma \right) = c \right\}.
\]

Suppose $\sigma \in I_2'$. W.l.o.g.\ let us assume that $a$ is ranked higher than $b$ by $\sigma$, and therefore $f\left( \sigma \right) = b$, since $\sigma \in I_2'$. Then we can get from $\sigma$ to $\sigma'$ by first bubbling up $a$ to the top, and then bubbling up $b$ to just below $a$. Since $f\left( \sigma \right) = b$ and $f\left( \sigma' \right) = a$, there must be a 2-manipulation point $\hat{\sigma}$ along the way, which is equal to $\sigma$ except perhaps the positions of $a$ and $b$ are arbitrarily shifted.

Now suppose that $\sigma \in I_{2,c}$ for some $c \notin \left\{a,b \right\}$. We distinguish two cases: either $c$ is ranked above both $a$ and $b$ in $\sigma$, or it is not.

If not, then say $a$ is ranked above $c$ in $\sigma$. Bubble $a$ all the way to the top, and then bubble $b$ as well, all the way to the top, just below $a$. Since $f\left( \sigma \right) = c$ and $f\left( \sigma' \right) = a$, there must be a 2-manipulation point $\hat{\sigma}$ along the way, which is equal to $\sigma$ except perhaps the positions of $a$ and $b$ are arbitrarily shifted.

If $c$ is ranked above both $a$ and $b$ in $\sigma$, then the argument is similar. First bubble up $a$ and $b$ to just below $c$, and denote this ranking profile by $\tilde{\sigma}$, then permute these three alternatives arbitrarily, and then bubble $a$ and $b$ to the top. It is not hard to think through that either there is a 2-manipulation $\hat{\sigma}$ along the way, which is then equal to $\sigma$ except perhaps the positions of $a$ and $b$ are arbitrarily shifted, or else $\tilde{\sigma} \in \LD^{\left\{a,b,c\right\}}$.

Combining these cases we see that either \eqref{eq:2manip_again_1vot} or \eqref{eq:loc_dict_again_1vot} must hold.
\end{proof}

So if \eqref{eq:2manip_again_1vot} holds then we are done, and if \eqref{eq:loc_dict_again_1vot} holds, then we refer back to Section \ref{sec:loc_dict_1vot}, where we deal with the case of local dictators on three alternatives.

%%%%%%%%%%%%%%%%%%%%%%%%%%%%%%%%%%%%%%%%%%%%%%%%%%%%%%%%%%%%%%%%%%%%%%%%%
\subsection{Proof of Theorem \ref{thm:quant_GS_1voter} concluded} %%%
%%%%%%%%%%%%%%%%%%%%%%%%%%%%%%%%%%%%%%%%%%%%%%%%%%%%%%%%%%%%%%%%%%%%%%%%%

\begin{proof}[Proof of Theorem \ref{thm:quant_GS_1voter}]
Our starting point is Lemma~\ref{lem:big_ab_bdry}, which implies that \eqref{eq:B_1vot} holds (unless we already have many 2-manipulation points). We then consider two cases, as indicated in Section \ref{sec:cases_1vot}.

We deal with the small fiber case---when \eqref{eq:sm_fbr_case} holds---in Section \ref{sec:sm_fbr_1vot}. First, Lemma \ref{lem:comparable_bdries_1vot}, Corollary \ref{cor:big_bdry_of_B}, Lemma \ref{lem:bdry_of_B} and Corollary \ref{cor:manip_or_loc_dict} show that either there are many 3-manipulation points, or there are many local dictators on three alternatives. We then deal with the case of many local dictators in Section \ref{sec:loc_dict_1vot}. Lemma \ref{lem:loc_dict_abc_to_top_1vot}, Corollary \ref{cor:many_dict_at_top_1vot}, Lemmas \ref{lem:lg_bdry_for_loc_dict_1vot} and \ref{lem:manip_pts_on_bdry_of_loc_dict_1vot}, Corollary \ref{cor:manip_by_bdry_loc_dict_1vot}, and Lemmas \ref{lem:cond_on_top_1vot}, \ref{lem:d_notin_K_1vot} and \ref{lem:final_loc_dict_1vot} together show that there are many 3-manipulation points if there are many local dictators on three alternatives, and the SCF is $\eps$-far from the family of nonmanipulable functions.

We deal with the large fiber case---when \eqref{eq:lg_fbr_case} holds---in Section \ref{sec:lg_fbr_case_1vot}. Here Lemma \ref{lem:lg_fbr_final_1vot} shows that either we have many 2-manipulation points, or we have many local dictators on three alternatives. In this latter case we refer back to Section \ref{sec:loc_dict_1vot} to conclude the proof.
\end{proof}

%%%%%%%%%%%%%%%%%%%%%%%%%%%%%%%%%%%%%%%%%%%%%%%%%%%%%%%%%%%%%%%%%%%%%%%%%%%%%%%%%%%%%%%%%%%%%%%%%%%
\section{Inverse polynomial manipulability for any number of alternatives}\label{sec:manip_ref} %%%
%%%%%%%%%%%%%%%%%%%%%%%%%%%%%%%%%%%%%%%%%%%%%%%%%%%%%%%%%%%%%%%%%%%%%%%%%%%%%%%%%%%%%%%%%%%%%%%%%%%

In this section we prove the theorem below, which is the same as our main theorem, Theorem~\ref{cor:k_refined_truenonmanip}, except that the condition of $\Dist \left( f, \NONMANIP \right) \geq \eps$ from  Theorem~\ref{cor:k_refined_truenonmanip} is replaced with the stronger condition $\Dist \left( f, \overline{\NONMANIP} \right) \geq \eps$.

\begin{theorem}\label{thm:k_refined}
Suppose we have $n \geq 2$ voters, $k \geq 3$ alternatives, and a SCF $f : S_k^n \to \left[k\right]$ satisfying $\Dist \left( f, \overline{\NONMANIP} \right) \geq \eps$. Then
\begin{equation}\label{eq:manip_refined2}
\p\left( \sigma \in M \left( f \right) \right)\geq \p \left( \sigma \in M_4 \left( f \right) \right)\geq p \left( \eps, \frac{1}{n}, \frac{1}{k} \right),
\end{equation}
for some polynomial $p$, where $\sigma \in S_k^n$ is selected uniformly. In particular, we show a lower bound of $\frac{\eps^5}{10^9 n^7 k^{46}}$.

An immediate consequence is that
\[
\p \left( \left( \sigma, \sigma' \right) \text{ is a manipulation pair for } f \right) \geq q \left( \eps, \frac{1}{n}, \frac{1}{k} \right),
\]
for some polynomial $q$, where $\sigma \in S_k^n$ is uniformly selected, and $\sigma'$ is obtained from $\sigma$ by uniformly selecting a coordinate $i \in \left\{1, \dots, n \right\}$, uniformly selecting $j \in \left\{1, \dots, n-3 \right\}$, and then uniformly randomly permuting the following four adjacent alternatives in $\sigma_i$: $\sigma_i \left( j \right), \sigma_i \left( j + 1 \right), \sigma_i \left( j + 2 \right)$, and $\sigma_i \left( j + 3 \right)$. In particular, the specific lower bound for $\p \left( \sigma \in M_4 \left( f \right) \right)$ implies that we can take $q \left( \eps, \frac{1}{n}, \frac{1}{k} \right) = \frac{\eps^5}{10^{11} n^8 k^{47}}$.
\end{theorem}

For the remainder of the section, let us fix the number of voters $n \geq 2$, the number of alternatives $k \geq 3$, and the SCF $f$, which satisfies $\Dist \left( f, \overline{\NONMANIP} \right) \geq \eps$. Accordingly, we typically omit the dependence of various sets (e.g.\ boundaries between two alternatives) on $f$.

%%%%%%%%%%%%%%%%%%%%%%%%%%%%%%%%%%%%%%%%%%%%%%%%%%%%%%%%%%%%%
\subsection{Division into cases}\label{sec:cases_refined} %%%
%%%%%%%%%%%%%%%%%%%%%%%%%%%%%%%%%%%%%%%%%%%%%%%%%%%%%%%%%%%%%

Our starting point in proving Theorem \ref{thm:k_refined} is Lemma \ref{lem:boundaries2}. Clearly if \eqref{eq:neutralPairs3ManipProb} holds then we are done, so in the rest of Section \ref{sec:manip_ref} we assume that this is not the case. Then Lemma \ref{lem:boundaries2} tells us that \eqref{eq:k_inf_ref} holds, and w.l.o.g.\ we may assume that the two boundaries that the lemma gives us have $i=1$ and $j=2$. I.e.\ we have
\[
\p \left( \sigma \text{ on } B_1^{a,b;\left[a:b\right]} \right) \geq \frac{4\eps}{n k^7} \quad \text{and} \quad \p \left( \sigma \text{ on } B_2^{c,d;\left[c:d\right]} \right) \geq \frac{4\eps}{n k^7},
\]
where recall that $\sigma$ is on $B_1^{a,b;\left[a:b\right]}$ if $f\left( \sigma \right) = a$ and $f\left( \left[a:b\right]_1 \sigma \right) = b$. If $\sigma$ is on $B_1^{a,b;\left[a:b\right]}$ and $b \stackrel{\sigma_1}{>} a$, then $\sigma$ is a 2-manipulation point, so if this happens in more than half of the cases when $\sigma$ is on $B_1^{a,b;\left[a:b\right]}$, then we have
\[
\p\left( \sigma \in M_2 \right) \ge \frac{2\eps}{n k^7},
\]
and we are done. Similarly in the case of the boundary between $c$ and $d$ in coordinate 2. So we may assume from now on that
\[
\p \left( \sigma \in \cup_{z_{-1}^{a,b}} B_1 \left( z_{-1}^{a,b} \right) \right) \geq \frac{2\eps}{n k^7} \quad \text{and} \quad \p \left( \sigma \in \cup_{z_{-2}^{c,d}} B_2 \left( z_{-2}^{c,d} \right) \right) \geq \frac{2\eps}{n k^7}.
\]
The following lemma is an immediate corollary.
\begin{lemma}\label{lem:cases_ref}
Either
\begin{equation}\label{eq:sm_fbr_ref}
\p \left( \sigma \in \Sm \left( B_1^{a,b;\left[a:b\right]} \right) \right) \geq \frac{\eps}{n k^7}
\end{equation}
or
\begin{equation}\label{eq:lg_fbr_ref}
\p \left( \sigma \in \Lg \left( B_1^{a,b;\left[a:b\right]} \right) \right) \geq \frac{\eps}{n k^7},
\end{equation}
and the same can be said for the boundary $B_2^{c,d;\left[c:d\right]}$.
\end{lemma}

We distinguish cases based upon this: either \eqref{eq:sm_fbr_ref} holds, or \eqref{eq:sm_fbr_ref} holds for the boundary $B_2^{c,d;\left[c:d\right]}$, or \eqref{eq:lg_fbr_ref} holds for both boundaries. We only need one boundary for the small fiber case, and we need both boundaries only in the large fiber case. So in the large fiber case we must differentiate between two cases: whether $d\in \left\{ a, b \right\}$ or $d \notin \left\{a, b \right\}$. First of all, in the $d\notin \left\{a, b \right\}$ case the problem of finding a manipulation point with not too small (i.e.\ inverse polynomial in $n$, $k$ and $\eps^{-1}$) probability has already been solved by Isaksson, Kindler and Mossel \cite{isaksson2010geometry}, so we are primarily interested in the $d\in \left\{a,b \right\}$ case. But moreover, we will see that our method of proof works in both cases.

In the rest of the section we first deal with the small fiber case, and then with the large fiber case.

%%%%%%%%%%%%%%%%%%%%%%%%%%%%%%%%%%%%%%%%%%%%%%%%%%%%%%%
\subsection{Small fiber case}\label{sec:sm_fbr_ref} %%%
%%%%%%%%%%%%%%%%%%%%%%%%%%%%%%%%%%%%%%%%%%%%%%%%%%%%%%%

We now deal with the case when \eqref{eq:sm_fbr_ref} holds. We formalize the ideas of the outline in a series of statements.

First, we want to formalize that the boundaries of the boundaries are big in this refined graph setting as well, when we are on a small fiber. The proof uses the canonical path method, as successfully adapted to this setting by Isaksson, Kindler and Mossel \cite{isaksson2010geometry}, and is very similar to the proof of Lemma \ref{lem:comparable_bdries_1vot}, with some necessary modifications due to the fact that we now have $n$ coordinates.
\begin{lemma}\label{lem:comparable_bdries_ref}
Fix a coordinate and a pair of alternatives---for simplicity we choose coordinate 1 and alternatives $a$ and $b$, but we note that this lemma holds in general, we do not assume anything special about these choices. Let $z_{-1}^{a,b}$ be such that $B_1 \left( z_{-1}^{a,b} \right)$ is a small fiber for $B_1^{a,b;\left[a:b\right]}$. Then, writing $B \equiv B_1 \left( z_{-1}^{a,b} \right)$ for simplicity, we have
\begin{equation}\label{eq:comp_bdries_ref}
\p \left( \sigma \in \partial \left( B \right) \right) \geq \frac{\gamma}{2nk^5} \p \left( \sigma \in B \right).
\end{equation}
\end{lemma}

\begin{proof}
Let $B^c = \bar{F} \left( z_{-1}^{a,b} \right) \setminus B$. For every $\left( \sigma, \sigma' \right) \in B \times B^c$, we define a canonical path from $\sigma$ to $\sigma'$, which has to pass through at least one edge in $\partial_e \left( B \right)$. Then if we show that every edge in $\partial_e \left( B \right)$ lies on at most $r$ canonical paths, then it follows that $\left| \partial_e \left( B \right) \right| \geq \left|B \right| \left|B^c \right| / r$.

So let $\left( \sigma, \sigma' \right) \in B \times B^c$. We define a path from $\sigma$ to $\sigma'$ by applying a path construction in each coordinate one by one, and then concatenating these paths: first in the first coordinate we get from $\sigma_1$ to $\sigma'_1$, while leaving all other coordinates unchanged, then in the second coordinate we get from $\sigma_2$ to $\sigma'_2$, while leaving all other coordinates unchanged, and so on, finally in the last coordinate we get from $\sigma_n$ to $\sigma'_n$. In the first coordinate we apply the path construction of \cite[Proposition 6.4.]{isaksson2010geometry}, but considering the block formed by $a$ and $b$ as a single element; in all other coordinates we apply the path construction of \cite[Proposition 6.6.]{isaksson2010geometry}. Since this path goes from $\sigma$ (which is in $B$) to $\sigma'$ (which is in $B^c$), it must pass through at least one edge in $\partial_e \left( B \right)$.

For a given edge $\left( \pi, \pi' \right) \in \partial_e \left( B \right)$, at most how many possible $\left( \sigma, \sigma' \right) \in B \times B^c$ pairs are there such that the canonical path between $\sigma$ and $\sigma'$ defined above passes through $\left( \pi, \pi' \right)$? Let us differentiate two cases.

Suppose $\pi$ and $\pi'$ differ in the first coordinate. Then coordinates 2 through $n$ of $\sigma$ must agree with the respective coordinates of $\pi$, while coordinates 2 through $n$ of $\sigma'$ can be anything (up to the restriction given by $\sigma' \in B^c \subseteq \bar{F} \left( z_{-1}^{a,b} \right)$), giving $\left( \frac{k!}{2} \right)^{n-1}$ possibilities. Now fixing all coordinates except the first, \cite[Proposition 6.4.]{isaksson2010geometry} tells us that there are at most $\left( k - 1 \right)^2 \left( k -1 \right)! / 2 < k^2 \left( k - 1 \right)!$ possibilities for the pair $\left( \sigma_1, \sigma'_1 \right)$. So altogether there are at most $k^2 \left( k- 1 \right)! \left( \frac{k!}{2} \right)^{n-1}$ paths that pass through a given edge in $\partial_e \left( B \right)$ in this case.

Suppose now that $\pi$ and $\pi'$ differ in the $i^{\text{th}}$ coordinate, $i \neq 1$. Then the first $i-1$ coordinates of $\sigma'$ must agree with the first $i-1$ coordinates of $\pi$, while coordinates $i+1, \dots, n$ of $\sigma$ must agree with the respective coordinates of $\pi$. The first $i-1$ coordinates of $\sigma$, and coordinates $i+1, \dots, n$ of $\sigma'$ can be anything (up to the restriction given by $\sigma, \sigma' \in \bar{F} \left( z_{-1}^{a,b} \right)$), giving $\left( k - 1 \right)! \left( \frac{k!}{2} \right)^{n-2}$ possibilities. Now fixing all coordinates except the $i^{\text{th}}$ coordinate, \cite[Proposition 6.6.]{isaksson2010geometry} tells us that there are at most $k^4 k!$ possibilities for the pair $\left( \sigma_i, \sigma'_i \right)$. So altogether there are at most $2k^4 \left( k - 1 \right)! \left( \frac{k!}{2} \right)^{n-1}$ paths that pass through a given edge in $\partial_e \left( B \right)$ in this case.

So in any case, there are at most $2k^4 \left( k - 1 \right)! \left( \frac{k!}{2} \right)^{n-1}$ paths that pass through a given edge in $\partial_e \left( B \right)$.

Recall that $\left| \bar{F} \left( z_{-1}^{a,b} \right) \right| = \left( k - 1 \right)! \left( \frac{k!}{2} \right)^{n-1}$, and also $\left| B^c \right| \geq \gamma \left( k - 1 \right)! \left( \frac{k!}{2} \right)^{n-1}$ since $B$ is a small fiber. Therefore
\[
\left| \partial_e \left( B \right) \right| \geq \frac{\left| B \right| \left| B^c \right|}{ 2 k^4 \left( k - 1 \right)! \left( \frac{k!}{2} \right)^{n-1}} \geq \frac{\gamma}{2 k^4} \left| B \right|.
\]
Now in $G\left( z_{-1}^{a,b} \right)$ every ranking profile has no more than $n k$ neighbors, which implies \eqref{eq:comp_bdries_ref}.
\end{proof}

\begin{corollary}\label{cor:comparable_bdries_ref}
If \eqref{eq:sm_fbr_ref} holds, then
\[
\p \left( \sigma \in \bigcup_{z_{-1}^{a,b}} \partial \left( B_1 \left( z_{-1}^{a,b} \right) \right) \right) \geq \frac{\gamma \eps}{2 n^2 k^{12}}.
\]
\end{corollary}

\begin{proof}
Using the previous lemma and \eqref{eq:sm_fbr_ref} we have
\begin{align*}
\p \left( \sigma \in \bigcup_{z_{-1}^{a,b}} \partial \left( B_1 \left( z_{-1}^{a,b} \right) \right) \right) &= \sum_{z_{-1}^{a,b}} \p \left( \sigma \in \partial \left( B_1 \left( z_{-1}^{a,b} \right) \right) \right)\\
&\geq \sum_{z_{-1}^{a,b} : B_1 \left( z_{-1}^{a,b} \right)\subseteq \Sm \left( B_1^{a,b;\left[a:b\right]} \right) } \p \left( \sigma \in \partial \left( B_1 \left( z_{-1}^{a,b} \right) \right) \right)\\
&\geq \sum_{z_{-1}^{a,b} : B_1 \left( z_{-1}^{a,b} \right) \subseteq \Sm \left(  B_1^{a,b;\left[a:b\right]} \right) } \frac{\gamma}{2 n k^5} \p \left( \sigma \in  B_1 \left( z^{a,b} \right) \right)\\
&= \frac{\gamma}{2 n k^5} \p \left( \sigma \in \Sm \left( B_1^{a,b} \right) \right) \geq \frac{\gamma \eps}{2 n^2 k^{12}}. \qedhere
\end{align*}
\end{proof}

Next, we want to find manipulation points on the boundaries of boundaries.

Before we do this, let us divide the boundaries of the boundaries according to which direction they are in. If $\sigma \in \partial \left( B_1 \left( z_{-1}^{a,b} \right) \right)$ for some $z_{-1}^{a,b}$, then we know that there exists a ranking profile $\pi$ such that $\left( \sigma, \pi \right) \in \partial_{e} \left( B_1 \left( z_{-1}^{a,b} \right) \right)$. We know that $\sigma$ and $\pi$ differ in exactly one coordinate, say coordinate $j$; in this case we say that $\sigma$ is on the boundary of $B_1 \left( z_{-1}^{a,b} \right)$ in direction $j$, and we write $\sigma \in \partial_j \left( B_1 \left( z_{-1}^{a,b} \right) \right)$. (This notation should not be confused with that of the edge boundary.)

We can write the boundary of $B_1 \left( z_{-1}^{a,b} \right)$ as a union of boundaries in the different directions:
\[
\partial \left( B_1 \left( z_{-1}^{a,b} \right) \right) = \cup_{j=1}^{n} \partial_j \left( B_1 \left( z_{-1}^{a,b} \right) \right),
\]
but note that this is not (necessarily) a disjoint union, as a ranking profile $\sigma$ for which $\sigma \in \partial \left( B_1 \left( z_{-1}^{a,b} \right) \right)$ might lie on the boundary in multiple directions.

In particular, we differentiate between the boundary in direction 1 and the boundary in all other directions. To this end we introduce the notation
\[
\partial_{-1} \left( B_1 \left( x_{-1}^{a,b} \right) \right) := \cup_{j=2}^{n} \partial_j \left( B_1 \left( x_{-1}^{a,b} \right) \right).
\]
With this notation we have the following corollary of Corollary \ref{cor:comparable_bdries_ref}.

\begin{corollary}\label{cor:cases_bdry_dir_ref}
If \eqref{eq:sm_fbr_ref} holds, then either
\begin{equation}\label{eq:lg_bdry_in_dir!=1}
\p \left( \sigma \in \cup_{z_{-1}^{a,b}} \partial_{-1} \left( B_1 \left( z_{-1}^{a,b} \right) \right) \right) \geq \frac{\gamma \eps}{4n^2 k^{12}}
\end{equation}
or
\begin{equation}\label{eq:lg_bdry_in_dir=1}
\p \left( \sigma \in \cup_{z_{-1}^{a,b}} \partial_{1} \left( B_1 \left( z_{-1}^{a,b} \right) \right) \right) \geq \frac{\gamma \eps}{4n^2 k^{12}}.
\end{equation}
\end{corollary}

\begin{lemma}\label{lem:manip_on_bdry_of_bdry_in_dir!=1_ref}
Suppose the ranking profile $\sigma$ is on the boundary of a fiber for $B_1^{a,b;\left[a:b\right]}$ in direction $j \neq 1$, i.e.\
\[
\sigma \in \cup_{z_{-1}^{a,b}} \partial_{-1} \left( B_1 \left( z_{-1}^{a,b} \right) \right).
\]
Then there exists a 3-manipulation point $\hat{\sigma}$ which agrees with $\sigma$ in all coordinates except perhaps coordinate 1 and some coordinate $j \neq 1$; furthermore $\hat{\sigma}_1$ is equal to $\sigma_1$ or $\left[a:b\right]\sigma_1$, except that the position of a third alternative $c$ might be shifted arbitrarily, and $\hat{\sigma}_j$ is equal to $\sigma_j$ or $z \sigma_j$ for some adjacent transposition $z \in T$, except the position of $b$ might be shifted arbitrarily.
\end{lemma}

\begin{proof}
Suppose $x_{-1}^{a,b} \left( \sigma \right) = z_{-1}^{a,b}$.
Since $\sigma \in \partial \left( B_1 \left( z_{-1}^{a,b} \right) \right) \subseteq B_1 \left( z_{-1}^{a,b} \right)$, we know that $f \left( \sigma \right) = a$, and if $\sigma' = \left[a:b\right]_1 \sigma$, then $f\left( \sigma' \right) = b$.

Let $\pi = \left( \pi_j, \sigma_{-j} \right)$ denote the ranking profile such that $\left( \sigma, \pi \right) \in \partial_{e} \left( B_1 \left( z_{-1}^{a,b} \right) \right)$. Let $\pi' := \left[a:b\right]_1 \pi$. Since $\pi \notin B_1 \left( z_{-1}^{a,b} \right)$, $\left( f\left( \pi \right), f \left( \pi' \right) \right) \neq \left( a, b \right)$. Then, by Lemma \ref{lem:nonManipBoundary}, if $f\left( \pi \right) \neq f\left( \pi' \right)$, then one of $\pi$ and $\pi'$ is a 2-manipulation point.

So let us suppose that $f\left( \pi \right) = f \left( \pi' \right)$. If $f\left( \pi' \right) = a$, then one of $\sigma'$ and $\pi'$ is a 2-manipulation point by Lemma \ref{lem:nonManipBoundary}, since $\pi' = z_j \sigma'$ for some adjacent transposition $z \neq \left[a:b\right]$. If $f\left( \pi \right) = b$, then similarly one of $\sigma$ and $\pi$ is a 2-manipulation point.

Finally let us suppose that $f\left(\pi \right) = c$ for some $c\notin \left\{ a, b \right\}$. In this case Lemma \ref{lem:nonManipTriple} tells us that there exists an appropriate 3-manipulation point $\hat{\sigma}$.
\end{proof}

\begin{corollary}\label{cor:lg_bdry_in_dir!=1_gives_manip}
If \eqref{eq:lg_bdry_in_dir!=1} holds, then
\begin{equation}\label{eq:lg_bdry_in_dir!=1_gives_manip}
\p \left( \sigma \in M_3 \right) \geq \frac{\gamma \eps}{8 n^3 k^{16}}.
\end{equation}
\end{corollary}

\begin{proof}
Lemma \ref{lem:manip_on_bdry_of_bdry_in_dir!=1_ref} tells us that for every ranking profile $\sigma$ which is on the boundary of a fiber for $B_1^{a,b;\left[a:b\right]}$ in some direction $j\neq 1$, there is a 3-manipulation point $\hat{\sigma}$ ``nearby''; the lemma specifies what ``nearby'' means.

How many ranking profiles $\sigma$ may give the same $\hat{\sigma}$? At most $2 n k^4$, which comes from the following: $\sigma$ and $\hat{\sigma}$ agree in all coordinates except maybe two, one of which is the first coordinate; there are $n-1 < n$ possibilities for the other coordinate; in the first coordinate, $\hat{\sigma}_1$ is either $\sigma_1$ or $\left[a:b\right] \sigma_1$ (giving 2 possibilities), while some alternative $c$ ($k-2 < k$ possibilities) might be shifted arbitrarily (at most $k$ possibilities); in the other coordinate $j \neq 1$, $\hat{\sigma}_j$ is equal to $\sigma_j$ or $z \sigma_j$ for some adjacent transposition $z\in T$ (at most $k$ possibilities), except $b$ might be shifted arbitrarily ($k$ possibilities).

So putting this result from Lemma \ref{lem:manip_on_bdry_of_bdry_in_dir!=1_ref} together with \eqref{eq:lg_bdry_in_dir!=1} yields \eqref{eq:lg_bdry_in_dir!=1_gives_manip}.
\end{proof}

The remaining case we have to deal with is when \eqref{eq:lg_bdry_in_dir=1} holds.

\begin{lemma}\label{lem:bdry_of_bdry_in_coord1_ref}
Suppose the ranking profile $\sigma$ is on the boundary of a fiber for $B_1^{a,b;\left[a:b\right]}$ in direction 1, i.e.\
\[
\sigma \in \cup_{z_{-1}^{a,b}} \partial_{1} \left( B_1 \left( z_{-1}^{a,b} \right) \right).
\]
Then either $\sigma \in \LD_1 \left( a, b \right)$, or there exists a 3-manipulation point $\hat{\sigma}$ which agrees with $\sigma$ in all coordinates except perhaps in coordinate 1; furthermore $\hat{\sigma}_1$ is equal to $\sigma_1$, or $\left[a:b\right] \sigma_1$ except that the position of a third alternative $c$ might be shifted arbitrarily.
\end{lemma}

\begin{proof}
Just like the proof of Lemma \ref{lem:bdry_of_B}.
\end{proof}

The following corollary then tells us that either we have found many 3-manipulation points, or we have many local dictators on three alternatives in coordinate 1.

\begin{corollary}\label{cor:sm_fbr_ref_last_cor}
Suppose \eqref{eq:lg_bdry_in_dir=1} holds. Then either
\begin{equation}\label{eq:many_loc_dict}
\sum_{c\notin\left\{a,b\right\}} \p \left( \sigma \in \LD_1^{\left\{a,b,c\right\}} \right) = \p\left( \sigma \in \LD_1 \left( a, b \right) \right) \geq \frac{\gamma \eps}{8 n^2 k^{12}}
\end{equation}
or
\[
\p \left( \sigma \in M_3 \right) \geq \frac{\gamma \eps}{16 n^2 k^{14}}.
\]
\end{corollary}

%%%%%%%%%%%%%%%%%%%%%%%%%%%%%%%%%%%%%%%%%%%%%%%%%%%%%%%%%%%%%%%%%%%%
\subsubsection{Dealing with local dictators}\label{sec:loc_dict} %%%
%%%%%%%%%%%%%%%%%%%%%%%%%%%%%%%%%%%%%%%%%%%%%%%%%%%%%%%%%%%%%%%%%%%%

So the remaining case we have to deal with in this small fiber case is when \eqref{eq:many_loc_dict} holds, i.e.\ we have many local dictators in coordinate 1.

\begin{lemma}\label{lem:loc_dict_abc_to_top}
Suppose $\sigma \in \LD_1^{\left\{a,b,c\right\}}$ for some alternative $c\notin \left\{a,b\right\}$. Define $\sigma' := \left( \sigma'_1, \sigma_{-1} \right)$ by letting $\sigma'_1$ be equal to $\sigma_1$ except that the block of $a$, $b$ and $c$ is moved to the top of the coordinate. Then
\begin{itemize}
\item either $\sigma' \in \LD_1^{\left\{a,b,c\right\}}$,
\item or there exists a 3-manipulation point $\hat{\sigma}$ which agrees with $\sigma$ in all coordinates except perhaps in coordinate 1; furthermore $\hat{\sigma}_1$ is equal to $\sigma_1$ except that the position of $a$, $b$ and $c$ might be shifted arbitrarily.
\end{itemize}
\end{lemma}

\begin{proof}
Just like the proof of Lemma \ref{lem:loc_dict_abc_to_top_1vot}.
\end{proof}

\begin{corollary}\label{cor:many_dict_at_top}
If \eqref{eq:many_loc_dict} holds, then either
\begin{equation}\label{eq:loc_dict_abc_on_top}
\sum_{c \notin \left\{ a, b \right\}} \p \left( \sigma \in LD_1^{\left\{a,b,c\right\}}, \left\{ \sigma_1 \left( 1 \right), \sigma_1 \left( 2 \right), \sigma_1 \left( 3 \right) \right\} = \left\{a,b,c\right\} \right) \geq \frac{\gamma \eps}{16 n^2 k^{13}}
\end{equation}
or
\[
\p \left( \sigma \in M_3 \right) \geq \frac{\gamma \eps}{16 n^2 k^{15}}.
\]
\end{corollary}

\begin{proof}
Just like the proof of Corollary \ref{cor:many_dict_at_top_1vot}.
\end{proof}

Now \eqref{eq:loc_dict_abc_on_top} is equivalent to
\begin{equation}\label{eq:loc_dict_abc_on_top2}
\sum_{c \notin \left\{ a, b \right\}} \p \left( \sigma \in LD_1^{\left\{a,b,c\right\}}, \left( \sigma_1 \left( 1 \right), \sigma_1 \left( 2 \right), \sigma_1 \left( 3 \right) \right) = \left(a,b,c\right) \right) \geq \frac{\gamma \eps}{96 n^2 k^{13}}.
\end{equation}
We know that
\[
\p \left( \left( \sigma_1 \left( 1 \right), \sigma_1 \left( 2 \right), \sigma_1 \left( 3 \right) \right) = \left(a,b,c\right) \right) = \frac{1}{k \left( k-1 \right) \left( k-2 \right)} \leq \frac{6}{k^3},
\]
and so \eqref{eq:loc_dict_abc_on_top2} implies (recall Definition \ref{def:cond})
\begin{equation}\label{eq:loc_dict_abc_on_top3}
\sum_{c \notin \left\{ a, b \right\}} \p_1^{\left(a,b,c\right)} \left( \sigma \in LD_1^{\left\{a,b,c\right\}}  \right) \geq \frac{\gamma \eps}{576 n^2 k^{10}}.
\end{equation}

Now fix an alternative $c \notin \left\{a,b\right\}$ and define the graph $G_{\left( a, b, c \right)} = \left( V_{\left( a, b, c \right)}, E_{\left( a, b, c \right)} \right)$ to have vertex set
\[
V_{\left( a, b, c \right)} := \left\{ \sigma \in S_k^n : \left( \sigma_1 \left( 1 \right), \sigma_1 \left( 2 \right), \sigma_1 \left( 3 \right) \right) = \left(a,b,c\right) \right\}
\]
and for $\sigma, \pi \in V_{\left( a, b, c \right)}$ let $\left( \sigma, \pi \right) \in E_{\left( a, b, c \right)}$ if and only if $\sigma$ and $\pi$ differ in exactly one coordinate, and by an adjacent transposition in this coordinate. So $G_{\left( a, b, c \right)}$ is the subgraph of the refined rankings graph induced by the vertex set $V_{\left( a, b, c \right)}$.

Let
\[
T_1 \left( a, b, c \right) := V_{\left( a, b, c \right)} \cap LD_1^{\left\{a,b,c\right\}},
\]
and let $\partial_e \left( T_1 \left( a, b, c \right) \right)$ and $\partial \left( T_1 \left( a, b, c \right) \right)$ denote the edge and vertex boundary of $T_1 \left( a, b, c \right)$ in $G_{\left( a, b, c \right)}$, respectively.

The next lemma shows that unless $T_1 \left( a, b, c \right)$ is almost all of $V_{\left( a, b, c \right)}$, the size of the boundary $\partial \left( T_1 \left( a, b, c \right) \right)$ is comparable to the size of $T_1 \left( a, b, c \right)$. The proof uses a canonical path argument, just like in Lemma \ref{lem:comparable_bdries_ref}, and is very similar to the proof of Lemma \ref{lem:lg_bdry_for_loc_dict_1vot}.

\begin{lemma}\label{lem:lg_bdry_for_loc_dict}
Let $c \notin \left\{a,b \right\}$ be arbitrary. Write $T \equiv T_1 \left( a, b, c \right)$ for simplicity. If $\p_1^{\left(a,b,c\right)} \left( \sigma \in T \right) \leq 1 - \delta$, then
\begin{equation}\label{eq:lg_bdry_for_loc_dict}
\p_1^{\left( a, b, c \right)} \left( \sigma \in \partial \left( T \right) \right) \geq \frac{\delta}{n k^3} \p_1^{\left( a, b, c \right)} \left( \sigma \in T \right).
\end{equation}
\end{lemma}

\begin{proof}
The proof is essentially the same as the proof of Lemma \ref{lem:lg_bdry_for_loc_dict_1vot}, with a slight modification to deal with $n$ coordinates. Let $T^c = V_{\left( a, b, c \right)} \setminus T \left( a, b, c \right)$. For every $\left( \sigma, \sigma' \right) \in T \times T^c$ we define a canonical path from $\sigma$ to $\sigma'$ by applying a path construction in each coordinate one by one, and then concatenating these paths. In all coordinates we apply the path construction of \cite[Proposition 6.4.]{isaksson2010geometry}, but in the first coordinate we only apply it to alternatives $\left[k \right] \setminus \left\{ a, b, c \right\}$.

The analysis of this construction is done in exactly the same way as in Lemma \ref{lem:comparable_bdries_ref}; in the end we get that $\left| \partial_e \left( T \right) \right| \geq \frac{\delta}{k^2} \left| T \right|$. Now every vertex in $V_{\left(a,b,c\right)}$ has no more than $nk$ neighbors, which implies \eqref{eq:lg_bdry_for_loc_dict}.
\end{proof}

The next lemma tells us that if $\sigma$ is on the boundary of a set of local dictators on $\left\{a,b,c \right\}$ for some alternative $c \notin \left\{a,b\right\}$ in coordinate 1, then there is a 4-manipulation point $\hat{\sigma}$ which is close to $\sigma$. The proof is similar to that of Lemma \ref{lem:manip_pts_on_bdry_of_loc_dict_1vot}, but we have to take care of all $n$ coordinates.

\begin{lemma}\label{lem:manip_pts_on_bdry_of_loc_dict}
Suppose $\sigma \in \partial \left( T_1 \left( a,b, c \right) \right)$ for some $c\notin \left\{a,b\right\}$. We distinguish two cases, based on the number of alternatives.

If $k=3$, then there exists a (3-)manipulation point $\hat{\sigma}$ which differs from $\sigma$ in at most two coordinates, one of them being the first coordinate.

If $k \geq 4$, then there exists a 4-manipulation point $\hat{\sigma}$ which differs from $\sigma$ in at most two coordinates, one of them being the first coordinate; furthermore, $\hat{\sigma}_1$ is equal to $\sigma_1$ except that the order of the block of $a$, $b$ and $c$ might be rearranged and an additional alternative $d$ might be shifted arbitrarily; and in the other coordinate, call it $j$, $\hat{\sigma}_j$ is equal to $\sigma_j$ except perhaps $a$, $b$ and $c$ are shifted arbitrarily.
\end{lemma}

\begin{proof}
Let $\pi$ be the ranking profile such that $\left( \sigma, \pi \right) \in \partial_e \left( T_1 \left( a, b , c \right) \right)$, let $j$ be the coordinate in which they differ, and let $z$ be the adjacent transposition in which they differ, i.e.\ $\pi = z_j \sigma$. Since $\pi \notin T_1 \left( a, b, c \right)$, there exists a reordering of the block of $a$, $b$, and $c$ at the top of $\pi_1$ such that the outcome of $f$ is not the top ranked alternative in coordinate 1. Call the resulting vector $\pi'_1$, and let $\pi':= \left( \pi'_1, \pi_{-1} \right)$. W.l.o.g.\ let us assume that $\pi'_1 \left( 1 \right) = a$. Let us also define $\sigma' := z_j \pi'$. We distinguish two cases: $j=1$ and $j\neq 1$.

If $j=1$ (in which case we must have $k \geq 5$), $\pi'$ is a 2-manipulation point, since $f\left( \sigma' \right) = a$.

If $j\neq 1$, then there are various cases to consider. If the adjacent transposition $z$ does not move $a$, then either $\pi'$ or $\sigma'$ is a 2-manipulation point. So let us suppose that $z = \left[a:d\right]$ for some $d \neq a$.

Clearly we must have $f\left( \pi' \right) = d$, or else $\pi'$ or $\sigma'$ is a 2-manipulation point. Suppose first that $d\in \left\{ b, c \right\}$. W.l.o.g.\ suppose that $d=b$.

Then take alternative $c$ in coordinate $j$ of both $\sigma'$ and $\pi'$, and bubble it to the block of $a$ and $b$ simultaneously in the two ranking profiles. If along the way the value of the outcome of the SCF $f$ changes from $a$ or $b$, respectively, then we have a 2-manipulation point by Lemma \ref{lem:nonManipBoundary}. Otherwise, we now have $a$, $b$, and $c$ adjacent in both coordinates 1 and $j$. Now rearranging the order of the blocks of $a$, $b$, and $c$ in these two coordinates (which can be done using adjacent transpositions), we either get a 2-manipulation point by Lemma \ref{lem:nonManipBoundary}, or we can define a new SCF on two voters and three alternatives, $a$, $b$, and $c$. This SCF takes on three values and it is also not hard to see that the outcome is not only a function of the first coordinate, so by the Gibbard-Satterthwaite theorem we know that this SCF has a manipulation point, which is a 3-manipulation point of the original SCF $f$.

Now let us look at the case when $d \notin \left\{b,c\right\}$. In this case we do something similar to what we just did in the previous paragraph. In both $\sigma'$ and $\pi'$, first bubble up alternative $d$ in coordinate 1 up to the block of $a$, $b$, and $c$, and then bubble $b$ and $c$ in coordinate $j$ to the block of $a$ and $d$. All of this using adjacent transpositions. If the value of the outcome of the SCF $f$ changes from $a$ or $d$, respectively, at any time along the way, then we have a 2-manipulation point by Lemma \ref{lem:nonManipBoundary}. Otherwise, we now have $a$, $b$, $c$ and $d$ adjacent in both coordinates 1 and $j$, and we can apply the same trick to find a 4-manipulation point, using the Gibbard-Satterthwaite theorem.
\end{proof}

The next corollary puts together Corollary \ref{cor:many_dict_at_top} and Lemmas \ref{lem:lg_bdry_for_loc_dict} and \ref{lem:manip_pts_on_bdry_of_loc_dict}.

\begin{corollary}\label{cor:manip_by_bdry_loc_dict}
Suppose \eqref{eq:loc_dict_abc_on_top} holds. Then if for every $c \notin \left\{a,b\right\}$ we have $\p_1^{\left(a,b,c\right)} \left( \sigma \in T_1 \left( a,b, c \right) \right) \leq 1 - \frac{\eps}{100 k}$, then
\[
\p \left( \sigma \in M_4 \right) \geq \frac{\gamma \eps^2}{345600 n^4 k^{22}}.
\]
\end{corollary}

\begin{proof}
We know that \eqref{eq:loc_dict_abc_on_top} implies
\[
\sum_{c \notin \left\{a,b \right\}} \p_1^{a,b,c} \left( \sigma \in T_1 \left( a,b,c\right) \right) \geq \frac{\gamma \eps}{576 n^2 k^{10}}.
\]
Now then using the assumptions, Lemma \ref{lem:lg_bdry_for_loc_dict} with $\delta = \frac{\eps}{100k}$ and Lemma \ref{lem:manip_pts_on_bdry_of_loc_dict}, we have
\begin{align*}
\p \left( \sigma \in M_4 \right) &\geq \sum_{c \neq \left\{a,b \right\}} \frac{1}{k^3} \p_1^{\left(a,b,c \right)} \left( \sigma \in M_4 \right) \geq \sum_{c\notin \left\{a,b\right\}} \frac{1}{6nk^8} \p_1^{\left( a, b, c \right)} \left( \sigma \in \partial \left( T_1 \left( a,b,c \right) \right) \right)\\
&\geq \sum_{c \notin \left\{a,b\right\}} \frac{\eps}{600 n^2 k^{12}} \p_1^{\left( a, b , c \right)} \left( \sigma \in T_1 \left( a, b, c \right) \right) \geq \frac{\gamma \eps^2}{345600 n^4 k^{22}}. \qedhere
\end{align*}
\end{proof}

So again we are left with one case to deal with: if there exists an alternative $c \notin \left\{ a, b \right\}$ such that $\p_1^{\left(a,b,c\right)} \left( \sigma \in T_1 \left( a,b, c \right) \right) > 1 - \frac{\eps}{100 k}$. Define a subset of alternatives $K \subseteq \left[k \right]$ in the following way:
\[
K := \left\{a,b \right\} \cup \left\{ c \in \left[k \right] \setminus \left\{a,b \right\} : \p_1^{\left(a,b,c\right)} \left( \sigma \in T_1 \left( a,b, c \right) \right) > 1 - \frac{\eps}{100 k} \right\}.
\]
In addition to $a$ and $b$, $K$ contains those alternatives that whenever they are at the top of coordinate 1 with $a$ and $b$, they form a local dictator with high probability.

So our assumption now is that $\left| K \right| \geq 3$.

Our next step is to show that unless we have many manipulation points, for any alternative $c \in K$, conditioned on $c$ being at the top of the first coordinate, the outcome of $f$ is $c$ with probability close to 1.

\begin{lemma}\label{lem:cond_on_top}
Let $c \in K$. Then either
\begin{equation}\label{eq:cond_on_top1}
\p_1^{\left( c \right)} \left( f \left( \sigma \right) = c \right) \geq 1 - \frac{\eps}{50 k},
\end{equation}
or
\begin{equation}\label{eq:else_2manip}
\p \left( \sigma \in M_2 \right) \geq \frac{\eps}{100 k^4}.
\end{equation}
\end{lemma}

\begin{proof}
Just like the proof of Lemma \ref{lem:cond_on_top_1vot}.
\end{proof}

We now deal with alternatives that are not in $K$: either we have many manipulation points, or for any alternative $d \notin K$, the outcome of $f$ is \emph{not} $d$ with probability close to 1.

\begin{lemma}\label{lem:d_notin_K_ref}
Let $d \notin K$. If $\p \left( f \left( \sigma \right) = d \right) \geq \frac{\eps}{4k}$, then
\[
\p \left( \sigma \in M_4 \right) \geq \frac{\eps^2}{10^6 n^2 k^{13}}.
\]
\end{lemma}

\begin{proof}
The proof is very similar to that of Lemma \ref{lem:d_notin_K_1vot}: we do the same steps in the first coordinate as done in the proof of Lemma \ref{lem:d_notin_K_1vot}, and the fact that we have $n$ coordinates only matters at the very end.

Let $\sigma$ be such that $f\left( \sigma \right) = d$. We will keep coordinates 2 through $n$ to be fixed as $\sigma_{-1}$ throughout the proof. By bubbling alternatives $d$, $a$, and $b$ in the first coordinate, we can define $\sigma'$, $\sigma^{\left(d, b, a \right)}$, $\sigma^{\left(d, a, b \right)}$, $\sigma^{\left( a, b, d \right)}$, $\sigma^{\left( a, d, b \right)}$, $\sigma^{\left( b, a, d \right)}$, and $\sigma^{\left( b, d, a \right)}$ just as in  the proof of Lemma \ref{lem:d_notin_K_1vot}.  Again, we can show that either
\[
\p \left( \sigma \in M_2 \right) \geq \frac{\eps}{1600 k^3},
\]
in which case we are done, or
\begin{equation}\label{eq:loc_dict_abd}
\p_1^{\left( a, b, d \right)} \left( \sigma^{\left(a,b,d\right)} \in LD_1^{\left\{a,b,d\right\}} \right) = \p \left( \sigma: \sigma^{\left( a, b, d \right)} \in LD_1^{\left\{a,b,d\right\}} \right) \geq \frac{\eps}{1600 k}.
\end{equation}
Define $G_{\left(a,b,d\right)}$ and $T_{\left( a, b, d \right)}$ analogously to $G_{\left( a, b, c \right)}$ and $T_{\left( a, b, c \right)}$, respectively.

Suppose that \eqref{eq:loc_dict_abd} holds. We also know that $d\notin K$, so Lemma \ref{lem:lg_bdry_for_loc_dict} applies, and then Lemma \ref{lem:manip_pts_on_bdry_of_loc_dict} shows us how to find manipulation points. We can put these arguments together, just like in the proof of Corollary \ref{cor:manip_by_bdry_loc_dict}, to show what we need:
\begin{align*}
\p \left( \sigma \in M_4 \right) &\geq \frac{1}{k^3} \p_1^{\left(a,b,d \right)} \left( \sigma \in M_4 \right) \geq \frac{1}{6nk^8} \p_1^{\left( a, b, d \right)} \left( \sigma \in \partial \left( T_1 \left( a,b,d \right) \right) \right)\\
&\geq \frac{\eps}{600 n^2 k^{12}} \p_1^{\left( a, b , d \right)} \left( \sigma \in T_1 \left( a, b, d \right) \right) \geq \frac{\eps^2}{10^6 n^2 k^{13}}. \qedhere
\end{align*}
\end{proof}

Putting together the results of the previous lemmas, there is only one case to be covered, which is covered by the following final lemma. Basically, this lemma says that unless there are enough manipulation points, our function is close to a dictator in the first coordinate, on the subset of alternatives $K$.
\begin{lemma}\label{lem:final_loc_dict}
Recall that we assume that $\Dist \left( f, \overline{\NONMANIP} \right) \geq \eps$. Furthermore assume that $\left| K \right| \geq 3$, for every $c \in K$ we have
\begin{equation}\label{eq:dict_cond_on_top}
\p_1^{\left( c \right)} \left( f \left( \sigma \right) = c \right) \geq 1 - \frac{\eps}{50k},
\end{equation}
and for every $d \notin K$ we have
\[
\p \left( f\left( \sigma \right) = d \right) \leq \frac{\eps}{4k}.
\]
Then
\begin{equation}\label{eq:final_manip}
\p \left( \sigma \in M_2 \right) \geq \frac{\eps}{4k^2}.
\end{equation}
\end{lemma}

\begin{proof}
Just like the proof of Lemma \ref{lem:final_loc_dict}.
\end{proof}

To conclude the proof in the small fiber case, inspect all the lower bounds for $\p \left( \sigma \in M_4 \right)$ obtained in Section \ref{sec:sm_fbr_ref}, and recall that $\gamma = \frac{\eps^3}{10^3 n^3 k^{24}}$.

%%%%%%%%%%%%%%%%%%%%%%%%%%%%%%%%%%%%%%%%%%%%%%%%%%%%%%%
\subsection{Large fiber case}\label{sec:lg_fbr_ref} %%%
%%%%%%%%%%%%%%%%%%%%%%%%%%%%%%%%%%%%%%%%%%%%%%%%%%%%%%%

We now deal with the large fiber case, when \eqref{eq:lg_fbr_ref} holds for both boundaries, i.e.\ when
\[
\p \left( \sigma \in \Lg \left( B_1^{a,b;\left[a:b\right]} \right) \right) \geq \frac{\eps}{nk^7}
\]
and
\[
\p \left( \sigma \in  \Lg \left( B_2^{c,d;\left[c:d\right]} \right)  \right) \geq \frac{\eps}{nk^7}.
\]
We differentiate between two cases: whether $d \in \left\{a,b\right\}$ or $d \notin \left\{a,b\right\}$. As mentioned before, the $d \notin \left\{a,b\right\}$ case has already been solved by Isaksson, Kindler and Mossel \cite{isaksson2010geometry}. But moreover, we will see that our method of proof works in both cases.

\subsubsection{Case 1}
Suppose $d \in \left\{a,b \right\}$, in which case w.l.o.g.\ we may assume that $d = a$. That is, in the rest of this case we may assume that
\begin{equation}\label{eq:lg_fbr_d=a_1}
\p \left( \sigma \in  \Lg \left( B_1^{a,b;\left[a:b\right]} \right) \right) \geq \frac{\eps}{nk^7}
\end{equation}
and
\begin{equation}\label{eq:lg_fbr_d=a_2}
\p \left( \sigma \in \Lg \left( B_2^{a,c;\left[a:c\right]} \right) \right) \geq \frac{\eps}{nk^7}.
\end{equation}

First, let us look at only the boundary between $a$ and $b$ in direction 1. Let us fix a vector $z_{-1}^{a,b}$ which gives a large fiber $B_1 \left( z_{-1}^{a,b} \right)$ for the boundary $B_1^{a,b;\left[a:b\right]}$, i.e.\ we know that
\begin{equation}\label{eq:lg_fbr_def2}
\p \left( \sigma \in B_1 \left( z_{-1}^{a,b} \right) \, \middle| \, \sigma \in \bar{F} \left( z_{-1}^{a,b} \right) \right) \geq 1 - \gamma.
\end{equation}

Our basic goal in the following will be to show that conditional on the ranking profile $\sigma$ being in the fiber $F \left( z_{-1}^{a,b} \right)$ (but not necessarily in $\bar{F} \left( z_{-1}^{a,b} \right)$), with high probability the outcome of the vote is $\tp_{\left\{a,b\right\}} \left( \sigma_1 \right)$, or else we have a lot of 2-manipulation points or local dictators on three alternatives in coordinate 1.

Our first step towards this is the following.
\begin{lemma}\label{lem:cond_ab_top_still_lg_fbr}
Suppose $z_{-1}^{a,b}$ gives a large fiber $B_1 \left( z_{-1}^{a,b} \right)$ for the boundary $B_1^{a,b;\left[a:b\right]}$. Then
\begin{equation}\label{eq:cond_ab_top_still_lg_fbr}
\p_1^{\left(a,b\right)} \left( \sigma \in B_1 \left( z_{-1}^{a,b} \right) \, \middle| \, \sigma \in F\left( z_{-1}^{a,b} \right) \right) \geq 1 - k \gamma.
\end{equation}
\end{lemma}

\begin{proof}
We know that
\[
\p \left( \left( \sigma_1 \left( 1 \right), \sigma_1 \left( 2 \right) \right) = \left( a, b \right) \, \middle| \, \sigma \in \bar{F} \left( z_{-1}^{a,b} \right) \right) = \frac{1}{k-1},
\]
and so
\begin{align*}
\p_1^{\left(a,b\right)} &\left( \sigma \notin B_1 \left( z_{-1}^{a,b} \right)\, \middle| \, \sigma \in F \left( z_{-1}^{a,b} \right) \right) = \p_1^{\left(a,b\right)} \left( \sigma \notin B_1 \left( z_{-1}^{a,b} \right)\, \middle| \, \sigma \in \bar{F}\left( z_{-1}^{a,b} \right) \right)\\
&= \left( k - 1 \right) \p \left( \sigma \notin B_1 \left( z_{-1}^{a,b} \right), \left( \sigma_1 \left( 1 \right), \sigma_1 \left( 2 \right) \right) = \left( a, b \right) \, \middle| \, \sigma \in \bar{F} \left( z_{-1}^{a,b} \right) \right) \leq \left( k - 1 \right) \gamma < k \gamma. \qedhere
\end{align*}
\end{proof}

The next lemma formalizes our goal mentioned above.
\begin{lemma}\label{lem:towards_inv_hyp_contr}
Suppose $z_{-1}^{a,b}$ gives a large fiber $B_1 \left( z_{-1}^{a,b} \right)$ for the boundary $B_1^{a,b;\left[a:b\right]}$. Then either
\begin{equation}\label{eq:f=top_ab_1}
\p \left( f\left( \sigma \right) = \tp_{\left\{a,b \right\}} \left( \sigma_1 \right) \, \middle| \, \sigma \in F \left( z_{-1}^{a,b} \right) \right) \geq 1 - 2 k \gamma
\end{equation}
or
\begin{equation}\label{eq:2manip_again}
\p \left( \sigma \in M_2 \, \middle| \, \sigma \in F \left( z_{-1}^{a,b} \right) \right) \geq \frac{\gamma}{2k}
\end{equation}
or
\begin{equation}\label{eq:loc_dict_again}
\p \left( \sigma \in LD_1 \left(a,b\right) \, \middle| \, \sigma \in F \left( z_{-1}^{a,b} \right) \right) \geq \frac{\gamma}{2k}.
\end{equation}
\end{lemma}

\begin{proof}
The proof of this lemma is essentially the same as that of Lemma \ref{lem:lg_fbr_final_1vot}, there are only two slight differences. First, we use Lemma \ref{lem:cond_ab_top_still_lg_fbr} to know that \eqref{eq:cond_ab_top_still_lg_fbr} holds. Second, we take $\sigma \in F \left( z_{-1}^{a,b} \right)$ to be uniform, and we stay on the fiber $F \left( z_{-1}^{a,b} \right)$ throughout the proof: we modify only the first coordinate throughout the proof, in the same way as we did for Lemma \ref{lem:lg_fbr_final_1vot}. We omit the details.
\end{proof}

Now this lemma holds for all vectors $z_{-1}^{a,b}$ which give a large fiber $B_1 \left( z_{-1}^{a,b} \right)$ for the boundary $B_1^{a,b;\left[a:b\right]}$. By \eqref{eq:lg_fbr_d=a_1} we know that
\[
\p \left( \sigma : B_1 \left( x_{-1}^{a,b} \left( \sigma \right) \right) \text{ is a large fiber} \right)\geq \frac{\eps}{nk^7}.
\]
Now if \eqref{eq:2manip_again} holds for at least a third of the vectors $z_{-1}^{a,b}$ that give a large fiber $B_1 \left( z_{-1}^{a,b} \right)$, then it follows that
\[
\p \left( \sigma \in M_2 \right) \geq \frac{\gamma \eps}{6n k^8}
\]
and we are done. If \eqref{eq:loc_dict_again} holds for at least a third of the vectors $z_{-1}^{a,b}$ that give a large fiber $B_1 \left( z_{-1}^{a,b} \right)$, then similarly we have
\[
\p \left( \sigma \in LD_1 \left( a, b \right) \right) \geq \frac{\gamma \eps}{6nk^8},
\]
which means that \eqref{eq:many_loc_dict} also holds, and so we are done by the argument in Section \ref{sec:loc_dict}.

So the remaining case to consider is when \eqref{eq:f=top_ab_1} holds for at least a third of the vectors $z_{-1}^{a,b}$ that give a large fiber $B_1 \left( z_{-1}^{a,b} \right)$.

We can go through this same argument for the boundary between $a$ and $c$ in direction 2 as well, and either we are done because
\[
\p \left( \sigma \in M_2 \right) \geq \frac{\gamma \eps}{6n k^8}
\]
or
\[
\p \left( \sigma \in LD_2 \left( a, c \right) \right) \geq \frac{\gamma \eps}{6nk^8},
\]
or for at least a third of the vectors $z_{-2}^{a,c}$ that give a large fiber $B_2 \left( z_{-2}^{a,c} \right)$ we have
\[
\p \left( f \left( \sigma \right) = \tp_{\left\{a,c \right\}} \left( \sigma_2 \right) \, \middle| \, \sigma \in F \left( z_{-2}^{a,c} \right) \right) \geq 1 - 2 k \gamma.
\]

So basically our final case is if
\begin{equation}\label{eq:lg_fbr_final_case_1}
\p \left( \sigma \in F_1^{a,b} \right) \geq \frac{\eps}{3 n k^7}
\end{equation}
and also
\begin{equation}\label{eq:lg_fbr_final_case_2}
\p \left( \sigma \in F_2^{a,c} \right) \geq \frac{\eps}{3 n k^7}.
\end{equation}
Notice that being in the set $F_1^{a,b}$ only depends on the vector $x^{a,b} \left( \sigma \right)$ of preferences between $a$ and $b$, and similarly being in the set $F_2^{a,c}$ only depends on the vector $x^{a,c} \left( \sigma \right)$ of preferences between $a$ and $c$. We know that $\left\{\left( x_i^{a,b} \left( \sigma \right), x_i^{a,c} \left( \sigma \right) \right) \right\}_{i=1}^{n}$ are independent, and for any given $i$ we know that $\left| \E \left( x_i^{a,b} \left( \sigma \right) x_i^{a,c} \left( \sigma \right) \right) \right| = \frac{1}{3}$. Hence we can apply reverse hypercontractivity (Lemma \ref{lem:invhypcontr}), to get the following result.

\begin{lemma}\label{lem:applying_inv_hyp_contr}
If \eqref{eq:lg_fbr_final_case_1} and \eqref{eq:lg_fbr_final_case_2} hold, then also
\begin{equation}\label{eq:result_of_inv_hyp_contr}
\p \left( \sigma \in F_1^{a,b} \cap F_2^{a,c} \right) \geq \frac{\eps^3}{27 n^3 k^{21}}.
\end{equation}
\end{lemma}

\begin{proof}
See above.
\end{proof}

The next and final lemma then concludes that we have lots of manipulation points.
\begin{lemma}\label{lem:inv_hyp_contr_end}
Suppose \eqref{eq:result_of_inv_hyp_contr} holds. Then
\begin{equation}\label{eq:final_manip2}
\p \left( \sigma \in M_3 \right) \geq \frac{\eps^3}{54n^3 k^{27}} - \frac{9 \gamma}{k^3}.
\end{equation}
\end{lemma}

\begin{proof}
First let us define two events:
\begin{align*}
I_1 &:= \left\{ \sigma : f \left( \sigma \right) = \tp_{\left\{a,b \right\}} \left( \sigma_1 \right) \right\}\\
I_2 &:= \left\{ \sigma : f \left( \sigma \right) = \tp_{\left\{a,c \right\}} \left( \sigma_2 \right) \right\}.
\end{align*}
Using similar estimates as previously in Lemma \ref{lem:lg_fbr_1}, we have
\begin{multline*}
\p \left( \sigma \in I_1 \cap I_2 \cap F_1^{a,b} \cap F_2^{a,c} \right) \geq \p \left( \sigma \in F_1^{a,b} \cap F_2^{a,c} \right)\\
 - \p \left( \sigma \notin I_1, \sigma \in F_1^{a,b} \cap F_2^{a,c} \right) - \p \left( \sigma \notin I_2, \sigma \in F_1^{a,b} \cap F_2^{a,c} \right).
\end{multline*}
The first term is bounded below via \eqref{eq:result_of_inv_hyp_contr}, while the other two terms can be bounded using the definition of $F_1^{a,b}$ and $F_2^{a,c}$, respectively:
\[
\p \left( \sigma \notin I_1, \sigma \in F_1^{a,b} \cap F_2^{a,c} \right) \leq \p \left( \sigma \notin I_1, \sigma \in F_1^{a,b} \right) \leq \p \left( \sigma \notin I_1 \, \middle| \, \sigma \in F_1^{a,b} \right) \leq 2 k \gamma,
\]
and similarly for the other term. Putting everything together gives us
\[
\p \left( \sigma \in I_1 \cap I_2 \cap F_1^{a,b} \cap F_2^{a,c} \right) \geq \frac{\eps^3}{27 n^3 k^{21}} - 4k \gamma.
\]
If $\sigma \in I_1 \cap I_2 \cap F_1^{a,b} \cap F_2^{a,c}$, then clearly we must have $f\left( \sigma \right) = a$, and therefore $x_1^{a,b} \left( \sigma \right) = 1$ and $x_2^{a,c} \left( \sigma \right) = 1$. Now define $\sigma'$ from $\sigma$ by bubbling up $b$ in coordinate 1 to just below $a$, and bubbling up $c$ in coordinate 2 to just below $a$. Either we encounter a 2-manipulation point along the way, or the outcome is still $a$: $f\left( \sigma' \right) = a$. If we encounter a 2-manipulation point along the way for at least half of such ranking profiles, then we are done:
\[
\p \left( \sigma \in M_2 \right) \geq \frac{1}{k^2} \left( \frac{\eps^3}{54 n^3 k^{21}} - 2k \gamma \right) = \frac{\eps^3}{54 n^3 k^{23}} - \frac{2 \gamma}{k}.
\]
Otherwise, we may assume that
\[
\p \left( \sigma \in I_1 \cap I_2 \cap F_1^{a,b} \cap F_2^{a,c}, f\left( \sigma' \right) = a  \right) \geq \frac{\eps^3}{54 n^3 k^{21}} - 2k \gamma.
\]
In this case define $\tilde{\sigma}' := \left[a:b\right]_1 \sigma'$ and $\tilde{\sigma}'' := \left[a:c\right]_2 \sigma'$.
%\begin{align*}
%\tilde{\sigma}' &:= \left[a:b\right]_1 \sigma'\\
%\tilde{\sigma}'' &:= \left[a:c\right]_2 \sigma'.
%\end{align*}
If $f \left( \tilde{\sigma}' \right) \notin \left\{a,b \right\}$ or $f \left( \tilde{\sigma}'' \right) \notin \left\{a,c \right\}$, then we automatically have that one of $\sigma', \tilde{\sigma}', \tilde{\sigma}''$ is a 2-manipulation point. If $f \left( \tilde{\sigma}' \right) = b$ and $f\left( \tilde{\sigma}'' \right) = c$, then by Lemma \ref{lem:nonManipTriple} we know that there exists a 3-manipulation point $\hat{\sigma}$ which agrees with $\sigma$ except perhaps $a$, $b$, and $c$ could be arbitrarily shifted in the first two coordinates. The final case is when $a \in \left\{ f\left( \tilde{\sigma}' \right), f \left( \tilde{\sigma}'' \right) \right\}$. But we now show that this has small probability, and therefore \eqref{eq:final_manip2} follows.

First let us look at the case of $f \left( \tilde{\sigma}' \right) = a$. We have
\begin{multline*}
\p \left( \sigma \in I_1 \cap I_2 \cap F_1^{a,b} \cap F_2^{a,c}, f\left( \sigma' \right) = a, f\left( \tilde{\sigma}' \right) = a  \right)\\
\begin{aligned}
&= \sum_{z_{-1}^{a,b} : F \left( z_{-1}^{a,b} \right) \subseteq F_1^{a,b}} \p \left( \sigma \in I_1 \cap I_2 \cap F \left( \left( 1, z_{-1}^{a,b} \right) \right) \cap F_2^{a,c}, f\left( \sigma' \right) = a, f\left( \tilde{\sigma}' \right) = a  \right)\\
&= \sum_{z_{-1}^{a,b} : F \left( z_{-1}^{a,b} \right) \subseteq F_1^{a,b}} \p \left( \sigma \in I_1 \cap I_2 \cap F_2^{a,c}, f\left( \sigma' \right) = a, f\left( \tilde{\sigma}' \right) = a  \, \middle| \,  \sigma \in F \left( \left( 1, z_{-1}^{a,b} \right) \right) \right) \p \left( \sigma \in  F \left( \left( 1, z_{-1}^{a,b} \right) \right) \right)\\
&\leq \sum_{z_{-1}^{a,b} : F \left( z_{-1}^{a,b} \right) \subseteq F_1^{a,b}} \p \left( \sigma :  f\left( \tilde{\sigma}' \right) = a  \, \middle| \,  \sigma \in F \left( \left( 1, z_{-1}^{a,b} \right) \right) \right) \p \left( \sigma \in  F \left( \left( 1, z_{-1}^{a,b} \right) \right) \right).
\end{aligned}
\end{multline*}
Now we know that $\tilde{\sigma}' \in F \left( \left( -1, z_{-1}^{a,b} \right) \right) \subseteq F_1^{a,b}$, and we also know that
\[
\p \left( f \left( \sigma \right) \neq b \, \middle| \, \sigma \in F \left( \left( -1, z_{-1}^{a,b} \right) \right) \right) \leq 4k \gamma.
\]
The number of $\sigma$'s that give the same $\tilde{\sigma}'$ is at most $k^2$, and so we can conclude that
\[
\p \left( \sigma \in I_1 \cap I_2 \cap F_1^{a,b} \cap F_2^{a,c}, f\left( \sigma' \right) = a, f\left( \tilde{\sigma}' \right) = a  \right) \leq 4 k^3 \gamma,
\]
and similarly
\[
\p \left( \sigma \in I_1 \cap I_2 \cap F_1^{a,b} \cap F_2^{a,c}, f\left( \sigma' \right) = a, f\left( \tilde{\sigma}'' \right) = a  \right) \leq 4 k^3 \gamma,
\]
which shows that
\[
\p \left( \sigma \in M_3 \right) \geq \frac{1}{k^6} \left( \frac{\eps^3}{54 n^3 k^{21}} - 2k \gamma - 8 k^3 \gamma \right) \geq \frac{\eps^3}{54 n^3 k^{27}} - \frac{9 \gamma}{k^3}. \qedhere
\]
\end{proof}
To conclude the proof in this case, recall that we have chosen $\gamma = \frac{\eps^3}{10^3 n^3 k^{24}}$.

\subsubsection{Case 2}

First, as in the previous case, we can look at simply the boundary between $a$ and $b$ in direction 1, and conclude that either there are many manipulation points, or there are many local dictators, or \eqref{eq:lg_fbr_final_case_1} holds. This holds similarly for the boundary between $c$ and $d$ in direction 2. Finally, just as in Section \ref{sec:lg_fbr_gen_case2}, we can show that \eqref{eq:lg_fbr_final_case_1} and \eqref{eq:lg_fbr_final_case_2} cannot hold at the same time. We omit the details.

%%%%%%%%%%%%%%%%%%%%%%%%%%%%%%%%%%%%%%%%%%%%%%%%%%%%%%%%%%%%%%%
\subsection{Proof of Theorem \ref{thm:k_refined} concluded} %%%
%%%%%%%%%%%%%%%%%%%%%%%%%%%%%%%%%%%%%%%%%%%%%%%%%%%%%%%%%%%%%%%

\begin{proof}[Proof of Theorem \ref{thm:k_refined}]

Our starting point is Lemma~\ref{lem:boundaries2}, which directly implies Lemma~\ref{lem:cases_ref} (unless there are many 2-manipulation points, in which case we are done). We then consider two cases, as indicated in Section \ref{sec:cases_refined}.

We deal with the small fiber case in Section \ref{sec:sm_fbr_ref}. First, Lemmas \ref{lem:comparable_bdries_ref}, \ref{lem:manip_on_bdry_of_bdry_in_dir!=1_ref}, and \ref{lem:bdry_of_bdry_in_coord1_ref}, and Corollaries \ref{cor:comparable_bdries_ref}, \ref{cor:cases_bdry_dir_ref}, \ref{cor:lg_bdry_in_dir!=1_gives_manip},  and \ref{cor:sm_fbr_ref_last_cor} imply that either there are many 3-manipulation points, or there are many local dictators on three alternatives in coordinate 1. We then deal with the case of many local dictators in Section \ref{sec:loc_dict}. Lemma \ref{lem:loc_dict_abc_to_top}, Corollary \ref{cor:many_dict_at_top}, Lemmas \ref{lem:lg_bdry_for_loc_dict}, \ref{lem:manip_pts_on_bdry_of_loc_dict}, Corollary \ref{cor:manip_by_bdry_loc_dict}, and Lemmas \ref{lem:cond_on_top}, \ref{lem:d_notin_K_ref}, and \ref{lem:final_loc_dict} together show that there are many 4-manipulation points if there are many local dictators on three alternatives, and the SCF is $\eps$-far from the family of nonmanipulable functions.

We deal with the large fiber case in Section \ref{sec:lg_fbr_ref}. Here Lemmas \ref{lem:cond_ab_top_still_lg_fbr}, \ref{lem:towards_inv_hyp_contr}, \ref{lem:applying_inv_hyp_contr}, and \ref{lem:inv_hyp_contr_end} show that if there are not many local dictators on three alternatives, then there are many 3-manipulation points. In the case when there are many local dictators, we refer back to Section \ref{sec:loc_dict} to conclude the proof.
\end{proof}

%%%%%%%%%%%%%%%%%%%%%%%%%%%%%%%%%%%%%%%%%%%%%%%%%%%%%%%%%%%%%%%%%%%%%%%%%%%%%%%%%%%%%%%%%%
\section{Reduction to distance from truly nonmanipulable SCFs}\label{sec:TRUENONMANIP} %%%
%%%%%%%%%%%%%%%%%%%%%%%%%%%%%%%%%%%%%%%%%%%%%%%%%%%%%%%%%%%%%%%%%%%%%%%%%%%%%%%%%%%%%%%%%%

In this section we prove Theorem \ref{thm:TRUENONMANIP}, which says that if our SCF $f$ is close to $\overline{\NONMANIP}$, then it is also close to $\NONMANIP$, or else we have lots of manipulation points. Consequently, this means that if we can prove a quantitative Gibbard-Satterthwaite theorem with distance measured from $\overline{\NONMANIP}$, then we can also prove a quantitative Gibbard-Satterthwaite theorem with distance measured from $\NONMANIP$. In particular, Theorem \ref{cor:k_refined_truenonmanip} is what we get when we combine Theorems \ref{thm:k_refined} and \ref{thm:TRUENONMANIP}.

\begin{proof}[Proof of Theorem \ref{thm:TRUENONMANIP}]
Our assumption means that there exists a SCF $g \in \overline{\NONMANIP}$ such that $\Dist \left( f, g \right) \leq \alpha$. We distinguish two cases: either $g$ is a function of one coordinate, or $g$ takes on at most two values.

\textbf{Case 1. $g$ is a function of one coordinate.} In this case we can assume w.l.o.g.\ that $g$ is a function of the first coordinate, i.e.\ there exists a SCF $h : S_k \to \left[k\right]$ on one coordinate such that for every ranking profile $\sigma$, we have $g \left( \sigma \right) = h \left( \sigma_1 \right)$.

We know from the quantitative Gibbard-Satterthwaite theorem for one voter that for any $\beta$ either $\Dist \left( h, \NONMANIP \left( 1, k \right) \right) \leq \beta$, or $\p \left( \sigma \in M_3 \left( h \right) \right) \geq \frac{\beta^3}{10^5 k^{16}}$.

In the former case, we have that 
\[
\Dist \left( g, \NONMANIP \left( n, k \right) \right) \leq \Dist \left( h, \NONMANIP \left(1, k \right) \right) \leq \beta,
\]
and so consequently
\[
\Dist \left( f, \NONMANIP \left(n, k \right) \right) \leq \alpha + \beta.
\]

In the latter case, we have that
\[
\p \left( \sigma \in M_3 \left( g \right) \right) = \p \left( \sigma \in M_3 \left( h \right) \right) \geq \frac{\beta^3}{10^5 k^{16}},
\]
and so consequently
\[
\p \left( \sigma \in M_3 \left( f \right) \right) \geq \frac{\beta^3}{10^5 k^{16}} - 6 n k \alpha,
\]
since changing the outcome of a SCF at one ranking profile can change the number of 3-manipulation points by at most $6nk$. Now choosing $\beta = 100 n k^6 \alpha^{1/3}$ shows that either \eqref{eq:true_NONMANIP} or \eqref{eq:many_manip} holds.

\textbf{Case 2. $g$ is a function which takes on at most two values.} W.l.o.g.\ we may assume that the range of $g$ is $\left\{a,b\right\} \subset \left[k\right]$, i.e.\ for every ranking profile $\sigma \in S_k^n$ we have $g\left( \sigma \right) \in \left\{a,b\right\}$.

There is one thing we have to be careful about: even though $g$ takes on at most two values, it is not necessarily a Boolean function, since the value of $g \left( \sigma \right)$ does not necessarily depend only on the Boolean vector $x^{a,b} \left( \sigma \right)$.

We now define a function $h : S_k^n \to \left\{a,b \right\}$ that is close in some sense to $g$ and which can be viewed as a Boolean function $h : \left\{a,b \right\}^n \to \left\{a,b \right\}$ because $h\left( \sigma \right)$ depends on $\sigma$ only through $x^{a,b} \left( \sigma \right)$. (The vector $x^{a,b} \left( \sigma \right) \in \left\{ -1, 1 \right\}^{n}$ encodes which of $a$ and $b$ is preferred in each coordinate, and a vector in $\left\{a,b\right\}^{n}$ can encode the same information.) For a given ranking profile $\sigma$, let us consider the fiber on which it is on, $F \left( x^{a,b} \left( \sigma \right) \right)$, and let us define $g|_{F \left( x^{a,b} \left( \sigma \right) \right)}$ to be the restriction of $g$ to ranking profiles in the fiber $F \left( x^{a,b} \left( \sigma \right) \right)$. Then define (see Definition \ref{def:maj})
\[
h \left( \sigma \right) := \Maj \left( g|_{F\left( x^{a,b} \left( \sigma \right) \right)} \right).
\]
By definition, $h\left( \sigma \right)$ depends on $\sigma$ only through $x^{a,b} \left( \sigma \right)$, so we may also view $h$ as a Boolean function $h : \left\{a,b \right\}^n \to \left\{a,b \right\}$.

For any given $0 < \delta < 1$, we either have $\Dist \left( g, h \right) \leq \delta$, in which case $\Dist \left( f, h \right) \leq \alpha + \delta$, or if $\Dist \left( g, h \right) > \delta$, then we show presently that
\begin{equation}\label{eq:many_manip_h}
\p \left( \sigma \in M_2 \left( f \right) \right) \geq \frac{\delta}{4 n k^5} - n k \alpha.
\end{equation}
Choosing $\delta = 8 n^2 k^6 \alpha$ then shows that either \eqref{eq:many_manip} holds, or $\Dist \left( f, h \right) \leq 9 n^2 k^6 \alpha$.

Let us now show \eqref{eq:many_manip_h}. We use a canonical path argument again, but first we divide the ranking profiles according to the fibers with respect to preference between $a$ and $b$.

Let us consider an arbitrary fiber $F\left( z^{a,b} \right)$, and divide it into two disjoint sets: into those ranking profiles for which the outcome of $g$ and $h$ agree, and those for which these outcomes are different. I.e.\
\[
F\left( z^{a,b} \right) = F^{\maj} \left( z^{a,b} \right) \cup F^{\min} \left( z^{a,b} \right),
\]
where
\begin{align*}
F^{\maj} \left( z^{a,b} \right) &= \left\{ \sigma \in F \left( z^{a,b} \right) : g \left( \sigma \right) = h \left( \sigma \right) \right\},\\
F^{\min} \left( z^{a,b} \right) &= \left\{ \sigma \in F \left( z^{a,b} \right) : g \left( \sigma \right) \neq h \left( \sigma \right) \right\}.
\end{align*}
By construction, we know that
\[
\left| F^{\min} \left( z^{a,b} \right) \right| \leq \frac{1}{2} \left| F \left( z^{a,b} \right) \right| = \frac{1}{2} \left( \frac{k!}{2} \right)^n.
\]
Now for every pair of profiles $\left( \sigma, \sigma' \right) \in F^{\min} \left( z^{a,b} \right) \times F^{\maj} \left( z^{a,b} \right)$ define a canonical path from $\sigma$ to $\sigma'$ by applying a path construction in each coordinate one by one, and then concatenating these paths. In each coordinate we apply the path construction of \cite[Proposition 6.6.]{isaksson2010geometry}: we bubble up everything except $a$ and $b$, and then finally bubble up the last two alternatives as well.

For a given edge $\left( \pi, \pi' \right) \in F^{\min} \left( z^{a,b} \right) \times F^{\maj} \left( z^{a,b} \right)$ there are at most $2 k^4 \left( \frac{k!}{2} \right)^n$ possible pairs $\left( \sigma, \sigma' \right) \in F^{\min} \left( z^{a,b} \right) \times F^{\maj} \left( z^{a,b} \right)$ such that the canonical path between $\sigma$ and $\sigma'$ defined above passes through $\left( \pi, \pi' \right)$. (This can be shown just like in the previous lemmas, e.g.\ Lemma \ref{lem:comparable_bdries_ref}.) Consequently we have
\[
\left| \partial_e \left( F^{\min} \left( z^{a,b} \right) \right) \right| \geq \frac{\left| F^{\min} \left( z^{a,b} \right) \right| \left| F^{\maj} \left( z^{a,b} \right) \right|}{2 k^4 \left( \frac{k!}{2} \right)^n} \geq \frac{\left| F^{\min} \left( z^{a,b} \right) \right|}{4 k^4},
\]
where the edge boundary $\partial_e \left( F^{\min} \left( z^{a,b} \right) \right)$ is defined via the refined rankings graph restricted to the fiber $F \left( z^{a,b} \right)$. Summing this over all fibers we have that
\begin{equation}\label{eq:sum_fbrs}
\sum_{z^{a,b}} \left| \partial_e \left( F^{\min} \left( z^{a,b} \right) \right) \right| \geq \sum_{z^{a,b}} \frac{\left| F^{\min} \left( z^{a,b} \right) \right|}{4 k^4} \geq \frac{\delta}{4 k^4} \left( k! \right)^n,
\end{equation}
using the fact that $\Dist \left( g, h \right) > \delta$.

Now it is easy to see that if $\left( \sigma, \sigma' \right) \in \partial_e \left( F^{\min} \left( z^{a,b} \right) \right)$ for some $z^{a,b}$, then either $\sigma$ or $\sigma'$ is a 2-manipulation point for $g$. In the refined rankings graph every vertex (ranking profile) has $n \left( k - 1 \right) < nk$ neighbors, so each 2-manipulation point can be counted at most $nk$ times in the sum on the left hand side of \eqref{eq:sum_fbrs}, showing that
\[
\p \left( \sigma \in M_2 \left( g \right) \right) \geq \frac{\delta}{4nk^5},
\]
from which \eqref{eq:many_manip_h} follows immediately, since changing the outcome of a SCF at one ranking profile can change the number of 2-manipulation points by at most $nk$.

So either we are done because \eqref{eq:many_manip} holds, or $\Dist \left( f, h \right) \leq 9 n^2 k^6 \alpha$; suppose the latter case. Our final step is to look at $h$ as a Boolean function, and use a result on testing monotonicity \cite{goldreich2000testing}.

Denote by $\tilde{\Dist}$ the distance of $h$ when viewed as a Boolean function from the set of monotone Boolean functions. Let $0 < \eps < 1$ be arbitrary. Then either $\tilde{\Dist} \leq \eps$, in which case $\Dist \left(h, \NONMANIP \right) \leq \tilde{\Dist} \leq \eps$ and therefore $\Dist \left( f, \NONMANIP \right) \leq 9 n^2 k^6 \alpha + \eps$, or $\tilde{\Dist} > \eps$. In the latter case we show that then
\begin{equation}\label{eq:manip_eps}
\p \left( \sigma \in M_2 \left( f\right) \right) \geq \frac{2\eps}{nk} - 9 n^3 k^7 \alpha.
\end{equation}
Choosing $\eps = 5 n^4 k^8 \alpha$ then shows that either \eqref{eq:true_NONMANIP} or \eqref{eq:many_manip} holds.

Let us now show \eqref{eq:manip_eps}. Let us view $h$ as a Boolean function, and denote by $p\left( h \right)$ the fraction of pairs of strings, differing on one coordinate, that violate the monotonicity condition. Goldreich, Goldwasser, Lehman, Ron, and Samorodnitsky showed in \cite[Theorem 2]{goldreich2000testing} that $p \left( h \right) \geq \frac{\tilde{\Dist}}{n}$.

Now going back to viewing $h$ as a SCF on $k$ alternatives, this tells us that there are at least $\frac{\eps}{2} 2^n$ pairs of fibers, which differ on one coordinate, that violate monotonicity. For each such pair of fibers, whenever $a$ and $b$ are adjacent in the coordinate where the two fibers differ, we get a 2-manipulation point. Such a 2-manipulation point can be counted at most $n$ times in this way (since there are $n$ coordinates where $a$ and $b$ can be adjacent). Consequently, we have
\[
\left| M_2 \left( h \right) \right| \geq \frac{\eps}{2} \cdot 2^n \cdot 2 \left( k - 1 \right)! \left( \frac{k!}{2} \right)^{n-1} \cdot \frac{1}{n} = \frac{2\eps}{nk}\left( k! \right)^n,
\]
i.e.\
\[
\p \left( \sigma \in M_2 \left( h \right) \right) \geq \frac{2 \eps}{nk},
\]
from which \eqref{eq:manip_eps} follows immediately, since changing the outcome of a SCF at one ranking profile can change the number of 2-manipulation points by at most $nk$.
\end{proof}

\begin{proof}[Proof of Theorem \ref{cor:k_refined_truenonmanip}]
First we argue without specific bounds. Suppose on the contrary that our SCF $f$ does not have many 4-manipulation points. Then $f$ is close to $\overline{\NONMANIP}$ by Theorem \ref{thm:k_refined}. Consequently, by Theorem \ref{thm:TRUENONMANIP}, $f$ is close to $\NONMANIP$, which is a contradiction.

Now we argue with specific bounds. Assume on the contrary that
\[
\p \left( \sigma \in M_4 \left( f \right) \right) < \frac{\eps^{15}}{10^{39} n^{67} k^{166}}.
\]
Then by Theorem \ref{thm:k_refined} we have that $\Dist \left( f, \overline{\NONMANIP} \right) < \frac{\eps^3}{10^6 n^{12} k^{24}}$, and consequently by Theorem \ref{thm:TRUENONMANIP} we have $\Dist \left( f, \NONMANIP \right) < \eps$, which is a contradiction.
\end{proof}

%%%%%%%%%%%%%%%%%%%%%%%%%%%%%%%%%%%%%%%%%%%
\section{Open problems}\label{sec:open} %%%
%%%%%%%%%%%%%%%%%%%%%%%%%%%%%%%%%%%%%%%%%%%

We conclude with a few open problems that arise naturally, some of which have already been asked by Isaksson, Kindler and Mossel~\cite{isaksson2010geometry}.
\begin{itemize}
 \item In Section~\ref{sec:discussion} we mentioned that our techniques do not lead to tight bounds. It would be interesting to find the correct tight bounds. When discussing tight bounds there are various different ways to measure the manipulability of a function: in terms of the probability of having manipulating voters, in terms of the expected number of manipulating voters, in terms of the number of manipulative edges (either in the refined or non-refined graph), etc.

 \item A related question is to find, in some natural subsets of functions, the one that minimizes manipulation. For example, among anonymous SCFs, which function minimizes the expected number of manipulating voters? For example, for plurality, the probability that a ranking profile is manipulable is  $\Theta \left( 1 / \sqrt{n} \right)$, and if it is manipulable, then $\Theta \left( n \right)$ voters can manipulate individually, so consequently the expected value of the number of voters who can manipulate individually is $\Theta \left( \sqrt{n} \right)$. Is it true that for all anonymous SCFs, this expectation is $\Omega \left( \sqrt{n} \right)$?
% 
%  \item Among specific classes of SCFs, say anonymous SCFs, which function is the best? I.e., which function minimizes the probability of a ranking profile being manipulable?
% 
%  \item It is interesting to look at more detailed information about manipulation properties. For instance, given a manipulable ranking profile, how many voters can manipulate individually? What is the expected value of the number of voters who can manipulate individually? For plurality, the probability that a ranking profile is manipulable is $\Theta \left( 1 / \sqrt{n} \right)$, and if it is manipulable, then $\Theta \left( n \right)$ voters can manipulate, so consequently this expected value is $\Theta \left( \sqrt{n} \right)$. Is it true that for all anonymous SCFs, this expectation is $\Omega \left( \sqrt{n} \right)$?

 \item What if the distribution over rankings is not i.i.d.\ uniform? It would be interesting to consider a quantitative Gibbard-Satterthwaite theorem, and also the questions asked above, in this setting.
\end{itemize}

%%%%%%%%%%%%%%%%%%%%%%%%
%%% Acknowledgements %%%
%%%%%%%%%%%%%%%%%%%%%%%%

\section*{Acknowledgments}

We thank Ariel Procaccia, Gireeja Ranade and Piyush Srivastava for helpful discussions, and Joe Neeman for valuable comments on a draft of the paper. We also thank anonymous referees for useful comments.

%%%%%%%%%%%%%%%%%%%%
%%% Bibliography %%%
%%%%%%%%%%%%%%%%%%%%

\bibliographystyle{plain}
\bibliography{manip}

%%%%%%%%%%%%%%%%%%
\end{document}